\documentclass[leqno,12pt]{amsart}
\usepackage{amssymb}

\newtheorem{theorem}{\indent\sc Theorem}[section]
\newtheorem{lemma}[theorem]{\indent\sc Lemma}
\newtheorem{corollary}[theorem]{\indent\sc Corollary}
\newtheorem{proposition}[theorem]{\indent\sc Proposition}

\theoremstyle{definition}
\newtheorem{definition}[theorem]{\indent\sc Definition}
\newtheorem{remark}[theorem]{\indent\sc Remark}
\newtheorem{example}[theorem]{\indent\sc Example}

\begin{document}

\title[Mollifier Smoothing of tensor fields on differentiable
manifolds]{Mollifier Smoothing of tensor fields on differentiable
manifolds and applications to Riemannian Geometry}

\author[Ryuichi Fukuoka]{Ryuichi Fukuoka}

\dedicatory{Dedicated to Professor Francesco Mercuri on his
sixtieth birthday}

\subjclass[2000]{Primary 53B21; Secondary 41A35, 53A45, 53B20.}
%
\keywords{Mollifier smoothing, Non-regular Riemannian metrics,
Levi-Civita connection, Riemannian curvature tensor,
Lipschitz-Killing curvature measures.}

\address{
Universidade Estadual de Maring\'a \endgraf Departamento de
Matem\'atica
\endgraf Av. Colombo 5790 \endgraf Maring\'a - PR 87020-900 \endgraf Brazil }
\email{rfukuoka@uem.br, ryuichifukuoka@yahoo.com.br}


\begin{abstract}
Let $M$ be a differentiable manifold. We say that a tensor field
$g$ defined on $M$ is non-regular if $g$ is in some local $L^p$
space or if $g$ is continuous. In this work we define a mollifier
smoothing $g_\varepsilon$ of $g$ that has the following feature:
If $g$ is a Riemannian metric of class $C^2$, then the Levi-Civita
connection and the Riemannian curvature tensor of $g_\varepsilon$
converges to the Levi-Civita connection and to the Riemannian
curvature tensor of $g$ respectively as $\varepsilon$ converges to
zero. Therefore this mollifier smoothing is a good starting point
in order to generalize objects of the classical Riemannian
geometry to non-regular Riemannian manifolds. Finally we give some
applications of this mollifier smoothing. In particular, we
generalize the concept of Lipschitz-Killing curvature measure for
some non-regular Riemannian manifolds.
\end{abstract}

\maketitle

\section{Introduction}
\label{introducao}

The concept of Riemannian manifold can be generalized in several
ways. Finsler manifolds, sets with positive reach in Euclidean
spaces and Alexandrov spaces are examples of generalizations of
some classes of Riemannian manifolds. A good starting point for
the reader who is interested in this subject is M. Berger's book
\cite{2} and references therein. In this work we introduce the
$L^p_{\mathrm {loc}}$ tensor fields and the continuous tensor
fields on differentiable manifolds, which we call {\em non-regular
tensor fields}. We are particularly interested to study
non-regular Riemannian metrics $\widehat g$ on $n$-dimensional
differentiable manifolds $M^n$. We call a pair $(M^n,\widehat g)$
by {\em non-regular Riemannian manifold}.

As an example of an $L^p_{\mathrm{loc}}$ Riemannian manifold (that
is, a differentiable manifold with an $L^p_{\mathrm{loc}}$
Riemannian metric), fix a convex set $D \subset \mathbb R^n$ and
suppose that it is closed and bounded. Its boundary $\partial D$
is called a convex surface. It is well known that $\partial D$
have a tangent space almost everywhere and we can define an inner
product from the ambient space on each tangent space. It is not
difficult to see that $\partial D$ endowed with this field of
inner products is a non-regular Riemannian manifold. Eqs.
(\ref{exemplometrica1}), (\ref{exemplometrica2}) and
(\ref{exemplometrica3}) give other examples of
$L^p_{\mathrm{loc}}$ Riemannian manifolds.

Non-regular Riemannian manifolds can not be studied through
differentiation as in classical Riemannian geometry. However the
space of non-regular tensor fields can be topologized in the same
way as the space of non-regular functions. The main idea of this
work is to study non-regular Riemannian manifolds through
approximations by (smooth) Riemannian manifolds. Let $(M,\widehat
g)$ be a non-regular Riemannian manifold and suppose that
$\{(M,\widehat g(\varepsilon)), \varepsilon > 0 \}$ is a
one-parameter family of Riemannian manifolds such that $\widehat
g(\varepsilon)$ converges to $\widehat g$ in
$L^p_{\mathrm{loc}}(M)$ (or in $C^0_{\mathrm{loc}}(M)$) as
$\varepsilon$ goes to zero. We can try to define some geometric
object in $(M,\widehat g)$ studying the behavior of the
correspondent object in the family $\{(M,\widehat
g(\varepsilon)),\varepsilon > 0\}$ as $\varepsilon$ goes to zero.
As an illustration of this approach we can imagine the round
sphere of radius $1$ in $\mathbb R^3$ deforming gradually and
converging to a cube. Observe that the total curvature during the
process is equal to $4\pi$. Moreover the curvature concentrates
around the vertices of the cube while the deformation follows,
because the other points are becoming flat. If the deformation is
made ``symmetrically'', then the total curvature in a neighborhood
of each vertex is equal to $\pi/2$. Therefore the Gaussian
curvature of a vertex $p$ can be thought as $(\pi/2)$. In this
way, we ``were able to generalize'' the concept of Gaussian
curvature for the cube.

The distance between two points in a non-regular Riemannian
manifold $(M,\widehat g)$ is also defined using approximation by
(smooth) Riemannian manifolds (See Definition \ref{distancialp}
and Remark \ref{distanciac0}). If we identify the points with zero
distance, then the resulting identification space is a metric
space that is not necessarily homeomorphic to a differentiable
manifold. We can see Eqs. (\ref{exemplometrica1}) and
(\ref{exemplometrica3}) again to have an idea of which kind of
metric spaces a non-regular Riemannian manifold can generate.

Let us explain how this work is organized. Meanwhile, we cite some
related works. In Section \ref{preliminares}, we fix some
notations and present some classical results that will be used
afterwards. In Section \ref{espacoslpec0}, we establish the
foundations of non-regular tensor fields on differentiable
manifolds. In Section \ref{mollifiergeral}, we introduce the
mollifier smoothing of a non-regular tensor field $\widehat T$ in
a Riemannian manifold $(M,\widetilde g)$. A smooth Riemannian
metric $\widetilde g$, which is called the {\em background
metric}, is introduced as an ``auxiliary'' metric which is
necessary to define the mollifier smoothing. We are particularly
interested when $\widehat T$ is a non-regular Riemannian metric.

In Section \ref{suavizacaosubordinadaparticao} we adapt the
mollifier smoothing introduced in Section \ref{mollifiergeral} in
the following way: Take a locally finite covering of $M$ by open
sets with Euclidean background metric. Besides, take a partition
of the unity subordinated to this covering. We call this covering
with this partition of the unity by $\mathcal P$. Consider the
mollifier smoothing defined in Section \ref{mollifiergeral} on
each open set of the covering. The sum of these mollifier
smoothings weighted by the partition of the unity is called {\em
mollifier smoothing with respect to $\mathcal P$}. It has the
following interesting property: If the non-regular Riemannian
metric $\widehat g$ is of class $C^2$, then the Levi-Civita
connection and the Riemannian curvature tensor of the mollifier
smoothing $\widehat g_\varepsilon$ with respect to $\mathcal P$
converges to the Levi-Civita connection and the Riemannian
curvature tensor of $\widehat g$ respectively as $\varepsilon$
goes to zero. Therefore this mollifier smoothing provides a
natural starting point in order to generalize objects of the
classical Riemannian Geometry to non-regular Riemannian manifolds.

It is interesting to present other mollifier smoothings that we
can find in the literature: The mollifier smoothing of
distributions (See the regularization ``$R$'' in \cite{8}) is
quite similar to the mollifier smoothing with respect to $\mathcal
P$. The main difference is that the ``gluing'' of local smoothings
is done there by composition of the local smoothings. Another work
where we find a mollifier smoothing is in Nash's celebrated work
about embeddings of Riemannian manifolds into Euclidean spaces
(see \cite{14}). There he defines a smoothing of tensor fields on
Riemannian manifolds. He embeds the manifold into an Euclidean
space and he makes the convolution of the tensor field with a
mollifier. This mollifier is defined on the ambient space and it
decays quickly as it goes to infinity. This convolution eliminates
the ``high frequencies'' of the original tensor field, what is
typical in Fourier analysis. Our definition of mollifier smoothing
is intrinsic and our kernel has compact support. Finally another
work that is worth mentioning (although it is not directly related
with our work) is Karcher's paper \cite{11}. There the author
introduces a mollifier smoothing of a map between Riemannian
manifolds using the center of mass.

The last sections are devoted to some applications. We do not want
to be extensive there. Our aim is to convince the reader that the
mollifier smoothings of tensor fields given in Definitions
\ref{definicaosuavizacao} and \ref{suavizacaoemrelacaoaP} can be
useful in order to study non-regular Riemannian metrics and its
singularities. The author hopes that this theory can provide
useful examples to the people who study the convergence of
sequences of Riemannian metrics, although the author's lack of
knowledge in this field does not allow him to say anything deeper
about this subject.

In Section \ref{distancianaoregular} we generalize the concept of
distance between two points for non-regular Riemannian manifolds.
There are other works that study the concept of distance for
non-regular manifolds. For instance in \cite{3} and \cite{4},
Cecco and Palmieri generalize the concept of distance on Lipschitz
manifolds with a Riemannian metric. There the non-regularity is
imposed by a Lipschitz atlas, which is an atlas with Lipschitz
change of coordinates. Particularly interesting is the definition
of distance between two points given in \cite{4}, which depends on
an integral over $(M,\widehat g)$.

In Section \ref{transporteparalelonaoregular} we generalize the
parallel transport for some non-regular Riemannian manifolds using
mollifier smoothing with respect to $\mathcal P$. We prove that it
is well defined for some curves on piecewise smooth
two-dimensional Riemannian manifolds (see Definition
\ref{definicaosuperficiesuaveporpartes}). For the reader who is
interested in other works about parallel transport in less regular
spaces, Nikolaev and Petrunin study the parallel transport in
Alexandrov spaces in \cite{15} and \cite{16}.

In Section \ref{curvaturalknaoregular} we define the
Lipschitz-Killing curvature measures for closed non-regular
Riemannian manifolds $(M,\widehat g)$ (see Definition
\ref{definicaomedidadecurvatura}). It is a measure that represents
the total Lipschitz-Killing curvature on Borel sets in $M$. The
existence of the Lipschitz-Killing curvature measures depends on
the existence of a Lipschitz-Killing curvature measure generator,
which is essentially a family of ``regular'' open sets of $M$ that
generates the topology of $M$ (See Definition
\ref{geradorcurvaturemeasure} for more details). The
Lipschitz-Killing curvature measure for these generators is
defined using the mollifier smoothing with respect to $\mathcal
P$. Then we extend the Lipschitz-Killing curvature measure to
Borel subsets of $M$.

In \cite{5}, Cheeger, M\"uller and Schrader study sequences of
piecewise flat manifolds that converges to a smooth Riemannian
manifold $(M,g)$. Although this seems to be the opposite we are
doing here, the general idea is the same: to approximate a
Riemannian manifold by simpler objects. There the authors prove
that the Lipschitz-Killing curvature of some sequences of
piecewise flat manifolds converges, in the sense of measures, to
the Lipschitz-Killing curvatures of the smooth Riemannian
manifold. A similar convergence holds for the mean curvatures of
the boundary.

The concept of curvature measures is also used, for instance, in
the Federer's theory of subsets with positive reach in $\mathbb
R^n$ (See \cite{10}).  Let $A \subset \mathbb R^n$ be a subset
with positive reach and let $B \subset A$ be a Borel subset of
$\mathbb R^n$. The coefficients of the Steiner polynomial of $B
\subset A$ are essentially the curvature measures of $B \subset A$
and some of these coefficients are essentially generalizations of
the Lipschitz-Killing curvature measures for the case of smooth
hypersurfaces.

Another instance where the concept of curvature measure appear is
in Alexandrov theory of convex surfaces (See \cite{1}). There the
Gaussian curvature measure is defined in the following way: in
open geodesic triangles it is defined as $2\pi$ minus the sum of
the internal angles; in open geodesic arcs it is defined as zero;
at the points it is defined as $2\pi$ minus the total angle at the
point. From these curvature measures, we are able to extend the
total Gaussian curvature to Borel subsets of the surface.

In Section \ref{gaussbonnetnaoregular} we study piecewise smooth
two-dimensional Riemannian manifolds more deeply. They are
essentially polyhedras with non-flat faces and edges. We get a
Gaussian curvature measure generator for these kind of surfaces
and we prove that the Gaussian curvature measure has the
``expected'' geometrical values. Meanwhile we get an alternative
proof of a (well known) generalization of the Gauss-Bonnet theorem
for these kind of surfaces (See Theorem
\ref{curvaturasuavizacaoconverge}).

It is completely natural to conjecture that the Lipschtz-Killing
curvature measures we define here coincides with the classical
definitions if the geometrical object lies in the intersection of
both theories.

In Section \ref{secaocurvaturadimensaofracionaria}, we present an
odd example of a two-dimensional sphere with a Riemannian metric
of class $C^1$ such that it is flat outside a subset with
Hausdorff dimension $(1+(\ln 2/\ln 3))$. In other words, its
curvature is ``concentrated'' in a subset with Hausdorff dimension
$(1+(\ln 2/\ln 3))$, in the same way that the curvature of a cube
is concentrated on its vertices (zero dimensional subsets). This
sphere shows a strange behaviour that a curvature can assume in
non-regular Riemannian manifolds. It can also be useful in order
to study the ``dimensional character'' of the curvature at a
point, like the zero-dimensional curvature at the vertices, the
one dimensional curvature at the edges and the two-dimensional
curvature at the faces of a piecewise smooth two-dimensional
Riemmannian manifold.

\section{Preliminaries} \label{preliminares}

In this section we fix some notations and remember some classical
theorems that will be used afterwards. The material exposed here
can be found in any good textbook of Riemannian Geometry and
Measure Theory, for instance \cite{6}, \cite{7}, \cite{9},
\cite{12}, \cite{17}, \cite{18} and \cite{20}. The exception is
Theorem \ref{transporteparalelosuave}, which proof is given here.
In Riemannian geometry, there exist two definitions of curvature
tensor which differ by a sign. We adopt the convention given in
\cite{12}.

\subsection{Riemannian geometry}
\label{preliminargeometriariemanniana}

Let $M^n$ be an $n$-dimensional differentiable manifold (The
superscript $n$ is omitted whenever there is not any possibility
of misunderstandings). Denote the tangent bundle over $M$ by $TM$,
the cotangent bundle over $M$ by $T^*M$ and the tensor bundle of
type $(m,s)$ over $M$ by $T^{m,s}M$. The respective fibers over
$x\in M$ is denoted by $T_xM$, $T^*_xM$ and $T^{m,s}_xM$. If
$\varphi_1,\ldots,\varphi_m\in TM$, $v_1,\ldots, v_s\in T^*M$ and
$T\in T^{m,s}M$, then
$T(\varphi_1,\ldots,\varphi_m,v_1,\ldots,v_s)$ is the contraction
of these tensor fields.

Let $(M,g)$ be an $n$-dimensional Riemannian manifold with smooth
metric $g$. As usual we use the notation
$\left<u,v\right>_g:=g(u,v)$ for the scalar product of the vectors
$u$ and $v$ and the subscript $g$ will be omitted whenever there
is not any possibility of misunderstandings. The symbol $\nabla$
denote the Levi-Civita connection and $R$ denote the curvature
tensor, where $R(u,v) = \nabla_u \nabla_v - \nabla_v \nabla_u -
\nabla_{[u,v]}$. The Ricci tensor is defined by
$Ric(u,v)=\sum\limits^n_{i=1}\left<R(w_i,u)v,w_i\right>$ and the
scalar curvature by
$S=\sum\limits_{i,j=1}^n\left<R(w_i,w_j)w_j,w_i\right>$, where
$\{w_1,\ldots,w_n\}$ is an orthonormal basis of $T_xM$.

Put a coordinate system $(x_1,\ldots,x_n)$ in a neighborhood of
$x\in M$. The coordinate vector fields is denoted by
$\partial/\partial x_i,i=1,\ldots,n$. The components of the metric
with respect to this coordinate system is denoted by $g_{ij}=$
$g(\partial/\partial x_i,\partial/\partial x_j),$ $i,j=1,\ldots,
n$. The Christoffel symbols are denoted by $\Gamma^k_{ij}$ and
they are defined implicitly as
\[
\nabla_{\frac{\partial}{\partial x_i}}\frac{\partial}{\partial
x_j}=\sum_{k=1}^n \Gamma^k_{ij}\frac{\partial}{\partial x_k}.
\]
It is well known that
\begin{equation}
\label{levicivitacoordenadas}
\Gamma^k_{ij}=\frac{1}{2}\sum_{k=1}^n \left\{
\frac{\partial}{\partial x_i} g_{jm}+\frac{\partial}{\partial x_j}
g_{mi}-\frac{\partial}{\partial x_m} g_{ij}\right\}g^{km}
\end{equation}
where $g^{km}$ are the components of the inverse matrix of $g$.

The components of the curvature tensor are given implicitly by
\begin{equation}
\label{tensordecurvatura} R\left(\frac{\partial}{\partial
x_i},\frac{\partial}{\partial x_j}\right)\frac{\partial}{\partial
x_k}=\sum_{l=1}^n R^l_{kij}\frac{\partial}{\partial x_l}.
\end{equation}
It is well known that
\begin{equation}
\label{curvaturacoordenadas} R^l_{kij}=\left(\frac{\partial
\Gamma^l_{jk}}{\partial x_i}-\frac{\partial
\Gamma^l_{ik}}{\partial x_j}\right)+\sum_{m=1}^n
\left(\Gamma^m_{jk}\Gamma^l_{im}-\Gamma^m_{ik}\Gamma^l_{jm}\right).
\end{equation}

Let $\{w_1,\ldots,w_n\}$ be an orthonormal moving frame for the
tangent bundle of $(M,g)$ and denote the dual frame by
$\{\varpi_1,\ldots, \varpi_n\}$. Define the curvature forms
$\Omega_{lk}$ for the orthonormal moving frame
$\{w_1,\ldots,w_n\}$ by
\begin{equation}
\label{formasdecurvatura} \sum_{l=1}^{n}\Omega_{lk}(w_i,w_j)w_l =
R(w_i,w_j)w_k.
\end{equation}

For $\kappa = 1,\ldots,n$, the $\kappa$-th Lipschitz-Killing
curvature measure $\mathcal R^\kappa$ is a measure defined on open
subsets $U \subset M$: if $\kappa$ is odd, then $\mathcal
R^\kappa$ is identically zero; if $\kappa$ is even then
\begin{equation}
\label{curvaturaintegralipschitzkilling} \mathcal
R^\kappa(U)=\int_U R^\kappa
\end{equation}
where
\begin{equation}
\label{lipschitzkillingform} R^\kappa =
\frac{(-1)^{\frac{\kappa}{2}}}{(n-\kappa)! 2^\kappa
\pi^{\frac{\kappa}{2}}(\kappa /2)!} \sum_\sigma (-1)^{ | \sigma |
}\Omega_{\sigma(1)\sigma(2)}\wedge \ldots \wedge
\Omega_{\sigma(\kappa - 1)\sigma(\kappa)}\wedge
\varpi_{\sigma(\kappa + 1)} \wedge \ldots \wedge
\varpi_{\sigma(n)},
\end{equation}
and the summation is over all permutations of $n$ elements. Notice
that $R^\kappa$ does not depend on the orthonormal moving frame we
choose. The curvature measure $\mathcal R^0$ is the volume and
$\mathcal R^2$ is proportional to the total scalar curvature.

The volume element of $(M,g)$ is denoted by $dV_g(\cdot)$. When
$(M,g)$ is the Euclidean space with its canonical metric, the
volume element is simply denoted by $dV(\cdot)$. $B_g(x,r):=\{y\in
M;\mathrm{dist}_g(x,y)<r\}$ denote the geodesic ball with center
$x$ and radius $r$.

The exponential map is denoted by $exp:\Psi M\subset TM\rightarrow
M$. The exponential map restricted to the tangent space $T_xM$ is
denoted by $exp_x:\Psi_x M \subset T_xM\rightarrow M$.

Suppose that there exist a unique minimizing geodesic $\gamma$
connecting $x$ and $y$ in $(M,g)$. We denote the parallel
transport between the tensor spaces $T^{m,s}_xM$ and $T^{m,s}_yM$
through $\gamma$ by ${\tau}_{x,y}$. The following theorem states
essentially that the parallel transport inside the injectivity
radius is smooth with respect to all its parameters. Although
natural, it is not usually found as stated here.

\begin{theorem}
\label{transporteparalelosuave} Let $(M^n,g)$ be an n-dimensional
Riemannian manifold with smooth Riemannian metric $g$. Let $q_1$
and $q_2$ be points in $M$ and $\gamma$ be the unique minimizing
geodesic joining them such that
$\exp_{q_1}(\gamma^\prime(0))=q_2$. Then there exists a
neighborhood $N_{q_1}$ of $q_1\in M$ and $N_{q_2}$ of $q_2\in M$
such that the map $\tau:TN_{q_1}\times N_{q_2}\rightarrow M\times
TM$ defined by $\tau((x,v),y)=(x,(y,\tau_{x,y}(v)))$ is smooth.
\end{theorem}

\begin{proof}
Set $TMM=\{(x,\dot x,v);x\in M, (\dot x,v)\in T_xM\times T_xM\}$.
Define $\rho_1:U\subset TMM\rightarrow M\times TM$ which is given
by $\rho_1(x,\dot x,v)=(x,(\exp(x,\dot x),\tau_{x,\exp(x,\dot
x)}v))$, where $U$ is chosen such that $\rho_1$ is well defined.

Parametrize $U$ by $(x_1,x_2,\ldots,x_n,\dot x_1,\ldots,\dot
x_n,v_1,\ldots, v_n)$. The function $\rho_1$ is smooth because
$\rho_1(x_0,\dot x_0,v_0)$ is the solution of the system of
ordinary differential equations

\begin{eqnarray*}
& & \frac{d\dot x_k}{dt}+\sum_{i,j}\Gamma_{ij}^k \dot x_i \dot
x_j=0
\hspace{1cm} k=1,\ldots,n \\
& & \frac{dx_d}{dt}=\dot x_d \hspace{3cm} d=1,\ldots,n \\
& & \frac{dv_l}{dt}+\sum_{i,j}\Gamma_{ij}^l \dot x_i v_j = 0
\hspace{1.2cm} l=1,\ldots,n
\end{eqnarray*}
with initial conditions $x(0)=x_0$, $\dot x(0)=\dot x_0$ and
$v(0)=v_0$, which is smooth with respect to the initial
conditions.

Let $q_1$, $q_2$ and $\gamma$ as in the hypothesis of the theorem.
Observe that the Jacobian matrix of $\rho_1$ at
$(q_1,\widetilde\gamma^\prime(0),v)$ is nonsingular for every
$v\in T_{q_1}M$. Therefore there exists a neighborhood $U_1$ of
$(q_1,\widetilde\gamma^\prime(0),v)$ such that
$\rho_1\vert_{U_1}:U_1\rightarrow \rho_1(U_1)$ is a diffeomorfism.
Observe that we can pick $U_1=\{(x,\dot x,v)\}$ such that every
$v\in T_xM$ is included.

Now define $\rho_2:U\subset TMM\rightarrow TM\times M$ which is
given by $\rho_2(x,\dot x,v)=((x,v),exp(x,\dot x))$ where $U$ is
chosen such that $\rho_2$ is well defined. Using the same argument
as before, if we take $q_1$, $q_2$ and $\gamma$ as in the
hypothesis of the theorem, then there exists a neighborhood $U_2$
of $(q_1,\gamma^\prime(0),v)$ such that
$\rho_2\vert_{U_2}:U_2\rightarrow \rho_2(U_2)$ is a diffeomorfism.
Observe that we can also pick $U_2=\{(x,\dot x,v)\}$ such that
every $v\in T_xM$ is included.

Finally observe that there exist neighborhoods $N_{q_1}$ of
$q_1\in M$ and $N_{q_2}$ of $q_2\in M$ such that
$\tau:=\rho_1\circ \rho_2^{-1}$ defined on $TN_{q_1}\times
N_{q_2}$ is smooth.

\end{proof}

\subsection{Measure theory}
\label{preliminarteoriadamedida}

\begin{theorem}
\label{lebesguedifferentiation} Let $f:\mathbb R^n \rightarrow
\mathbb R$ be a locally summable function. Then for a.e.
$x\in\mathbb{R}^n$, we have that
\begin{equation}
\label{formulalebesgue} \frac{\int_{B(x,r)}\mid f(y)-f(x)\mid
dV(y)}{\mathrm{Vol}(B(x,r))}\rightarrow 0
\end{equation}
as $r\rightarrow 0$. Such a point $x$ is called a Lebesgue point
of $f$.
\end{theorem}

\begin{definition}
\label{compactamentecontido} We say that an open subset $U \subset
M$ is compactly contained in $M$ if the closure $\bar U$ of $U$ is
compact. In this case, we denote $U \subset \subset M$.
\end{definition}

\begin{remark}
The Lebesgue's Differentiation Theorem is valid for smooth
Riemannian manifolds $(M^n,g)$. In fact, take a coordinate system
$\phi:U \rightarrow \mathbb R^n$ in $U \subset \subset M$ such
that $\phi(U) \subset \subset \mathbb R^n$. Consider three metrics
in $U$: $g_{ij}$, $c.\delta_{ij}$ and $C.\delta_{ij}$. Here
$\delta_{ij}$ is the Euclidean metric with respect to $\phi$ and
the constants $c$ and $C$ are chosen in such a way that $\Vert
v\Vert_{c\delta}\leq \Vert v\Vert_{g}\leq \Vert v\Vert_{C\delta}$
for every $v\in TU$. Suppose that $x\in U$ is a Lebesgue point of
$f$ in $(M,c\delta)$. Then

\[
\frac{\int_{B_g(x,r)}\mid f(y)-f(x)\mid
.dV_g(y)}{\int_{B_g(x,r)}.dV_g(y)}\leq \frac{\int_{B_g(x,r)}\mid
f(y)-f(x)\mid .dV_{C\delta}(y)}{\int_{B_g(x,r)}.dV_{c\delta}(y)}
\]
\[
\leq \frac{\int_{B_{c\delta}(x,r)}\mid f(y)-f(x)\mid
.dV_{C\delta}(y)}{\int_{B_{C\delta}(x,r)}.dV_{c\delta}(y)}\leq
\frac{C^n}{c^n} \frac{\int_{B_{c\delta}(x,r)}\mid f(y)-f(x)\mid
.dV_{c\delta}(y)}{\int_{B_{C\delta}(x,r)}.dV_{c\delta}(y)}
\]
\[
\leq \frac{C^{2n}}{c^{2n}} \frac{\int_{B_{c\delta}(x,r)}\mid
f(y)-f(x)\mid
.dV_{c\delta}(y)}{\int_{B_{c\delta}(x,r)}.dV_{c\delta}(y)}\rightarrow
0\hspace{5mm}\mathrm{as}\;\;r\rightarrow 0.
\]
Therefore every Lebesgue point of $f$ in $(U,c\delta)$ is a
Lebesgue point of $f$ in $(U,g)$, what proves the Lebesgue's
Differentiation Theorem for the Riemannian case.
\end{remark}

We will need the following variation of the Lebesgue
differentiation theorem:

\begin{proposition}
\label{lebesguemodificado} Let $(M,g)$ be a Riemannian manifold
and $f:M \rightarrow \mathbb R$ be a locally summable function.
Consider $x\in M$ a Lebesgue point of $f$ and let
$h:M\rightarrow\mathbb{R}$ be a continuous function. Then $x$ is
also a Lebesgue point of the product $h.f$.
\end{proposition}

\begin{proof}
\[
\frac{\int_{B_g(x,r)}\mid h(y)f(y)-h(x)f(x)\mid
dV_g(y)}{\mathrm{Vol}(B_g(x,r))}
\]
\[
\leq \frac{\int_{B_g(x,r)}\mid h(y)f(y)-h(y)f(x)\mid
dV_g(y)}{\mathrm{Vol}(B_g(x,r))}
\]
\[
+ \frac{\int_{B_g(x,r)}\mid h(y)f(x)-h(x)f(x)\mid
dV_g(y)}{\mathrm{Vol}(B_g(x,r))}
\]
which obviously goes to zero as $r$ goes to zero.
\end{proof}

Now we state Fatou's Lemma and Fubini's Theorem according to our
necessities.

\begin{theorem}
{\em (Fatou's Lemma)} Let $(M,g)$ be a Riemannian manifold and
consider a sequence of non-negative real valued functions
$\{f_i\}_{i \in \mathbb N}$. Assume that $\liminf\limits_{i
\rightarrow \infty}\Vert f_i\Vert_{L^1(M,g)}$ exists. Then
$\liminf\limits_{i\rightarrow \infty} f_i(x) $ exists for almost
every $x$, the function $\liminf\limits_{i \rightarrow \infty}
f_i$ is in $L^1(M,g)$, and we have that
\[
\int_M (\liminf\limits_{i \rightarrow \infty} f_i)(x) dV_g(x) \leq
\liminf\limits_{i \rightarrow \infty} \int_M f_i(x) dV_g(x) \leq
\liminf\limits_{i \rightarrow \infty} \Vert f_i \Vert_{L^1(M,g)}.
\]
\end{theorem}

\begin{theorem}
{\em (Fubini's Theorem)} Let $(M_1,g_1)$ and $(M_2,g_2)$ be two
smooth Riemannian manifolds and let $f:M_1\times M_2\rightarrow
\mathbb{R}$ be a Lebesgue summable function with respect to the
product Riemannian metric. Then
$f(x,\cdot):M_2\rightarrow\mathbb{R}$ are summable for almost
every $x\in M_1$ and $f(\cdot,y):M_1\rightarrow\mathbb{R}$ are
summable for almost every $y\in M_2$. Moreover
$\int_{M_2}f(x,y)dV_{g_2}(y)$ is summable on $M_1$,
$\int_{M_1}f(x,y)dV_{g_1}(x)$ is summable on $M_2$ and
\[
\int_{M_1\times M_2} f(x,y) dV_{g_1}(x). dV_{g_2}(y)
=\int_{M_1}\left(\int_{M_2}f(x,y)dV_{g_2}(y)\right)dV_{g_1}(x)
\]
\[
=\int_{M_2}\left(\int_{M_1}f(x,y)dV_{g_1}(x)\right)dV_{g_2}(y).
\]
\end{theorem}

Finally we remember the definition of a measure on a
sigma-algebra.

\begin{definition} \label{sigmaalgebra}
Let $M$ be an differentiable manifold. We say that a collection
$\mathcal L$ of subsets of $M$ is a {\em $\sigma$-algebra} in $M$
if

\begin{enumerate}

\item $M\in \mathcal L$.

\item $(M-A)\in\mathcal L$ whenever $A\in\mathcal L$.

\item If $\{A_i \in \mathcal L\}_{i \in \mathbb N}$ is a countable
family, then $\bigcup\limits_{i=1}^\infty A_i\in\mathcal L$.

\end{enumerate}

\end{definition}

\begin{theorem}
\label{menorsigmaalgebraexiste} For every collection $\mathcal M$
of subsets of $M$, there exists the smallest $\sigma$-algebra that
contains $\mathcal M$.
\end{theorem}

\begin{definition}
\label{borel} The smallest $\sigma$-algebra that contains all open
sets of $M$ is denoted by $\mathcal B(M)$ and its elements are
called {\em Borel sets} of $M$.
\end{definition}

\begin{definition}
\label{medidapositivaecomsinal} Let $\mathcal L$ be a
$\sigma$-algebra of $M$. A {\em positive measure} is a
non-negative function $\mathcal S:\mathcal L\rightarrow
[0,\infty]$ which is {\em countably additive}, that is, if $\{A_i
\in \mathcal L\}_{i \in \mathbb N}$ is a countable family of
pairwise disjoint subsets of $M$, then
\[
\mathcal S(\bigcup_{i=1}^\infty A_i,g)=\sum_{i=1}^\infty \mathcal
S(A_i,g).
\]

Finally we say that $\mathcal S:\mathcal L\rightarrow
[-\infty,\infty]$ is a {\em signed measure} if it is a difference
of two positive measures and if this difference is always well
defined (that is, the indetermination $\infty-\infty$ never
happens).

\end{definition}

\section{Non-regular tensor field spaces}
\label{espacoslpec0}

In this section, we develop the foundations of non-regular tensor
field spaces.

\begin{definition} \label{continuoint} We
say that a tensor field $T$ of type $(m,s)$ on a differentiable
manifold $M$ is continuous if the real valued function
$T(\varphi_1,\ldots,\varphi_m,v_1,\ldots,v_s):M\rightarrow
\mathbb{R}$ is continuous for every family
$\{\varphi_i:M\rightarrow T^*M;1\leq i\leq m\}$ of smooth 1-forms
and every family $\{v_i:M\rightarrow TM;1\leq i\leq s\}$ of smooth
vector fields. We denote the family of continuous tensor fields by
$C^0_{\mathrm{loc}}(M)$.
\end{definition}

Let $f:M\rightarrow \mathbb R$ be a measurable function. If we put
a smooth Riemannian metric $g$ on $M$, it makes sense to talk
about the space of functions such that $|f|^p$ is locally
integrable. But the local integrability of $\vert f \vert^p$ does
not depend on the choice of $g$. Therefore it makes sense to talk
about $L^p_{\mathrm{loc}}(M)$ spaces on a {\it differentiable}
manifold $M$. The same situation happens with tensor fields.

\begin{definition}
\label{localmentelpint} We say that a tensor field $T$ of type
$(m,s)$ on a differentiable manifold $M$ is in
$L^p_{\mathrm{loc}}(M)$ if the function
$T(\varphi_1,\ldots,\varphi_m,v_1,\ldots,v_s):M\rightarrow
\mathbb{R}$ is in $L^p_{\mathrm{loc}}(M)$ for every family
$\{\varphi_i:M\rightarrow T^*M;1\leq i\leq m\}$ of smooth 1-forms
and every family $\{v_i:M\rightarrow TM;1\leq i\leq s\}$ of smooth
vector fields. When $T \in L^1_{\mathrm{loc}}(M)$, we say that $T$
is locally summable.
\end{definition}

\begin{definition}
\label{tensornaoregular} We say that a tensor field is {\em
non-regular} if it is in either $L^p_{\mathrm{loc}}(M)$ or
$C^0_{\mathrm{loc}}(M)$.
\end{definition}

\begin{remark}
It is easy to see that every $C^0_{\mathrm{loc}}(M)$ tensor field
is an $L^p_{\mathrm{loc}}(M)$ tensor field. But we continue to
make the distinction between these spaces because they are
topologically different.
\end{remark}

We intend to define the mollifier smoothing $\widehat
T_\varepsilon$ of a non-regular tensor field $\widehat T$. In
order to do it, we put a smooth Riemannian metric $\widetilde g$
on $M$, which we call {\it the background metric}. Although most
of our results does not depend on $\widetilde g$, we have to
develop a theory about $L^p$ and $C^0$ tensor fields on a
Riemannian manifold $(M,\widetilde g)$. $\widetilde\nabla$ will
denote the Levi-Civita connection with respect to $\widetilde g$.

\begin{remark}
\label{notacaoparametrica} Symbols which comes with ``hat'', such
as $\widehat g$ and $\widehat T$, stand for non-regular tensor
fields. The symbol $\widetilde g$ stand for the background metric.
Symbols without any superscript, such as $g$ and $T$, stand for
general tensor fields.
\end{remark}

We denote the ``sup'' norm in several situations by $\Vert\cdot
\Vert_{C^0(M,\widetilde g)}$. When there is not any possibility of
misunderstandings, we simplify $\Vert\cdot \Vert_{C^0(M,\widetilde
g)}$ by $\Vert\cdot \Vert_{C^0}$. The following generalizations of
$L^p$ and $C^0$ spaces for the tensorial case are very natural:

\begin{definition}
\label{c0int} We say that the tensor field $T$ of type $(m,s)$ is
in $C^0(M,\widetilde g)$ if the function
$T(\varphi_1,\ldots,\varphi_m,v_1,\ldots,v_s):M\rightarrow
\mathbb{R}$ is continuous and bounded for every family of smooth
1-forms $\{\varphi_j:M\rightarrow T^*M;1\leq j\leq m,\Vert
\varphi_j\Vert_{C^0(M,\widetilde g)}\leq 1\}$ and every family of
smooth vector fields $\{v_i:M\rightarrow TM;1\leq i\leq s,\Vert
v_i\Vert_{C^0(M,\widetilde g)}\leq 1\}$.
\end{definition}

\begin{definition}
\label{espacolpint} Let $T$ be a tensor field of type $(m,s)$
defined on $(M,\widetilde g)$. We say that $T$ is in
$L^p(M,\widetilde g)$ if the function
$T(\varphi_1,\ldots,\varphi_m,v_1,\ldots,v_s):M\rightarrow
\mathbb{R}$ is in $L^p(M,\widetilde g)$ for every family of smooth
1-forms  $\{\varphi_j:M\rightarrow T^*M;1\leq j\leq m,\Vert
\varphi_j\Vert_{C^0(M,\widetilde g)}\leq 1\}$ and every family of
smooth vector fields $\{v_i:M\rightarrow TM;1\leq i\leq s,\Vert
v_i\Vert_{C^0(M,\widetilde g)}\leq 1\}$. Here we identify the
tensor fields that coincides almost everywhere. When $T \in
L^1(M,\widetilde g)$, we say that $T$ is summable.
\end{definition}

\begin{remark}
\label{condicaoequivalentelploc} Let $(M,\widetilde g)$ be a
differentiable manifold with background metric $\widetilde g$.
Observe that:

\begin{enumerate}

\item $T$ is continuous if and only if $T \in C^0(U, \widetilde
g\vert_U)$ for every $U \subset \subset M$.

\item $T \in L^p_{\mathrm{loc}}(M)$ if and only if $T \in L^p(U,
\widetilde g\vert_U)$ for every $U \subset \subset M$.

\end{enumerate}
The relationships $T \in C^0(U,\widetilde g\vert_U)$ and $T \in
L^p(U,\widetilde g\vert_U)$ do not depend on the background metric
$\widetilde g$ we choose. In particular, if $M$ is a closed
differentiable manifold, then the relationships $T \in
C^0(M,\widetilde g)$ and $T \in L^p(M,\widetilde g)$ do not depend
on the background metric $\widetilde g$ we choose.

\end{remark}

We will introduce complete norms on the tensor field spaces
$C^0(M,\widetilde g)$ and $L^p(M,\widetilde g)$, which generalize
the $C^0(M,\widetilde g)$ and $L^p(M,\widetilde g)$ norms for real
valued function spaces. The proof of the completeness of these
norms will be given in Theorems \ref{normac0} and \ref{normalp}
respectively. Although we are mainly interested to work with the
spaces $C^0_{\mathrm{loc}}(M)$ and $L^p_{\mathrm{loc}}(M)$,
Theorems \ref{normac0} and \ref{normalp} are interesting by
themselves and they are essential for the study of
$C^0_{\mathrm{loc}}(M)$ and $L^p_{\mathrm{loc}}(M)$ spaces.

Let us introduce a covering of $M$ by coordinate neighborhoods,
which will be useful afterwards. The general idea is that using
these coordinate systems, the local problem come out to be
essentially Euclidean. Afterwards we can glue each part using a
partition of the unity.

\begin{definition}
\label{almosteuclidean} Let $(M^n,\widetilde g)$ be a Riemannian
manifold, $U\subset M$ an open set and $\phi:U\rightarrow
\mathbb{R}^n$ a coordinate system. If $v$ is a smooth vector field
and $\varphi$ is a smooth 1-form defined on $M$, then their
restrictions to $U$ are written as
\begin{equation}
\label{defineai} \varphi\vert_U=\sum_{i=1}^n
a_i(x,\varphi).dx_i(x)
\end{equation}
and
\begin{equation}
\label{definebj} v\vert_U=\sum_{j=1}^n
b_j(x,v).\frac{\partial}{\partial x_j}(x)
\end{equation}
respectively. We say that $\phi$ is an almost Euclidean coordinate
system if
\begin{enumerate}
\item $\frac{1}{2}\leq \Vert dx_i(x)\Vert_{C^0(M,\widetilde g)}
\leq 2$ for every $i=1\ldots n$. \item $\frac{1}{2}\leq \Vert
\frac{\partial}{\partial x_j}(x)\Vert_{C^0(M,\widetilde g)} \leq
2$ for every $j=1\ldots n$. \item $\vert a_i(x,\varphi)\vert\leq
2$ for every 1-form such that $\Vert \varphi
\Vert_{C^0(M,\widetilde g)}\leq 1$. \item $\vert b_j(x,v)\vert\leq
2$ for every vector field such that $\Vert
v\Vert_{C^0(M,\widetilde g)}\leq 1$.
\end{enumerate}
\end{definition}

Observe that if we fix a point $x$, a normal coordinate system on
a sufficiently small neighborhood of $x$ is an almost Euclidean
coordinate system.

Let $\{U_y\subset\subset
M,\phi_y:U_y\rightarrow\mathbb{R}^n\}_{y\in M}$ be a family of
almost Euclidean coordinate system indexed by $y\in M$. Let
\begin{equation}
\label{coordenadaslocais} \{U_\omega\subset\subset
M,\phi_\omega:U_\omega\rightarrow\mathbb{R}^n\}_{\omega \in
\Lambda}
\end{equation}
be an almost Euclidean locally finite refinement of
$\{U_y,\phi_y\}_{y\in M}$, where $\Lambda$ is an index set. Denote
the cardinality of $\Lambda$ by $\vert\Lambda\vert$.

We denote the coordinate vector $\frac{\partial}{\partial x_i}$,
$1\leq i\leq n$, with respect to the coordinate system
$\{U_\omega,\phi_\omega\}$ at the point $x\in M$ by
$(\frac{\partial}{\partial x_i}(x))_\omega$. Its dual 1-form is
denoted by $(dx_i(x))_\omega$.

Take a partition of unity $\{\psi_\omega
:M\rightarrow\mathbb{R}\}_{\omega\in\Lambda}$ subordinated to
$\{U_\omega\}_{\omega\in\Lambda}$. Consider the 1-forms
$(\vartheta_{i})_\omega:M\rightarrow T^*M$ defined by
\begin{equation}
\label{varthetaiomega}
(\vartheta_i(x))_\omega=(\psi(x))_\omega.(dx_i(x))_\omega.
\end{equation}
Notice that $(\vartheta_i)_\omega$ is a smooth 1-form on $M$ for
every $i$ and $\omega$. In the same fashion, we can take smooth
vector fields $(u_i)_\omega:M\rightarrow T^*M$ defined by
\begin{equation}
\label{uiomega}
(u_i(x))_\omega=(\psi(x))_\omega.\left(\frac{\partial}{\partial
x_i}(x)\right)_\omega.
\end{equation}

\begin{proposition}
\label{formulaparav} Considering the notation given before, we
have the following formulas:
\[
\varphi(x)=\sum\limits_{\omega \in \Lambda} \sum\limits_{i=1}^n
(\psi(x))_\omega.(a_i(x,\varphi))_\omega.\left(dx_i(x)\right)_\omega
\]
\begin{equation}
\label{formulavarphi} =\sum\limits_{\omega \in \Lambda}
\sum\limits_{i=1}^n
(a_i(x,\varphi))_\omega.(\vartheta_i(x))_\omega
\end{equation}
and
\[
v(x)=\sum\limits_{\omega \in \Lambda} \sum\limits_{j=1}^n
(\psi(x))_\omega.(b_j(x,v))_\omega.\left(\frac{\partial}{\partial
x_j}\right)_\omega
\]
\begin{equation}
\label{formulav} =\sum\limits_{\omega \in \Lambda}
\sum\limits_{j=1}^n (b_j(x,v))_\omega.(u_j(x))_\omega.
\end{equation}
\end{proposition}

\begin{proof}
These formulas hold because the sum is locally finite.

\end{proof}

Let $T$ be a tensor field of type $(m,s)$ on $M$. Define the
functions
\begin{equation}
\label{Tijw} T_{i_1\ldots i_m j_1\ldots j_s\omega_1\ldots
\omega_{m+s}}=T((\vartheta_{i_1})_{\omega_1},\ldots,
(\vartheta_{i_m})_{\omega_m},(u_{j_1})_{\omega_{m+1}},
\ldots,(u_{j_s})_{\omega_{m+s}}),
\end{equation}
where the indexes $i_1,\ldots, i_m$ and $j_1,\ldots,j_s$ varies
between $1$ and $n$ and the indexes $\omega_1,\ldots,$ $
\omega_{m+s}$ are in $\Lambda$.

Proposition \ref{barT} define the function $\bar T$ which will be
useful afterwards.

\begin{proposition}
\label{barT} The function $\bar T:M\rightarrow \mathbb{R}$,
defined by
\begin{equation}
\label{equacaobarT} \bar
T(x)=2^{m+s}\sum_{(\omega_1,\ldots,\omega_{m+s})}
\sum_{(i_1,\ldots,i_m)} \sum_{(j_1,\ldots,j_s)} \vert T_{i_1\ldots
i_m j_1\ldots j_s\omega_1\ldots \omega_{m+s}}(x)\vert
\end{equation}
is greater than or equal to $T_{sup}:M\rightarrow\mathbb{R}$,
which is given by
\[
T_{sup}(x)=\sup_{\Vert \varphi_1\Vert_{C^0}\leq 1,\ldots,\Vert
u_s\Vert_{C^0}\leq 1} T(\varphi_1,\ldots,u_s)(x).
\]
\end{proposition}

\begin{proof}
The proof will be done for tensors of type $(1,1)$. The general
case follows similarly.
\[
T(\varphi,v)=
T\left(\sum_{\omega_1}\sum_{i}(a_i(x,\varphi))_{\omega_1}
(\vartheta_i(x))_{\omega_1},
\sum_{\omega_2}\sum_{j}(b_j(x,v))_{\omega_2}
(u_j(x))_{\omega_2}\right)
\]
\[
=\sum_{i}\sum_{j}\sum_{\omega_1,\omega_2}T
\left((a_i(x,\varphi))_{\omega_1}(\vartheta_i(x))_{\omega_1},
(b_j(x,v))_{\omega_2}(u_j(x))_{\omega_2}\right)
\]
\[
\leq 2^2 \sum_{i}\sum_{j}\sum_{\omega_1,\omega_2}\vert
T_{ij\omega_1\omega_2}\vert,
\]
because $\vert a_i(x,\varphi)\vert$ and $\vert
b_j(x,\varphi)\vert$ are less than or equal to 2.

\end{proof}

Now we define norms on $C^0(M,\widetilde g)$ and $L^p(M,\widetilde
g)$ which make them complete normed spaces.

\begin{theorem}
\label{normac0} Let $(M,\widetilde g)$ be a Riemannian manifold
and $C^0(M,\widetilde g)$ the space of $C^0$ tensor fields of type
$(m,s)$. Define $\Vert \cdot\Vert_{C^0(M,\widetilde
g)}:C^0(M,\widetilde g)\rightarrow \mathbb{R}$ by
\[
\Vert T\Vert_{C^0(M,\widetilde g)}=\sup\limits_{\Vert
\varphi_1\Vert_{C^0}\leq 1,\ldots,\Vert v_{s}\Vert_{C^0}\leq 1}
\Vert T(\varphi_1,\ldots,v_{s})\Vert_{C^0(M,\widetilde g)},
\]
where $\varphi_i$, $1\leq i\leq m$, and $v_j$, $1\leq j\leq s$,
are smooth $1$-forms and smooth vector fields defined on $M$
respectively. Then $\Vert \cdot\Vert_{C^0(M,\widetilde g)}$ is a
complete norm on $C^0(M,\widetilde g)$.
\end{theorem}

\begin{proof}
It is straightforward to see that $\Vert \cdot \Vert_{C^0}:= \Vert
\cdot \Vert_{C^0(M, \widetilde g)}$ is a norm on $C^0(M,\widetilde
g)$. Let us prove its completeness.

We will prove the completeness of $\Vert\cdot\Vert_{C^0}$ for
tensor spaces of type $(1,1)$. The general case follows similarly.
So let us begin with a Cauchy sequence $\{T_l\}_{l\in\mathbb{N}}$.
We will prove that there exists a tensor field $T$ such that
$\lim_{l\rightarrow\infty} \Vert T_l-T\Vert_{C^0}=0$.

Take a locally finite covering of $M$ as defined by
(\ref{coordenadaslocais}). Define the functions
$T_{ij\omega_1\omega_2}:M\rightarrow\mathbb{R}$ by
\begin{equation}
\label{tijw1w2a} T_{ij\omega_1\omega_2}=
\lim_{l\rightarrow\infty}T_l((\vartheta_i)_{\omega_1},(u_j)_{\omega_2})
\end{equation}
where $(\vartheta_i)_{\omega_1}$ and $(u_j)_{\omega_2}$ are given
by $(\ref{varthetaiomega})$ and $(\ref{uiomega})$ respectively.
Define
\begin{equation}
\label{formulaT} T(\varphi,v)=\sum_{i,j=1}^n\left(
\sum_{(\omega_1,\omega_2)\in\Lambda\times\Lambda}
(a_i)_{\omega_1}.(b_j)_{\omega_2}T_{ij\omega_1\omega_2}\right)
\end{equation}
where $\varphi$, $v$, $(a_i)_{\omega_1}$ and $(b_j)_{\omega_2}$
are related by Eqs. (\ref{defineai}) and (\ref{definebj}). It is
not difficult to see that $T$ is a tensor of type $(1,1)$. Let us
prove that $\lim_{l\rightarrow\infty}\Vert T_l-T\Vert_{C^0}=0$.

Let $\epsilon>0$. Then there exists a $N_\epsilon\in\mathbb N$
such that if $k,l>N_\epsilon$, then $\parallel T_k-T_l
\parallel_{C^0}<\epsilon/(32.n^2)$. This means that if we fix $i,j,\omega_1$
and $\omega_2$, then
\[
\vert T_{ij\omega_1\omega_2}-T_l((\vartheta_i(\cdot))_{\omega_1},
(u_j(\cdot))_{\omega_2}) \vert
\leq\frac{\epsilon}{8n^2}\psi_{\omega_1}\psi_{\omega_2}
\]
for every $l>N_\epsilon$.

Let $\varphi\in T^*M$ and $v\in TM$ be a smooth 1-form and a
smooth vector field respectively such that their $C^0$ norms are
less than or equal to one. Then
\[
\Vert T(\varphi,v)-T_l(\varphi,v) \Vert_{C^0}
\]
\[
=\left\Vert \sum_{i}\sum_{j}\sum_{\omega_1,\omega_2}
(T-T_l)\left((a_i(x,\varphi))_{\omega_1}(\vartheta_i(x))_{\omega_1},
(b_j(x,v))_{\omega_2}(u_j(x))_{\omega_2}\right)\right\Vert_{C^0}
\]
\[
\leq 4\left\Vert \sum_{i}\sum_{j}\sum_{\omega_1,\omega_2}
(T-T_l)\left((\vartheta_i(x))_{\omega_1},
(u_j(x))_{\omega_2}\right)\right\Vert_{C^0} \leq
\frac{\epsilon}{2}.
\]

Thus
\[
\Vert T-T_l \Vert_{C^0}=\sup_{\Vert \varphi \Vert_{C^0} \leq 1,
\Vert v\Vert_{C^0}\leq 1}\Vert (T-T_l)(\varphi,v)\Vert_{C^0}\leq
\frac{\epsilon}{2}<\epsilon
\]
whenever $l>N_\epsilon$ and $\lim\limits_{\varepsilon \rightarrow
0} \Vert T-T_l \Vert_{C^0} = 0$.

\end{proof}

\begin{theorem}
\label{normalp} Let $(M,\widetilde g)$ be a Riemannian manifold
and $L^p(M,\widetilde g)$ the space of $L^p$ tensor fields of type
$(m,s)$. Define $\Vert \cdot\Vert_{L^p(M,\widetilde
g)}:L^p(M,\widetilde g)\rightarrow \mathbb{R}$ by
\[
\Vert T\Vert_{L^p(M,\widetilde g)}:=\sup\limits_{\Vert
\varphi_1\Vert_{C^0}\leq 1,\ldots,\Vert v_s\Vert_{C^0}\leq
1}\left(\int_M \vert
T(\varphi_1,\ldots,v_s)(y)\vert^p.dV_{\widetilde
g}(y)\right)^{\frac{1}{p}},
\]
where $\varphi_i$, $1\leq i\leq m$, and $v_j$, $1\leq j\leq s$,
are smooth $1$-forms and smooth vector fields defined on $M$
respectively. Then $\Vert\cdot\Vert_{L^p(M,\widetilde g)}$ is a
complete norm on $L^p(M,\widetilde g)$.
\end{theorem}

\begin{proof} The proof will be done for tensor spaces of type $(1,1)$. The
general case follows similarly.

It is straightforward that the properties $\Vert
T_1+T_2\Vert_{L^p(M,\widetilde g)}\leq \Vert
T_1\Vert_{L^p(M,\widetilde g)}+\Vert T_2\Vert_{L^p(M,\widetilde
g)}$ and $\Vert c.T\Vert_{L^p(M,\widetilde g)}=\vert c\vert.\Vert
T\Vert_{L^p(M,\widetilde g)}$ holds for every $T, T_1, T_2\in
L^p(M,\widetilde g)$ and $c\in \mathbb{R}$.

Suppose that $\Vert T\Vert_{L^p(M,\widetilde g)}=0$. Then $T=0$
a.e.. In fact, take a locally finite covering
$\{U_\omega,\phi_\omega\}_{\omega\in\Lambda}$ of $M$ as defined by
(\ref{coordenadaslocais}) and consider $(\vartheta_i)_{\omega}$
and $(u_j)_{\omega}$ as defined in (\ref{varthetaiomega}) and
(\ref{uiomega}) respectively. Then
$T((\vartheta_i)_{\omega_1},(u_j)_{\omega_2})=0$ a.e. for every
$\{i,j,\omega_1,\omega_2\}\in \{1,\ldots,n\} \times \{1,\ldots,n\}
\times \Lambda \times \Lambda$. But using (\ref{formulavarphi})
and (\ref{formulav}), we can see that $T(\varphi,v)=0$ a.e. for
every smooth 1-form $\varphi$ and every smooth vector field $v$.
This implies that $T=0$ a.e.. Therefore
$\Vert\cdot\Vert_{L^p(M,\widetilde g)}$ is a norm defined on
$L^p(M,\widetilde g)$.

In order to prove the completeness of $\Vert
\cdot\Vert_{L^p(M,\widetilde g)}$, let us begin with a Cauchy
sequence $\{T_l\}_{l\in\mathbb{N}}$. We will prove that it
converges to the tensor field $T$ defined by (\ref{formulaT}).

\

{\it Claim 1:} $T \in L^p(M,\widetilde g)$, that is, $\Vert
T\Vert_{L^p(M,\widetilde g)} <\infty$.

\

The general idea here is to built a uniformly bounded sequence
$\{Q_d\}_{d \in \mathbb N}$ of tensor fields in $L^p(M,\widetilde
g)$ such that it converges pointwise to $T$ as $d$ goes to
infinity. Then using Fatou's Lemma, we will be able to prove that
$\Vert T\Vert_{L^p(M,\widetilde g)} <\infty$. The tensor field
$Q_d$ is essentially $T$ restricted to a compact subset $K_d
\subset M$, where $\{K_d\}_{d \in \mathbb N}$ is an increasing
sequence of compact sets (that is, $K_d \subset K_{d+1}$ for every
$d$) such that $\bigcup_{d\in \mathbb N} K_d = M$. The
decomposition $T = Q_d + (T - Q_d)$ allow us to split $T$ in a
``compact part'' $Q_d$ and in its complement $T-Q_d$, what is
typical in this kind of problem. Let us formalize the idea:

Let $\epsilon=\frac{1}{d}>0$, where $d\in\mathbb{N}$. Then there
exists a $N_d\in\mathbb{N}$ such that $\Vert
T_{l_1}-T_{l_2}\Vert_{L^p(M,\widetilde g)}<\frac{\epsilon}{2}$
whenever $l_1,l_2\geq N_d$.

Observe that there exists a compact subset $K_d\subset M$ such
that
\[
\left\Vert (T_{N_d})\mid_{(M-K_d)}\right\Vert_{L^p(M,\widetilde
g)}<\frac{\epsilon}{2}.
\]
In fact, if this is not the case, then it would not be difficult
to built a vector field $v$ and a 1-form $\varphi$ such that
$\Vert v \Vert_{C^0} \leq 1$, $\Vert \varphi \Vert_{C^0}\leq 1$
and $\left\Vert T_{N_d}(v,\varphi)\right\Vert_{L^p(M,\widetilde
g)} > \left\Vert T_{N_d}\right\Vert_{L^p(M,\widetilde g)}$, which
would give a contradiction. Then

\begin{equation}
\label{complementodecompacto} \left\Vert
(T_l)\mid_{(M-K_d)}\right\Vert_{L^p(M,\widetilde g)}<\epsilon
\end{equation}
for every $l\geq N_d$. Moreover we can assume that $\{K_d \}_{d\in
\mathbb N}$ is an increasing sequence of compact subsets of $M$
such that $\bigcup_{d\in \mathbb N} K_d = M$

Take a locally finite covering
$\{U_\omega,\phi_\omega\}_{\omega\in\Lambda}$ of $M$ as defined by
(\ref{coordenadaslocais}) and a partition of unity $\{\psi_\omega
:M\rightarrow\mathbb{R}\}_{\omega\in\Lambda}$ subordinated to
$\{U_\omega\}_{\omega\in\Lambda}$. Consider
$(\vartheta_i)_{\omega}$ and $(u_j)_{\omega}$ as defined in
(\ref{varthetaiomega}) and (\ref{uiomega}) respectively. Let
$\{U_{\widetilde\omega}\}_{\widetilde\omega\in \Lambda_d}\subset
\{U_\omega\}_{\omega\in\Lambda}$ be the set of coordinate
neighborhoods such that $U_{\widetilde\omega}\cap K_d \not
=\emptyset$. It is not difficult to see that
$\{U_{\widetilde\omega}\}_{\widetilde\omega\in \Lambda_d}$ is a
finite set, because $\{U_\omega\}_{\omega\in\Lambda}$ is a locally
finite covering of $M$. Observe also that $\Lambda_d \subset
\Lambda_{d+1}$.

Let $T_{ij\omega_1\omega_2}:M\rightarrow\mathbb{R}$ be the
functions defined by Eq. (\ref{tijw1w2a}). Define the sequence of
tensor fields $\{Q_d\}_{d\in\mathbb{N}}$ by
\[
Q_d(\varphi,v)=\sum_{i,j=1}^n \left(\sum_{(\omega_1,\omega_2)\in
\Lambda_d\times\Lambda_d}
(a_i)_{\omega_1}.(b_j)_{\omega_2}T_{ij\omega_1\omega_2}\right),
\]
where $\varphi$ is given by (\ref{formulavarphi}) and $v$ is given
by (\ref{formulav}). Notice that $Q_d\rightarrow T$ as
$d\rightarrow\infty$ everywhere.

The idea here is to prove that $\Vert Q_d\Vert_{L^p}\leq C$ for
some fixed constant $C$ and using the Fatou's Lemma we will be
able to prove that $\Vert T\Vert_{L^p}\leq C$. Fix a 1-form
$\varphi$ and a vector field $v$ such that $\Vert \varphi
\Vert_{C^0}, \Vert v \Vert_{C^0}\leq 1$. We have the following
estimate:
\[
\left\Vert T_{l}(\varphi,v)-Q_d(\varphi,v)
\right\Vert_{L^p(M,\widetilde g)}
\]
\[
\leq \left\Vert \sum_{i,j=1}^n
\left(\sum_{(\omega_1,\omega_2)\in\Lambda_d\times\Lambda_d}
(a_i)_{\omega_1}.(b_j)_{\omega_2}(T_l((\vartheta_i)_{\omega_1},
(u_j)_{\omega_2})-T_{ij\omega_1\omega_2})\right)\right\Vert_{L^p(M,\widetilde
g)}
\]
\begin{equation}
\label{tnd-sd} +\left\Vert\sum_{i,j=1}^n
\left(\sum_{(\omega_1,\omega_2)
\not\in\Lambda_d\times\Lambda_d}(a_i)_{\omega_1}.(b_j)_{\omega_2}
(T_l((\vartheta_i)_{\omega_1},(u_j)_{\omega_2}))\right)
\right\Vert_{L^p(M,\widetilde g)}.
\end{equation}
Due to (\ref{complementodecompacto}), the last term of the
right-hand-side of (\ref{tnd-sd}) is less than $\epsilon$ for
every $l\geq N_d$.

For the first term of the right-hand-side of (\ref{tnd-sd}), there
exists $N_d^\prime \geq N_d$ such that
\[
\Vert
(T_l((\vartheta_i)_{\omega_1},(u_j)_{\omega_2})-T_{ij\omega_1\omega_2})
\Vert_{L^p(M,\widetilde
g)}<\frac{\epsilon}{4.n^2.\vert\Lambda_d\vert^2}
\]
for every $(i,j,\omega_1,\omega_2)\in \{1,\ldots,n\} \times
\{1,\ldots,n\} \times \Lambda_d\times\Lambda_d$ and $l\geq
N_d^\prime$. Therefore the first term of the right-hand-side of
(\ref{tnd-sd}) is less than $\epsilon$ for $l\geq N_d^\prime$ what
implies that the right-hand-side of (\ref{tnd-sd}) is less than
$2\epsilon$.

These facts implies that
\[
\begin{array}{l}
\Vert \vert Q_d(\varphi,v) \vert^p \Vert^{1/p}_{L^1(M,\widetilde
g)}=\Vert Q_d(\varphi,v)\Vert_{L^p(M,\widetilde g)}\leq \sup_{l
\geq N_d^\prime}\Vert
T_l(\varphi,v)\Vert_{L^p(M,\widetilde g)} + 2\epsilon \leq \\
\sup_{l \in \mathbb N}\Vert T_l\Vert_{L^p(M,\widetilde g)} +
2<\infty
\end{array}
\]
for every $d\in\mathbb{N}$. Now we can use the Fatou's Lemma and
the pointwise convergence $\lim\limits_{d\rightarrow\infty} \vert
Q_d(\varphi,v) \vert^p$ $= \vert T(\varphi,v) \vert^p$ in order to
conclude that $\Vert \vert T \vert^p \Vert^{1/p}_{L^1(M,\widetilde
g)}=\Vert T\Vert_{L^p(M,\widetilde g)}$ is finite, what settles
Claim 1.

\

{\it Claim 2:} The sequence $\{T_l\}_{l\in\mathbb{N}}$ converges
to $T$ in $L^p(M,\widetilde g)$.

\

Some ideas that we used to prove the finiteness of $\Vert
T\Vert_{L^p(M,\widetilde g)}$ will appear here again.

For every $d^\prime \in \mathbb N$, set $\epsilon =
\frac{1}{d^\prime}>0$. Notice that $\{T-T_l\}_{l\in\mathbb{N}}$ is
a Cauchy sequence in $L^p(M,\widetilde g)$. In the same way we did
in Claim 1, there exist $N_{d^\prime} \in \mathbb N$ and a compact
set $K_{d^\prime} \subset M$ such that
\[
\Vert (T - T_l)\vert_{(M - K_{d^\prime})} \Vert_{L^p(M,\widetilde
g)} < \epsilon
\]
for every $l \geq N_{d^\prime}$ (Compare with
(\ref{complementodecompacto})). We can take the sequence
$\{K_d\}_{d \in \mathbb N}$ such that $K_d \subset K_{d+1}$ for
every $d \in \mathbb N$. Denote by
$\{U_{\widetilde\omega}\}_{\widetilde\omega\in
\Lambda_{d^\prime}}\subset \{U_\omega\}_{\omega\in\Lambda}$ the
set of coordinate neighborhoods such that
$U_{\widetilde\omega}\cap K_{d^\prime} \not =\emptyset$.

Let $\varphi$ be a smooth $1$-form such that $\Vert
\varphi\Vert_{C^0}\leq 1$ and let $v$ be a smooth vector field
such that $\Vert v\Vert_{C^0}\leq 1$. Then
\[
\left\Vert\sum_{i,j=1}^n
\left(\sum_{(\omega_1,\omega_2)\not\in\Lambda_{d^\prime}\times
\Lambda_{d^\prime}}(a_i)_{\omega_1}.(b_j)_{\omega_2}
((T-T_l)((\vartheta_i)_{\omega_1},(u_j)_{\omega_2}))\right)
\right\Vert_{L^p(M,\widetilde g)}<\epsilon.
\]

Now observe that
\[
\Vert T(\varphi,v)-T_l(\varphi,v) \Vert_{L^p(M,\widetilde g)}
\]
\[
\leq \left\Vert\sum_{i,j=1}^n
\left(\sum_{(\omega_1,\omega_2)\in\Lambda_{d^\prime}\times
\Lambda_{d^\prime}}(a_i)_{\omega_1}.(b_j)_{\omega_2}((T-T_l)
((\vartheta_i)_{\omega_1},(u_j)_{\omega_2}))\right)\right
\Vert_{L^p(M,\widetilde g)}
\]
\[
+ \left\Vert\sum_{i,j=1}^n
\left(\sum_{(\omega_1,\omega_2)\not\in\Lambda_{d^\prime}\times
\Lambda_{d^\prime}}(a_i)_{\omega_1}.(b_j)_{\omega_2}((T-T_l)
((\vartheta_i)_{\omega_1},(u_j)_{\omega_2}))\right)\right\Vert_{L^p(M,\widetilde
g)}.
\]

Here we repeat the procedure we used to control the first term of
the right-hand-side of (\ref{tnd-sd}): There exist
$N_{d^\prime}^\prime\in\mathbb{N}$ such that
\[
\left\Vert\sum_{i,j=1}^n \left(\sum_{(\omega_1,\omega_2)
\in\Lambda_{d^\prime}\times\Lambda_{d^\prime}}(a_i)_{\omega_1}.
(b_j)_{\omega_2}((T-T_l)((\vartheta_i)_{\omega_1},(u_j)_{\omega_2}))
\right)\right\Vert_{L^p(M,\widetilde g)}<\epsilon
\]
for $l\geq N_{d^\prime}^\prime$. Observe that this inequality does
not depend on $v$ and $\varphi$. Then
\[
\Vert T-T_l \Vert_{L^p(M,\widetilde g)}<2\epsilon
\]
for $l\geq N_{d^\prime}^\prime$. Therefore
$\lim\limits_{l\rightarrow\infty} T_l=T$ in $L^p(M,\widetilde g)$.

\end{proof}

The tensor field spaces $L^p_{\mathrm{loc}}(M)$ and
$C^0_{\mathrm{loc}}(M)$ can be topologized in the same way as
their correspondent function spaces:

\begin{definition}
\label{convergeemcoloc} Let $M$ be a differentiable manifold. We
say that a sequence $\{T_i\}_{i \in \mathbb N}$ of tensor fields
converges to $T$ in $C^0_{\mathrm{loc}}(M)$ if
\[
\lim_{i \rightarrow \infty}\Vert (T_i)\vert_U - (T) \vert_U
\Vert_{C^0(U,\widetilde g)} = 0
\]
for every background metric $\widetilde g$ and every open set
$U\subset \subset M$.

Analogously we say that a one-parameter family of tensor fields
$\{ T_\varepsilon\}_{\varepsilon >0}$ converges to $T$ in
$C^0_{\mathrm{loc}}(M)$ as $\varepsilon$ goes to zero if
\[
\lim_{\varepsilon \rightarrow 0}\Vert (T_\varepsilon)\vert_U - (T)
\vert_U \Vert_{C^0(U,\widetilde g)} = 0
\]
for every background metric $\widetilde g$ and every open set
$U\subset \subset M$.

Observe that these convergences do not depend on the choice of
$\widetilde g$.

\end{definition}

\begin{definition}
\label{convergeemlploc} Let $M$ be a differentiable manifold. We
say that a sequence $\{T_i\}_{i \in \mathbb N}$ of tensor fields
converges to $T$ in $L^p_{\mathrm{loc}}(M)$ if
\[
\lim_{i \rightarrow \infty}\Vert (T_i)\vert_U - (T) \vert_U
\Vert_{L^p(U,\widetilde g)} = 0
\]
for every background metric $\widetilde g$ and every open set
$U\subset \subset M$.

Analogously we say that a one-parameter family of tensor fields
$\{ T_\varepsilon\}_{\varepsilon >0}$ converges to $T$ in
$L^p_{\mathrm{loc}}(M)$ as $\varepsilon$ goes to zero if
\[
\lim_{\varepsilon \rightarrow 0}\Vert (T_\varepsilon)\vert_U - (T)
\vert_U \Vert_{L^p(U,\widetilde g)} = 0
\]
for every background metric $\widetilde g$ and every open set
$U\subset \subset M$.

Observe that these convergences do not depend on the choice of
$\widetilde g$.

\end{definition}

\section{The mollifier smoothing of a non-regular tensor field}
\label{mollifiergeral}

Let us remember the definition of mollifier smoothing of a locally
summable function $\widehat f:\widehat
U\subset\mathbb{R}^n\rightarrow\mathbb{R}$.

Let $\widetilde\eta:\mathbb{R}^n\rightarrow \mathbb{R}$ be the
$C^\infty$ function defined by
\[
\widetilde\eta(x):= \left\{
\begin{array}{ll}
C\exp\left(\frac{1}{\Vert x\Vert^2-1}\right) & \mathrm{if}\hspace{5mm}\Vert x\Vert <1 \\
0 & \mathrm{if}\hspace{5mm}\Vert x\Vert \geq 1,
\end{array}
\right.
\]
where $C$ is a constant such that
$\int_{\mathbb{R}^n}\widetilde\eta(x)dV=1$. Now define
\begin{equation}
\label{mollifiereuclidiano}
\bar\eta(x,y,\varepsilon):=\frac{1}{\varepsilon^n}\widetilde\eta
\left(\frac{x-y}{\varepsilon}\right).
\end{equation}

Let $U \subset \subset \widehat U$ and define the mollifier
smoothing $\widehat f_\varepsilon:U \rightarrow \mathbb{R}^n$ of
$\widehat f$ by
\[
\widehat f_\varepsilon(x)=\int_{\widehat
U}\bar\eta(x,y,\varepsilon)\widehat f(y)dV(y)
\]
where $\varepsilon < \mathrm{dist}(U,\mathbb R^n - \widehat U)$.

\begin{theorem}
\label{convergenciasmollifierfuncao}

\begin{enumerate}

\item $\widehat f_\varepsilon$ is $C^\infty(U)$.

\item $\widehat f_\varepsilon\rightarrow \widehat f$ a.e. as
$\varepsilon\rightarrow 0$.

\item If $\widehat f$ is a continuous function, then $\widehat
f_\varepsilon\rightarrow f$ uniformly on compact subsets of
$\widehat U$.

\item If $\widehat f\in L^p_{\mathrm{loc}}(\widehat U)$, then
$\widehat f_\varepsilon \rightarrow \widehat f$ in
$L^p_{\mathrm{loc}}(\widehat U)$.

\end{enumerate}
\end{theorem}

\begin{proof}
See \cite{7}.

\end{proof}

Let $M$ be a differentiable manifold. Let us define the mollifier
smoothing of a non-regular tensor field $\widehat T \in T^{m,s}M$.

When we define a mollifier smoothing of a function $\widehat
f:\widehat U\subset\mathbb{R}^n\rightarrow\mathbb{R}$, the
Euclidean metric on $\widehat U$ takes part in the process. In
order to define the mollifier smoothing of tensor fields on a
differentiable manifold, we introduce the background metric
$\widetilde g$ which is a smooth Riemannian metric on $M$.

Let $U \subset \subset M$ be an open set (This includes the case
$U=M$ if $M$ is a closed differentiable manifold). Define the
injectivity radius of $U$ by
\begin{equation}
\label{raiodeinjetividade} \mathrm{inj}(U,\widetilde g)=\inf_{x
\in (U,\widetilde g)} \mathrm{inj}(x)
\end{equation}
where $\mathrm{inj}(x)$ is the injectivity radius of $x \in
(M,\widetilde g)$. We say that a function $\eta:U\times M\times
(0,\mathrm{inj}(U,\widetilde g))$ is a mollifier if

\begin{enumerate}

\item $\eta$ is a smooth function.

\item $\eta(x,\cdot,\varepsilon):M\rightarrow \mathbb{R}$ has its
support in $\bar B(x,\varepsilon)$.

\item $\int_M \eta(x,y,\varepsilon).dV_{\widetilde g}(y)=1$ for
every $x\in U$ and $\varepsilon\in(0,\mathrm{inj}(U,\widetilde
g))$ fixed.

\end{enumerate}

Observe that $\bar\eta$ defined by (\ref{mollifiereuclidiano})
satisfies the conditions above for the Euclidean case. Moreover we
can adapt $\bar\eta$ for a Riemannian manifold in the following
fashion: Define
\[
\check\eta(x,y,\varepsilon)= \left\{
\begin{array}{l}
\exp\left(\frac{1}{\left(\frac{\mathrm{dist}(x,y)}{\varepsilon}\right)^2-1}\right)
\mathrm{\;\;if\;\;dist}(x,y)<\varepsilon < \mathrm{inj}(U,\widetilde g) \\
0 \mathrm{\;\;if\;\;dist}(x,y)\geq \varepsilon.
\end{array}
\right.
\]
The function $\check\eta$ is clearly $C^\infty$. Now define
\begin{equation}
\label{molifierpadrao}
\eta(x,y,\varepsilon)=\frac{\check\eta(x,y,\varepsilon)}{\int_M
\check\eta(x,y,\varepsilon) dV_{\widetilde g}(y)}.
\end{equation}
It is not difficult to see that the function $\eta$ defined above
satisfies the properties of a mollifier on $(M,\widetilde g)$ and
it will be called the {\em standard mollifier}. Moreover the
following proposition holds:

\begin{proposition}
\label{normadeeta} Let $\eta:U \times M \times
(0,\mathrm{inj}(U,\widetilde g)) \rightarrow \mathbb R$ be the
standard mollifier. Then there exist a positive constant $\Vert
\eta \Vert$ such that
\[
\eta(x,y,\varepsilon)\leq \frac{\Vert \eta \Vert}{\varepsilon^n},
\]
where $n$ is the dimension of $M$.
\end{proposition}

\begin{proof}
If we observe that

\begin{enumerate}

\item The proposition is trivially true for Euclidean spaces;

\item Every smooth Riemannian manifold is ``locally almost
Euclidean'';

\item It is enough to prove the proposition for small values of
$\varepsilon$;

\item $\bar U$ is compact,

\end{enumerate}
then it is not difficult to prove the proposition.

\end{proof}

Suppose that there exist a unique minimizing geodesic $\gamma$
connecting $x\in (M,\widetilde g)$ and $y \in (M,\widetilde g)$.
We denote the parallel transport between the tensor spaces
$T^{m,s}_xM$ and $T^{m,s}_yM$ through $\gamma$ by
$\widetilde{\tau}_{x,y}$. If we have $\xi\in T^{m,s}_xM$, then we
can define a smooth tensor field $\xi(x ;\widetilde g)$ on the
geodesic ball $B(x,\varepsilon)$ with $\varepsilon <
\mathrm{inj}(x)$ using parallel transport by
\[
\xi(x ;\widetilde g)(y):=\widetilde\tau_{x,y}(\xi),\; y \in
B(x,\varepsilon).
\]

Let $\widehat T$ be a non-regular tensor field in $T^{m,s}M$. We
define the mollifier smoothing $\widehat T_\varepsilon$ of
$\widehat T$ with respect to $\widetilde g$ as follows:

\begin{definition}
\label{definicaosuavizacao} Let $(M,\widetilde g)$ be a
differentiable manifold with a smooth background metric
$\widetilde g$. Consider an open set $U \subset \subset M$ and let
$\varepsilon < \mathrm{inj}(U, \widetilde g)$. We define the
mollifier smoothing $\widehat T_\varepsilon \in T^{m,s}U$ of a
non-regular tensor field $\widehat T\in T^{m,s}M$  with respect to
$\widetilde g$ by
\begin{equation}
\label{formulasuavizacao}
\begin{array}{c}
\widehat T_\varepsilon(\varphi_1,\ldots,\varphi_m,v_1,\ldots,v_s)(x) \\
=\int\limits_{(B(x,\varepsilon),\widetilde
g)}\eta\left(x,y,\varepsilon\right). \widehat
T(\varphi_1(x;\widetilde g),\ldots,\varphi_m(x;\widetilde g),
v_1(x;\widetilde g),\ldots,v_s(x;\widetilde g))(y).dV_{\widetilde
g}(y).
\end{array}
\end{equation}
where $\{\varphi_i \in T^*_x M ;1\leq i\leq m\}$ and $\{v_i\in T_x
M;1\leq i\leq s\}$.
\end{definition}

\begin{remark}
Sometimes it is more useful to see Definition
(\ref{formulasuavizacao}) as
\[
\begin{array}{c}
\widehat T_\varepsilon(\varphi_1,\ldots,\varphi_m,v_1,\ldots,v_s)(x) \\
=\int\limits_{(M,\widetilde g)}\eta\left(x,y,\varepsilon\right).
\widehat T(\varphi_1(x;\widetilde g),\ldots,\varphi_m(x;\widetilde
g), v_1(x;\widetilde g),\ldots,v_s(x;\widetilde
g))(y).dV_{\widetilde g}(y)
\end{array}
\]
(Notice that $\eta(x,\cdot,\varepsilon)$ has its support in
$B(x,\varepsilon)$).
\end{remark}

\begin{theorem}
\label{metricasuave} The function $\widehat T_\varepsilon$ defined
in Eq. (\ref{formulasuavizacao}) is a smooth tensor field of type
$(m,s)$ on $U$. Moreover, if $\widehat T$ is a non-regular
Riemannian metric on $M$, then $\widehat T_\varepsilon$ is a
smooth Riemannian metric on $M$.
\end{theorem}

\begin{proof}
In order to prove that $\widehat T_\varepsilon$ is a tensor field
of type $(m,s)$, it is enough to prove that it is linear on each
variable, which is straightforward. If $\widehat T$ is a
non-regular Riemannian metric on $M$, then we have to prove that
$\widehat T_\varepsilon$ is symmetric and positive definite, which
is also straightforward.

Let us prove the smoothness of $\widehat T_\varepsilon$. The
tensor field $\widehat T_\varepsilon$ is smooth if and only if the
function $\widehat
T_\varepsilon(\varphi_1,\ldots,\varphi_m,v_1,\ldots,v_s)$ is
smooth for every family of smooth 1-forms $\{\varphi_i:
M\rightarrow T^*M,1 \leq i \leq n\}$ and every family of smooth
vector fields $\{v_i: M\rightarrow TM,1 \leq i \leq n\}$. Consider
Definition (\ref{formulasuavizacao}) and notice that
\begin{itemize}
\item $\eta(\cdot,y,\varepsilon):M\rightarrow\mathbb{R}$ is a
smooth function. \item Fix $y\in M$. Due to Theorem
\ref{transporteparalelosuave}, $v_i(\cdot\: ;\widetilde g)(y)=
:B(y,\varepsilon)\rightarrow T_yM$ is smooth for every $1\leq
i\leq s$ (The smoothness of $\varphi_i(\cdot \: ;\widetilde g)(y)$
follows similarly). \item The partial derivatives with respect to
$x$ in (\ref{formulasuavizacao}) can be taken inside the integral
infinitely because $\widehat T$ is multilinear with respect to its
variables.
\end{itemize}

Therefore $\widehat T_\varepsilon$ is smooth, what settles the
theorem.

\end{proof}

We are interested to know when
$\lim\limits_{\varepsilon\rightarrow 0}\widehat
T_\varepsilon=\widehat T$ is true. The theorem below states that
this happens in several situations.

\begin{theorem}
\label{convergencias} Let $\widehat T$ be a locally summable
tensor field of type $(m,s)$ on a Riemannian manifold
$(M,\widetilde g)$. Then we have the following convergences:
\begin{enumerate}

\item $\widehat T_\varepsilon\rightarrow \widehat T$ a.e. as
$\varepsilon\rightarrow 0$.

\item If $\widehat T$ is continuous, then $\widehat
T_\varepsilon\rightarrow \widehat T$ in $C^0_{\mathrm{loc}}(M)$ as
$\varepsilon \rightarrow 0$.

\item If $1\leq p<\infty$ and $\widehat T\in
L^p_{\mathrm{loc}}(M)$, then $\widehat T_\varepsilon\rightarrow
\widehat T$ in $L^p_{\mathrm{loc}}(M)$ as $\varepsilon \rightarrow
0$.

\end{enumerate}

\end{theorem}

\begin{proof}

\

\

(1) $\widehat T_\varepsilon\rightarrow \widehat T$ a.e. as
$\varepsilon\rightarrow 0$.

\

Let $B(x_0,r) \subset \subset (M,\widetilde g)$ be a geodesic
ball, with $r$ less than the injectivity radius at $x_0$. Let
$\{\beta_1,\beta_2,\ldots,\beta_n\}$ be a family of smooth one
forms on $B(x_0,r)$ such that
$\{\beta_1(x),\beta_2(x),\ldots,\beta_n(x)\}$ is a basis of
$T^*_{x}M$ for every $x \in B(x_0,r)$. Analogously let
$\{w_1,w_2,\ldots,w_n\}$ be a family of smooth vector fields on
$B(x_0,r)$ such that $\{w_1(x),w_2(x),$ $\ldots,w_n(x)\}$ is a
basis of $T_{x}M$ for every $x \in B(x_0,r)$. Fix an $m$-tuple
$(\beta_{k_1},\beta_{k_2},\ldots,\beta_{k_m})$ and an $s$-tuple
$(w_{j_1},w_{j_2},\ldots,w_{j_s})$ (the elements of each family
can eventually appear more than once). The Lebesgue's
Differentiation Theorem states that
\begin{equation}
\label{lebesgue} \lim\limits_{\varepsilon\rightarrow
0}\frac{\int\limits_{B(x,\varepsilon)} \vert \widehat
T(\beta_{k_1},\ldots,w_{j_s})(y)-\widehat
T(\beta_{k_1},\ldots,w_{j_s})(x)\vert .dV_{\widetilde
g}(y)}{\mathrm{Vol}(B(x,\varepsilon))}=0
\end{equation}
for almost every $x \in B(x_0,r)$. Denote by $L$ the set of points
$x\in B(x_0,r)$ such that (\ref{lebesgue}) holds for every pair of
families $(\beta_{k_1},\ldots,\beta_{k_m})$ and
$(w_{j_1},\ldots,w_{j_s})$ ($L$ is the set of Lebesgue points of
$\widehat T\vert_{B(x_0,r)}$). We have that $B(x_0,r)-L$ has
measure zero.

We claim that $\widehat T_\varepsilon (x)\rightarrow \widehat
T(x)$ as $\varepsilon \rightarrow 0$ if $x\in L$. In fact,
consider a family $(\varphi_1,\ldots,\varphi_m)$ of smooth 1-forms
and a family $(v_1,\ldots,v_s)$ of smooth vector fields on $M$.
Let $\varepsilon>0$ such that $B(x,\varepsilon) \subset B(x_0,r)$.
It follows that
\begin{eqnarray}
& \left\vert \widehat
T_\varepsilon(\varphi_1,\ldots,v_s)(x)-\widehat
T(\varphi_1,\ldots,v_s)(x)\right\vert & \nonumber
\\
& =\left\vert \int_{M} \eta(x,y,\varepsilon)\left[\widehat
T(\varphi_1(x;\widetilde g),\ldots,v_s(x;\widetilde
g))(y)-\widehat T(\varphi_1,\ldots,v_s)(x)\right].dV_{\widetilde
g}(y) \right\vert. & \label{termo2T} \end{eqnarray}

For each point $y\in B(x,\varepsilon)$, we denote $\varphi_l(x;
\widetilde g)(y)=\sum_{k=1}^n h^*_{lk}(y)\beta_k(y)$ and
$v_i(x;\widetilde g)(y)=\sum_{j=1}^n h_{ij}(y)w_j(y)$, where
$1\leq l\leq m$ and $1\leq i\leq s$. Observe that $\widetilde
h^*_{lk}$ and $\widetilde h_{ij}$ are smooth functions for all the
indexes $l,k,i,j$. Then we have that

\begin{eqnarray}
& \widehat T(\varphi_1(x;\widetilde g),\ldots,
v_s(x;\widetilde g))(y) & \nonumber \\
& = \sum\limits_{j_1=1}^n\ldots\sum\limits_{k_m=1}^n
h^*_{1k_1}(y)\ldots h^*_{mk_m}(y) h_{1j_1}(y)\ldots
h_{sj_s}(y)\widehat T(\beta_{k_1},\ldots,v_{j_s})(y) &
\end{eqnarray}
for every $y\in B(x,\varepsilon)$ because $\widehat T$ is
multilinear. Therefore the integrand of the right-hand-side of Eq.
(\ref{termo2T}) is a sum of terms like
\begin{eqnarray}
&\eta(x,y,\varepsilon)h^*_{1k_1}(y)\ldots h^*_{mk_m}(y)h_{1j_1}(y)
\ldots h_{sj_s}(y)\widehat T(\beta_{k_1},\ldots,w_{j_s})(y) & \nonumber \\
& - \eta(x,y,\varepsilon)h^*_{1k_1}(x)\ldots
h^*_{mk_m}(x)h_{1j_1}(x)\ldots h_{sj_s}(x)\widehat
T(\beta_{k_1},\ldots,w_{j_s})(x), & \nonumber
\end{eqnarray}
which, for brevity, we write as
\[
\eta(x,y,\varepsilon)\widehat T^h(y)-\eta(x,y,\varepsilon)\widehat
T^h(x).
\]
Proposition \ref{lebesguemodificado} states that $x$ is a Lebesgue
point of $\widehat T^h$. Then (\ref{termo2T}) is less than or
equal to the sum of terms like
\[
\int_M \mid \eta(x,y,\varepsilon)[\widehat T^h(y)-\widehat
T^h(x)]\mid dV_{\widetilde g}(y) \leq \frac{\Vert
\eta\Vert}{\varepsilon^n}\int_{B(x,\varepsilon)} \mid\widehat
T^h(y)-\widehat T^h(x)\mid dV_{\widetilde g}(y)
\]
that goes to zero as $\varepsilon$ goes to zero. Therefore
$\widehat T_\varepsilon\rightarrow \widehat T$ for almost every $x
\in B(x_0,r)$, what implies that $\widehat
T_\varepsilon\rightarrow \widehat T$ a.e..

\

(2) If $\widehat T$ is continuous, then $\widehat
T_\varepsilon\rightarrow \widehat T$ in $C^0_{\mathrm{loc}}(M)$ as
$\varepsilon \rightarrow 0$.

\

We will prove the convergence for tensors $T$ of type $(1,1)$. The
general case is analogous.

Fix an open set $U \subset \subset M$. Let $\{U_\omega \subset
\subset M,\phi_\omega\}_{\omega \in \Lambda}$ be a locally finite
covering of $M$ by almost Euclidean coordinate systems. Take a
partition of unity $\{\psi_\omega
:M\rightarrow\mathbb{R}\}_{\omega\in\Lambda}$ subordinated to
$\{U_\omega\}_{\omega\in\Lambda}$. Consider the functions
$\widehat T_{ij\omega_1 \omega_2}$ as in (\ref{Tijw}). Every $x\in
M$ is a Lebesgue point of $\widehat T$ what implies that
$(\widehat T_\varepsilon)_{ij\omega_1 \omega_2}(x) \rightarrow
\widehat T_{ij\omega_1 \omega_2}(x)$ as $\varepsilon \rightarrow
0$ for every $x \in U$. Then $(\widehat T_\varepsilon)_{ij\omega_1
\omega_2}$ converges uniformly to $\widehat T_{ij\omega_1
\omega_2}$ in $U$.

Set $\Lambda_U:=\{\omega \in \Lambda; U_\omega \cap \bar U \neq
\emptyset \}$. Observe that $\Lambda_U$ is a finite index set. Fix
$\epsilon > 0$ and let $r_U$ be the injectivity radius of $U$ in
$(\bigcup\limits_{\omega \in \Lambda_U} U_\omega, \widetilde g)$.
For every $(i,j,\omega_1,\omega_2) \in \{1,\ldots,n\}\times
\{1,\ldots,n\}\times \Lambda_U \times \Lambda_U$, we can choose a
family of positive numbers $\varepsilon_{ij\omega_1\omega_2} <
r_U$ such that $\Vert (\widehat
T_{\widetilde\varepsilon_{ij\omega_1\omega_2}})_{ij\omega_1\omega_2}-\widehat
T_{ij\omega_1\omega_2}\Vert_{C^0(U,\widetilde
g)}<\frac{\epsilon}{4n^2\vert\Lambda_U \vert^2}$ for every
$\widetilde\varepsilon_{ij\omega_1\omega_2}<\varepsilon_{ij\omega_1\omega_2}$.
If we choose
$\widetilde\varepsilon=\min\limits_{ij\omega_1\omega_2}\varepsilon_{ij\omega_1\omega_2}$
we have that $\Vert \widehat T_{\varepsilon}-\widehat
T\Vert_{C^0(U,\widetilde g)}<\epsilon$ for every
$\varepsilon<\widetilde\varepsilon$ (See (\ref{formulavarphi}),
(\ref{formulav}) and inequalities $\vert
(a_i(x,\varphi))_{\omega_1}\vert\leq 2$, $\vert
(b_j(x,v))_{\omega_2}\vert\leq 2$, $\Vert
(\vartheta_i)_{\omega_1}\Vert_{C^0}\leq 1$ and $\Vert
(u_j)_{\omega_2}\Vert_{C^0}\leq 1$). Therefore $\widehat
T_\varepsilon \rightarrow \widehat T$ in $C^0_{\mathrm{loc}}(M)$
as $\varepsilon\rightarrow 0$.

\

(3) If $1\leq p<\infty$ and $\widehat T\in L^p_{\mathrm{loc}}(M)$,
then $\widehat T_\varepsilon\rightarrow \widehat T$ in
$L^p_{\mathrm{loc}}(M)$ as $\varepsilon \rightarrow 0$.

\

Fix an open set $U \subset \subset M$ and let $\{U_\omega \subset
\subset M, \phi_\omega\}_{\omega \in \Lambda}$ be a locally finite
covering of $M$ by almost Euclidean coordinate systems. Denote
$\Lambda_U = \{\omega \in \Lambda; U_\omega \cap \bar U\neq
\emptyset\}$. Let $U_1 = \bigcup\limits_{\omega \in \Lambda_U}
U_\omega$ and denote $\Lambda_{U_1} = \{\omega \in \Lambda;
U_\omega \cap \bar U_1 \neq \emptyset\}$.

\

{\it Claim 1:} $\Vert \widehat
T_\varepsilon\Vert_{L^p(U,\widetilde g)}\leq C.\Vert \widehat
T\Vert_{L^p(U_1,\widetilde g)}$ for some positive constant $C$
that does not depend on $\widehat T$ and $\varepsilon <
\mathrm{inj}(U,(U_1,\widetilde g))$.

\

Let $\{\varphi_1,\ldots,\varphi_m\}$ and $\{v_1,\ldots,v_s\}$ be
respectively a family of smooth 1-forms and a family of smooth
vector fields on $M$ such that their $\Vert
\cdot\Vert_{C^0(M,\widetilde g)}$ norm are less than or equal to
1. Let $1<p<\infty$ (the case $p=1$ is easier) and fix $x\in U$.
Considering $\varepsilon$ less than the injectivity radius of $U$
in $(U_1,\widetilde g)$, we have that
\[
\begin{array}{l} \left\vert \widehat
T_\varepsilon(\varphi_1,\ldots,v_s)(x)\right\vert =\left\vert
\int_{U_1}\eta(x,y,\varepsilon).\widehat T(\varphi_1(x;\widetilde
g),\ldots, v_s(x;\widetilde g))(y).dV_{\widetilde g}(y)\right\vert
\\
\leq \int_{U_1}[\eta(x,y,\varepsilon)]^{1-\frac{1}{p}}.
[\eta(x,y,\varepsilon)]^{\frac{1}{p}}. \vert \widehat
T(\varphi_1(x;\widetilde g),\ldots, v_s(x;\widetilde
g))(y)\vert.dV_{\widetilde g}(y) \\
\leq \left(\int_{U_1}\eta(x,y,\varepsilon)\vert \widehat
T(\varphi_1(x;\widetilde g),\ldots, v_s(x;\widetilde
g))(y)\vert^p.dV_{\widetilde g}(y).\right)^{\frac{1}{p}},
\end{array}
\]
where we used the H\"older inequality in the last step. Notice
that this inequality is trivially true when $p=1$. Integrating in
$x$ we have that
\[
\Vert \widehat
T_\varepsilon(\varphi_1,\ldots,v_s)\Vert_{L^p(U,\widetilde
g)}^p=\int_U \vert \widehat
T_\varepsilon(\varphi_1,\ldots,v_s)(x)\vert^p.dV_{\widetilde g}(x)
\]
\[
\leq \int_U \left(\int_{U_1}\eta(x,y,\varepsilon)\vert \widehat
T(\varphi_1(x;\widetilde g), \ldots,v_s(x;\widetilde
g))(y)\vert^p.dV_{\widetilde g}(y).\right) dV_{\widetilde g}(x)
\]
\[
\leq \int_U \left(\int_{U_1}\eta(x,y,\varepsilon)\vert \bar
T(y)\vert^p.dV_{\widetilde g}(y).\right) dV_{\widetilde g}(x) \leq
(*)
\]
where $\bar T$ is defined by (\ref{equacaobarT}).

Using the Fubini's Theorem we have that
\[
(*) \leq \int_{U_1} \left(\int_{U}\eta(x,y,\varepsilon)\vert \bar
T(y)\vert^p.dV_{\widetilde g}(x).\right) dV_{\widetilde g}(y)
\]
\[
\leq \int_{U_1} \left(\int_{B(y,\varepsilon) \cap
U}\eta(x,y,\varepsilon)\vert \bar T(y)\vert^p.dV_{\widetilde
g}(x).\right) dV_{\widetilde g}(y)
\]
\[
\leq
\frac{\Vert\eta\Vert}{\varepsilon^n}\int_{U_1}\left(\int_{B(y,\varepsilon)
\cap U} \vert \bar T(y)\vert^p dV_{\widetilde
g}(x)\right).dV_{\widetilde g}(y)
\]
\[
\leq \frac{\Vert\eta\Vert}{\varepsilon^n} \mathrm{Vol\:}(B_{(\min
K)}(\varepsilon))\int_{U_1}\vert \bar T(y)\vert^p .dV_{\widetilde
g}(y) \leq C^p.\Vert \bar T\Vert^p_{L^p(U_1,\widetilde g)}
\]
where $(\min K)$ is the minimum sectional curvature on
$(U_1,\widetilde g)$, $\mathrm{Vol}(B_{(\min K)}(\varepsilon))$ is
the volume of the geodesic ball with radius $\varepsilon$ and
constant curvature $(\min K)$ and $C$ is a constant that does not
depend on $\varepsilon<\mathrm{inj}(U,(U_1,\widetilde g))$. Then
\[
\Vert \widehat
T_\varepsilon(\varphi_1,\ldots,v_s)\Vert_{L^p(U,\widetilde g)}\leq
C.\Vert \bar T\Vert_{L^p(U_1,\widetilde g)}
\]
\[
\leq \left\Vert
C.2^{m+s}\sum_{(\omega_1,\ldots,\omega_{m+s})}\sum_{(i_1,\ldots,i_m)}
\sum_{(j_1,\ldots,j_s)} \vert \widehat T_{i_1\ldots i_m j_1\ldots
j_s\omega_1\ldots \omega_{m+s}}\vert
\right\Vert_{L^p(U_1,\widetilde g)}
\]
\[
\leq
C.2^{m+s}\sum_{(\omega_1,\ldots,\omega_{m+s})}\sum_{(i_1,\ldots,i_m)}
\sum_{(j_1,\ldots,j_s)}\Vert \widehat T_{i_1\ldots i_m j_1\ldots
j_s\omega_1\ldots \omega_{m+s}} \Vert_{L^p(U_1,\widetilde g)}
\]
\[
\leq C.2^{m+s}n^{m+s}\vert\Lambda_{U_1}\vert^{m+s}\left\Vert
\widehat T \right\Vert_{L^p(U_1,\widetilde g)}
\]
which settles Claim 1.

\

{\it Claim 2:} $\lim\limits_{\varepsilon\rightarrow 0}\left\Vert
\widehat T-\widehat T_\varepsilon\right\Vert_{L^p(U,\widetilde
g)}=0$.

\

We will prove Claim 2 for tensor fields of type $(1,1)$. The
general case is analogous.

For every $i,j \in \{1, \ldots, n\}$, $\omega_1,\omega_2 \in
\Lambda_{U_1}$ and $\epsilon >0$, there exist a continuous
function $Q_{ij\omega_1\omega_2}:M\rightarrow\mathbb{R}$ such that
\[
\left\Vert \widehat
T_{ij\omega_1\omega_2}-Q_{ij\omega_1\omega_2}\right\Vert_{L^p(M,\widetilde
g)}<\frac{\epsilon}{4\vert \Lambda_{U_1}\vert^2 n^2}.
\]

Define
\[
Q(\varphi,v):=\sum_{i,j=1}^n\left(\sum_{(\omega_1,\omega_2)\in
\Lambda_{U_1}\times\Lambda_{U_1}}
(a_i)_{\omega_1}(b_j)_{\omega_2}Q_{ij\omega_1\omega_2}\right),
\]
where $(a_i)_{\omega_1}$ and $(b_j)_{\omega_2}$ are defined by
Eqs. (\ref{defineai}) and (\ref{definebj}) respectively. It is
easy to see that $Q$ is a continuous tensor field, and it is
straightforward that
\[
\left\Vert \widehat T-Q\right\Vert_{L^p(U_1,\widetilde
g)}<\epsilon.
\]

Moreover we have that $\Vert \widehat
T_\varepsilon-Q_\varepsilon\Vert_{L^p(U,\widetilde g)}<C\epsilon$
when $\varepsilon$ is small enough (See Claim 1), what implies
\[
\left\Vert \widehat T-\widehat
T_\varepsilon\right\Vert_{L^p(U,\widetilde g)}\leq \left\Vert
\widehat T-Q\right\Vert_{L^p(U,\widetilde g)}+\left\Vert
Q-Q_\varepsilon \right\Vert_{L^p(U,\widetilde g)}+\left\Vert
Q_\varepsilon -\widehat T_\varepsilon\right\Vert_{L^p(U,\widetilde
g)}
\]
\[
\leq \epsilon + \left\Vert Q-Q_\varepsilon
\right\Vert_{L^p(U,\widetilde g)} + C.\epsilon.
\]
Finally making $\varepsilon$ even smaller, we have that
$\left\Vert Q-Q_\varepsilon \right\Vert_{L^p(U,\widetilde
g)}<\epsilon$ because $Q_\varepsilon \rightarrow Q$ uniformly on
$U$. Therefore we have that
\[
\lim_{\varepsilon\rightarrow 0}\left\Vert \widehat T-\widehat
T_\varepsilon\right\Vert_{L^p(U,\widetilde g)}=0,
\]
what settles the theorem.

\end{proof}

\section{The mollifier smoothing with respect to $\mathcal P$}
\label{suavizacaosubordinadaparticao}

Let $M$ be a differentiable manifold and consider a non-regular
tensor field $\widehat T$ defined on $M$. In this section we
define a mollifier smoothing $\widehat T_{\varepsilon,\mathcal P}$
of $\widehat T$ that is quite adequate for our purposes: If
$\widehat T$ is a Riemannian metric of class $C^2$, then the
Levi-Civita connection and the Riemannian curvature tensor of
$\widehat T_{\varepsilon,\mathcal P}$ converges to the Levi-Civita
connection and to the Riemannian curvature tensor of $\widehat T$
respectively as $\varepsilon$ converges to zero. The idea is to
take a locally finite covering $M$ by open sets, each of them
endowed with an Euclidean background metric. We make the mollifier
smoothing given by Eq. (\ref{formulasuavizacao}) on each open set
of the covering and sum them weighted by a partition of unity
subordinated to the covering. Let us formalize the idea.

Let $M$ be a differentiable manifold. Let $\{(O_\omega,\widetilde
e_\omega)\}_{\omega \in \Lambda}$ and $\{U_\omega \subset \subset
O_\omega \}_{\omega \in \Lambda}$ be families of open subsets of
$M$ indexed by $\omega \in \Lambda$ such that

\begin{enumerate}

\item $O_\omega$ is endowed with the Euclidean background metric
$\widetilde e_\omega$ for each $\omega \in \Lambda$.

\item $\{U_\omega\}_{\omega \in \Lambda}$ is a locally finite
covering of $M$.

\item $\{x\in O_\omega; \mathrm{dist}_{\widetilde e_\omega}
(x,U_\omega)<1\} \subset \subset O_\omega$.

\end{enumerate}
Denote this family by $\{(U_\omega \subset \subset
O_\omega,\widetilde e_\omega)\}_{\omega \in \Lambda}$. It is not
difficult to see that such a family always exist.

Consider a non-regular tensor field $\widehat T\in T^{m,s}M$. The
tensor field $\widehat T_{\omega\varepsilon} \in
C^\infty(U_\omega)$ is defined as mollifier smoothing of $\widehat
T \vert_{O_\omega}$ in $U_\omega$ with respect to $\widetilde
e_\omega$, that is

\[
\begin{array}{c}
\widehat T_{\omega\varepsilon}(\varphi_1,\ldots,\varphi_m,v_1,\ldots,v_s)(x) \\
=\int\limits_{B(x,\varepsilon)}\eta\left(x,y,\varepsilon\right).
\widehat T(\varphi_1(x;\widetilde
e_\omega),\ldots,v_s(x;\widetilde e_\omega))(y).dV_{\widetilde
e_\omega}(y),
\end{array}
\]
where $x \in U_\omega$ and $\varepsilon<1$. Now let us sum the
tensor fields $\widehat T_{\omega \varepsilon}$ using a partition
of unity.

Let $\{\psi_\omega\}_{\omega \in \Lambda}$ be a partition of unity
subordinated to $\{U_\omega\}_{\omega \in \Lambda}$. In order to
simplify the notation, we denote the locally finite covering
$\{(U_\omega \subset \subset O_\omega,\widetilde
e_\omega)\}_{\omega \in \Lambda}$ together with the partition of
unity $\{\psi_\omega\}_{\omega \in \Lambda}$ by $\mathcal P$.

\begin{definition}
\label{suavizacaoemrelacaoaP} Let $\widehat T$ be a non-regular
tensor field of type $(m,s)$ on a differentiable manifold $M$. Let
$\mathcal P$ be a locally finite covering $\{(U_\omega \subset
\subset O_\omega, \widetilde g)\}_{\omega \in \Lambda}$ of $M$
together with a partition of unity $\{\psi_\omega\}_{\omega \in
\Lambda}$ subordinated to $\{(U_\omega \subset \subset O_\omega,
\widetilde g)\}_{\omega \in \Lambda}$. The {\em mollifier
smoothing of $\widehat T\in T^{m,s}M$ with respect to $\mathcal
P$} is a smooth tensor field $\widehat T_{\varepsilon,\mathcal P}$
on $M$ defined by
\[
\widehat T_{\varepsilon,\mathcal P}(x)=\sum_{\omega \in \Lambda}
\psi_\omega(x) \widehat T_{\omega\varepsilon}(x),\;\;\; x\in M,
\]
where $\varepsilon < 1$. It will be denoted simply by $\widehat
T_\varepsilon$ if there is not any possibility of
misunderstandings.

\end{definition}

The mollifier smoothing with respect to $\mathcal P$ has also nice
convergence properties.

\begin{theorem}
\label{convergenciamollifiersmoothingpint} Let $\widehat T$ be a
locally summable tensor field of type $(m,s)$ on a differentiable
manifold $M$. Let $\mathcal P$ be a locally finite covering
$\{(U_\omega \subset \subset O_\omega, \widetilde g)\}_{\omega \in
\Lambda}$ of $M$ together with a partition of unity
$\{\psi_\omega\}_{\omega \in \Lambda}$ subordinated to
$\{(U_\omega \subset \subset O_\omega, \widetilde g)\}_{\omega \in
\Lambda}$. Then we have the following convergences:
\begin{enumerate}

\item $\widehat T_{\varepsilon, \mathcal P} \rightarrow \widehat
T$ a.e. as $\varepsilon\rightarrow 0$.

\item If $\widehat T$ is continuous, then $\widehat
T_{\varepsilon, \mathcal P}\rightarrow \widehat T$ in
$C^0_{\mathrm{loc}}(M)$ as $\varepsilon \rightarrow 0$.

\item If $1\leq p<\infty$ and $\widehat T\in
L^p_{\mathrm{loc}}(M)$, then $\widehat T_{\varepsilon, \mathcal
P}\rightarrow \widehat T$ in $L^p_{\mathrm{loc}}(M)$ as
$\varepsilon \rightarrow 0$.

\end{enumerate}

\end{theorem}

\begin{proof}
It is an immediate consequence of Theorem \ref{convergencias} and
the fact that each open subset $U \subset \subset M$ intercepts
only a finite number of elements of $\{U_\omega\}_{\omega \in
\Lambda}$.

\end{proof}

Now we prove that the mollifier smoothing with respect to
$\mathcal P$ behaves well in the smooth case, that is, if
$\widehat g$ is of class $C^2$, then the first derivatives and the
second derivatives of $\widehat g_{\varepsilon,\mathcal P}$
converges to the correspondent derivatives of $\widehat g$ as
$\varepsilon$ goes to zero. In particular, the Levi-Civita
connection and the Riemannian curvature tensor of $(M,\widehat
g_{\varepsilon,\mathcal P})$ converges to the Levi-Civita
connection and the Riemannian curvature tensor of $(M,\widehat g)$
as $\varepsilon$ goes to zero.

\begin{remark}
\label{conexaonaoetensor} The Levi-Civita connection is not a
tensor field. Then we can not apply fully the theory of
non-regular tensor fields developed in the former sections. But
some results about the convergence of Levi-Civita connections can
be obtained. For instance we can think about the pointwise
convergence of Levi-Civita connections as will be done in Theorem
\ref{suavizacao2boa}. Another application will be given in Section
\ref{transporteparalelonaoregular}, where the parallel transport
of a vector through a curve will be generalized for some
non-regular Riemannian manifolds.
\end{remark}

\begin{theorem}
\label{suavizacao2boa} Let $\widehat g$ be a $C^2$ Riemannian
metric defined on $M$ and fix $\mathcal P$ on $M$. Then the
Levi-Civita connection and the Riemannian curvature tensor with
respect to the metric $\widehat g_{\varepsilon,\mathcal P}$
converges everywhere to the the Levi-Civita connection and the
Riemannian curvature tensor of $(M,\widehat g)$ respectively as
$\varepsilon \rightarrow 0$.
\end{theorem}

\begin{proof}
Let $\mathcal P$ denote the locally finite covering $\{(U_\omega
\subset \subset O_\omega,\widetilde e_\omega)\}_{\omega \in
\Lambda}$ with a partition of unity $\{\psi_\omega\}_{\omega \in
\Lambda}$ subordinated to $\{U_\omega\}_{\omega \in \Lambda}$. We
begin studying the convergence of $\widehat g_{\omega
\varepsilon}$ to $\widehat g_{\omega}:=\widehat g\vert_{U_\omega}$
at $x \in M$ up to derivatives of order two as $\varepsilon$ goes
to zero.

Fix $\omega \in \Lambda$, $x \in U_\omega$ and let
$(x_1,\ldots,x_n)$ be a coordinate system in a neighborhood
$N_\omega \subset U_\omega$ of $x$ such that $(\widetilde
e_\omega)_{ij} = \delta_{ij}$ in this coordinate system.

Denote the representation of $\widehat g$ and $\widehat
g_{\omega\varepsilon}$ in this coordinate system by $(\widehat
g_\omega)_{ij}$ and $(\widehat g_{\omega\varepsilon})_{ij}$
respectively. Observe that
\[
(\widehat g_{\omega\varepsilon})_{ij}(x) = \widehat
g_{\omega\varepsilon} \left(\frac{\partial}{\partial x_i},
\frac{\partial}{\partial x_j}\right)(x)
\]
\[
=\int_M \eta(x,y,\varepsilon)\widehat
g_\omega\left(\frac{\partial}{\partial
x_i},\frac{\partial}{\partial x_j}\right)(y) dV_{\widetilde
e_\omega}(y)=\int_M \eta(x,y,\varepsilon)(\widehat
g_\omega)_{ij}(y) dV_{\widetilde e_\omega}(y).
\]
Therefore all the first and second order derivatives of $(\widehat
g_{\omega\varepsilon})_{ij}$ converges to the correspondent
derivatives of $(\widehat g_\omega)_{ij}$ because $N_{\omega}$ can
be considered as an open set in $\mathbb R^n$ (see \cite{7}).

Let $(z_1,\ldots,z_n)$ be another coordinate system in $N_\omega$.
The metric $\widehat g_{\omega \varepsilon}$ in this coordinate
system is given by
\begin{equation}
\label{mudancadecoordenadas} (\widehat g_{\omega
\varepsilon})_{\hat i\hat j}(x)= \sum_{i,j=1}^n(\widehat g_{\omega
\varepsilon})_{ij}(x) \frac{\partial x_i}{\partial z_{\hat
i}}\frac{\partial x_j}{\partial z_{\hat j}}= \sum_{i,j=1}^n
\frac{\partial x_i}{\partial z_{\hat i}}\frac{\partial
x_j}{\partial z_{\hat j}}\int_M \eta(x,y,\varepsilon)(\widehat
g_\omega)_{ij}(y) dV_{\widetilde e_\omega}(y)
\end{equation}
that converges to
\[
\sum_{i,j=1}^n \frac{\partial x_i}{\partial z_{\hat
i}}\frac{\partial x_j}{\partial z_{\hat j}}(\widehat
g_\omega)_{ij}(x) = (\widehat g_{\omega})_{\hat i\hat j}(x)
\]
for every $\hat i,\hat j \in \{1,\ldots,n\}$ as $\varepsilon$ goes
to $0$.

Now let us study the convergence of the derivatives of $(\widehat
g_{\omega \varepsilon})_{\hat i\hat j}$.

The coordinate vector fields satisfies the relationship
\[
\frac{\partial}{\partial z_{\hat k}}= \sum_{s=1}^n \frac{\partial
x_s}{\partial z_{\hat k}}\frac{\partial}{\partial x_s}.
\]
Calculating the derivative of (\ref{mudancadecoordenadas}) with
respect to $\partial/\partial z_k$, we have that
\[
\frac{\partial}{\partial z_{\hat k}}(\widehat g_{\omega
\varepsilon})_{\hat i\hat j}(x)
\]
\[
=\sum_{i,j=1}^n \left[\frac{\partial}{\partial z_{\hat k}} \left(
\frac{\partial x_i}{\partial z_{\hat i}}\frac{\partial
x_j}{\partial z_{\hat j}}\right) \right] \int_M
\eta(x,y,\varepsilon)(\widehat g_\omega)_{ij}(y) dV_{\widetilde
e_\omega}(y)
\]
\[
+\sum_{i,j,s=1}^n \left( \frac{\partial x_i}{\partial z_{\hat
i}}\frac{\partial x_j}{\partial z_{\hat j}}  \frac{\partial
x_s}{\partial z_{\hat k}}\right) \left(\frac{\partial}{\partial
x_s}\int_M \eta(x,y,\varepsilon)(\widehat g_\omega)_{ij}(y)
dV_{\widetilde e_\omega}(y) \right)
\]
\[
=\sum_{i,j=1}^n \left[\frac{\partial}{\partial z_{\hat k}}\left(
\frac{\partial x_i}{\partial z_{\hat i}}\frac{\partial
x_j}{\partial z_{\hat j}}\right) \right] \int_M
\eta(x,y,\varepsilon)(\widehat g_\omega)_{ij}(y) dV_{\widetilde
e_\omega}(y)
\]
\[
+\sum_{i,j,s=1}^n \left( \frac{\partial x_i}{\partial z_{\hat
i}}\frac{\partial x_j}{\partial z_{\hat j}} \frac{\partial
x_s}{\partial z_{\hat k}}\right) \left( \int_M
\eta(x,y,\varepsilon)\frac{\partial}{\partial x_s}(\widehat
g_\omega)_{ij}(y) dV_{\widetilde e_\omega}(y) \right)
\]
that converges to
\[
\sum_{i,j=1}^n \left[\frac{\partial}{\partial z_{\hat k}}\left(
\frac{\partial x_i}{\partial z_{\hat i}}\frac{\partial
x_j}{\partial z_{\hat j}}\right) \right] (\widehat
g_\omega)_{ij}(x) + \sum_{i,j,s=1}^n \left( \frac{\partial
x_i}{\partial z_{\hat i}}\frac{\partial x_j}{\partial z_{\hat j}}
\frac{\partial x_s}{\partial z_{\hat k}}\right)
\frac{\partial}{\partial x_s}(\widehat g_\omega)_{ij}(x)
\]
\[
=\frac{\partial}{\partial z_{\hat k}}(\widehat g_{\omega })_{\hat
i\hat j}(x)
\]
as $\varepsilon$ goes to zero.

Analogously we have that
\[
\frac{\partial}{\partial z_{\hat l}}\frac{\partial}{\partial
z_{\hat k}}(\widehat g_{\omega\varepsilon})_{\hat i\hat j}(x)
\mathrm{\;\;converges\;\;to\;\;}\frac{\partial}{\partial z_{\hat
l}}\frac{\partial}{\partial z_{\hat k}}(\widehat g_\omega)_{\hat
i\hat j}(x)
\]
for every $\hat i, \hat j, \hat l, \hat k\in\{1,\ldots,n\}$ as
$\varepsilon$ goes to zero.

Now let us study the convergence of $\widehat
g_{\varepsilon,\mathcal P}$ to $\widehat g$ at $x \in M$ up to
derivatives of order two. The mollifier smoothing of $\widehat g$
with respect to $\mathcal P$ is given by the locally finite sum
\[
\widehat g_{\varepsilon,\mathcal P}=\sum_\omega \psi_\omega
\widehat g_{\omega\varepsilon}.
\]
Take a coordinate system $(z_1,\ldots,z_n)$ in a neighborhood
$N_x$ of $x$. We have that
\[
(\widehat g_{\varepsilon,\mathcal P})_{\hat i\hat j}=\sum_\omega
\psi_\omega(x). (\widehat g_{\omega\varepsilon})_{\hat i\hat j}
\]
in $N_x$, and using the calculations made before, it is not
difficult to see that
\[
(\widehat g_{\varepsilon,\mathcal P})_{\hat i\hat j}(x)
\mathrm{\;\;converges\;\;to\;\;}\widehat g_{\hat i\hat j}(x),
\]
\[
\frac{\partial}{\partial z_{\hat k}}(\widehat
g_{\varepsilon,\mathcal P})_{\hat i\hat j}(x)
\mathrm{\;\;converges\;\;to\;\;}\frac{\partial}{\partial z_{\hat
k}}\widehat g_{\hat i\hat j}(x)
\]
and that
\[
\frac{\partial}{\partial z_{\hat l}}\frac{\partial}{\partial
z_{\hat k}}(\widehat g_{\varepsilon,\mathcal P})_{\hat i\hat j}(x)
\mathrm{\;\;converges\;\;to\;\;}\frac{\partial}{\partial z_{\hat
l}}\frac{\partial}{\partial z_{\hat k}}\widehat g_{\hat i\hat
j}(x)
\]
as $\varepsilon$ goes to zero.

Finally notice that the Christoffel symbols and the components
Riemannian curvature tensor can be written in terms of the metric
and its derivatives up to order two (See Eqs.
(\ref{levicivitacoordenadas}) and (\ref{curvaturacoordenadas})).
Therefore the Christoffel symbols of $\widehat
g_{\varepsilon,\mathcal P}$ converges to the Christoffel symbols
of $\widehat g$ and the components of the Riemannian curvature
tensor of $\widehat g_{\varepsilon,\mathcal P}$ converges to the
components of the Riemannian curvature tensor of $\widehat g$, as
we wanted to prove.

\end{proof}

We also have the following version of Theorem \ref{suavizacao2boa}
for the $\kappa$-th Lipschitz-Killing curvature measures.

\begin{corollary}
\label{versaoconvergenciaformasdecurvatura} Let $\widehat g$ be a
$C^2$ Riemannian metric defined on $M$ and fix $\mathcal P$ on
$M$. Then the $\kappa$-th Lipschitz Killing curvature measures of
$(M,\widehat g_{\varepsilon,\mathcal P})$ converges everywhere
 to the $\kappa$-th Lipschitz Killing
curvature measures of $(M,\widehat g)$.
\end{corollary}

\begin{proof}
For every $\varepsilon > 0$, choose an orthonormal moving frame
$\{w_{\varepsilon 1}, \ldots, w_{\varepsilon n}\}$ for $(M,
\widehat g_{\varepsilon,\mathcal P})$ such that $\{w_{\varepsilon
1}, \ldots,$ $w_{\varepsilon n}\}$ converges everywhere to an
orthonormal moving frame $\{w_1,\ldots,$ $w_n\}$ of $(M,\widehat
g)$ as $\varepsilon$ goes to zero. Consider the curvature forms of
$(M,\widehat g_{\varepsilon,\mathcal P})$ and $(M,\widehat g)$
with respect to these orthonormal frames. Then it follows from
Theorem \ref{suavizacao2boa} and (\ref{formasdecurvatura}),
(\ref{curvaturaintegralipschitzkilling}) and
(\ref{lipschitzkillingform}) that the curvature forms of
$(M,\widehat g_{\varepsilon,\mathcal P})$ converges everywhere to
the respective curvature forms of $(M,\widehat g)$, what settles
the corollary.

\end{proof}

Therefore if we want to generalize objects of the classical
Riemannian geometry to non-regular Riemannian manifolds we can
study how these objects behave when $(M,\widehat
g_{\varepsilon,\mathcal P})$ converge to $(M,\widehat g)$ as
$\varepsilon$ goes to zero. The following sections give
applications of this type.

\section{The definition of distance for non-regular Riemannian
manifolds} \label{distancianaoregular}

Let $M$ be a differentiable manifold. Two $L^p_{\mathrm{loc}}$
Riemannian metrics that differ on a subset of measure zero are
identical. Then the distance between $x$ and $y$ in $(M, \widehat
g)$ can not be defined in the classical fashion, that is, as the
infimum of the lengths of all piecewise regular curves that
connects $x$ and $y$. We should look for a more stable definition.

In this section we define the distance between two points in
$(M,\widehat g)$ using smooth approximations. But let us see some
typical examples first. The mollifier smoothing with respect to a
background metric $\widetilde g$ and the mollifier smoothing with
respect to $\mathcal P$ provide smooth approximations of a
non-regular Riemannian metric $\widehat g$. We use the former
mollifier smoothing in order to see what happens if we try to
define the distance between $x$ and $y$ in $(M,\widehat g)$ as
``$d_{\widehat g}(x,y)=\lim\limits_{\varepsilon \rightarrow
0}d_{\widehat g_\varepsilon}(x,y)$''.

\begin{example}
\label{metricanaopositiva} Let $(M,\widehat g)$ be a
two-dimensional non-regular Riemannian manifold. Let $U\subset M$
be an open set, $\phi: U\rightarrow
(-1,1)\times(-1,1)\subset\mathbb R^2$ a coordinate system and
define a background Riemannian metric $\widetilde g$ on $M$ such
that its restriction to $U$ is the canonical Euclidean metric with
respect to $\phi$. Suppose that the metric $\widehat g$ in this
coordinate system is given by
\begin{equation}
\label{exemplometrica1} \widehat g= \left\{
\begin{array}{ccc}
dx^2+x^2.dy^2 & \mathrm{if} & x\not = 0 \\
dx^2+dy^2 & \mathrm{if} & x=0.
\end{array}
\right.
\end{equation}
Observe that $\widehat g\vert_{U}$ is positive definite, but
\[
\lim_{\varepsilon \rightarrow 0} d_{\widehat g_\varepsilon}
((0,-1/2),(0,1/2))=0.
\]
We can see this fact representing $(U,\widehat g)$ as the domain

\[
\{(x,y)\in\mathbb R^2; x\in(-1,1),-|x|<y<|x|\} \cup \{(0,0)\}
\]
with the canonical Euclidean metric.

\end{example}

\begin{example}
\label{metricacomdistanciaoscilatoria} Let $U\subset M$,
$\widetilde g$ and $\phi$ as in Example \ref{metricanaopositiva}.
Suppose that the metric $\widehat g$ in the coordinate system
$\phi$ is given by
\begin{equation}
\label{exemplometrica2} \widehat g= \left\{
\begin{array}{ccc}
dx^2+dy^2 & \mathrm{if} & x=0;\\
dx^2+dy^2 & \mathrm{if} & \frac{1}{2^{2n+1}}<\vert x\vert
<\frac{1}{2^{2n}}
\;n\in\mathbb N\cup\{0\}; \\
2dx^2+2dy^2 & \mathrm{if} & \frac{1}{2^{2n+2}}<\vert x\vert
<\frac{1}{2^{2n+1}},\;n\in\mathbb N\cup\{0\}.
\end{array}
\right. \end{equation}

If we calculate the length of the straight line segment connecting
$(0,-1/2)$ and $(0,1/2)$ in $(M,\widehat g_\varepsilon)$, then
this length will oscillate as $\varepsilon$ goes to zero.
Therefore there is not any hope that the arc length of a curve
converges always to a limit as $\varepsilon$ goes to zero.
\end{example}

\begin{example}
\label{Metricainfinita} Let $U\subset M$, $\widetilde g$ and
$\phi$ as in Example \ref{metricanaopositiva}. Suppose that the
metric $\widehat g$ in this coordinate system is given by
\begin{equation}
\label{exemplometrica3} \widehat g= \left\{
\begin{array}{ccc}
\frac{1}{\sqrt{x^2+y^2}} (dx^2+dy^2) & \mathrm{if} & (x,y)\neq (0,0) \\
dx^2+dy^2 & \mathrm{if} & (x,y)=(0,0).
\end{array}
\right.
\end{equation}
Notice that $\widehat g\in L^p_{\mathrm{loc}}(M)$ for $p<2$. For
$P\neq (0,0)$, we have that $d_{\widehat g_\varepsilon}
((0,0),P)\rightarrow \infty$ as $\varepsilon$ goes to zero,
because $\widehat g_{\varepsilon}$ converges uniformly to
$\widehat g$ on compact subsets of $(U-(0,0))$ as $\varepsilon$
goes to zero. This example shows that the distance induced by an
$L^p_{\mathrm{loc}}(M)$ Riemannian metric can be infinity.
\end{example}

Now let us define the Riemannian distance. Let $(M,\widehat g)$ be
a non-regular Riemannian manifold. Denote by $\widehat
g(\varepsilon)$ a one-parameter family of smooth Riemannian
metrics parametrized by $\varepsilon > 0$ such that
$\lim\limits_{\varepsilon \rightarrow 0} \widehat
g(\varepsilon)=\widehat g$ in $L^p_{\mathrm{loc}}(M)$. Consider
\[
\mathcal G=\{\widehat g(\varepsilon); \lim_{\varepsilon
\rightarrow 0} \widehat g(\varepsilon)=\widehat g \mathrm{\;in\;}
L^p_{\mathrm{loc}}(M)\}
\]
the collection of all one-parameter families of smooth Riemannian
metrics param\-etrized by $\varepsilon >0$ that converges to
$\widehat g$ in $L^p_{\mathrm{loc}}(M)$.

\begin{definition}
\label{distancialp} Let $(M,\widehat g)$ be a non-regular
Riemannian manifold. The Riemannian distance is defined by
\[
d_{\widehat g}(x,y):=\sup\limits_{\widehat g(\varepsilon) \in
\mathcal G}[\limsup\limits_{\varepsilon \rightarrow 0} d_{\widehat
g(\varepsilon)}(x,y)].
\]
\end{definition}

\begin{remark}
\label{distanciac0} If $\widehat g \in C^0_{\mathrm{loc}}(M)$,
then our definition coincides with the classical definition of
distance.

\end{remark}

A non-regular Riemannian manifold $(M,\widehat g)$ with the
Riemannian distance $d_{\widehat g}$ is almost a metric space.
Eventually the distance between two points can be zero (See
Example \ref{metricanaopositiva}). But the sets of points with
zero distance are equivalence classes in $M$. The following
theorem states that $d_{\widehat g}$ induces a metric on the
identification space $M/\sim$.

\begin{theorem}
Let $(M,\widehat g)$ be a non-regular Riemannian manifold denote
the identification space of $M$ with the sets of points with zero
distance identified by $M/\sim$. Then $d_{\widehat g}$ induces a
metric $\widetilde d_{\widehat g}$ on $M/\sim$ by
\[
\widetilde d_{\widehat g}(\widetilde x,\widetilde y)= d_{\widehat
g}(x,y)
\]
where $x$ and $y$ are elements of the equivalence classes of
$\widetilde x$ and $\widetilde y$ respectively.
\end{theorem}

\begin{proof}
Immediate.

\end{proof}

Notice that $(M/\sim,\widetilde d_{\widehat g})$ can be a metric
space that is not homeomorphic to a differentiable manifold (See
Example \ref{metricanaopositiva}).

\section{The definition of Parallel Transport on non-regular Riemannian manifolds}
\label{transporteparalelonaoregular}

Let $(M,\widehat g)$ be a non-regular Riemannian manifold. Let $x$
and $y$ be two points in $M$ and consider a piecewise regular
curve $\gamma:[a,b]\rightarrow M$ connecting $x$ and $y$. Remember
that $\widehat g_{\varepsilon,\mathcal P}$ is the mollifier
smoothing of $\widehat g$ with respect to $\mathcal P$. We denote
by $\hat \tau^\gamma_{x,y,\varepsilon} : T^{m,s}_xM \rightarrow
T^{m,s}_yM$ the parallel transport through $\gamma$ with respect
to the metric $\widehat g_{\varepsilon,\mathcal P}$ (Here we do
not carry $\mathcal P$ in the notation of $\hat
\tau^\gamma_{x,y,\varepsilon}$ for the sake of simplicity).

Now we are ready to define the parallel transport in non-regular
Riemannian manifolds.

\begin{definition}
\label{definicaotransporteparalelo} Let $(M,\widehat g)$ be a
non-regular Riemannian manifold. Let $x$ and $y$ be two points in
$M$ and consider a piecewise regular curve
$\gamma:[a,b]\rightarrow M$ connecting $x$ and $y$. The parallel
transport $\hat \tau^\gamma_{x,y}:T^{m,s}_xM\rightarrow
T^{m,s}_yM$ through $\gamma$ is defined as
\[
\hat \tau^\gamma_{x,y}(T)=\lim_{\varepsilon\rightarrow 0}\hat
\tau^\gamma_{x,y,\varepsilon}(T)
\]
if the limit exists and it does not depend on $\mathcal P$.
\end{definition}

In this work, we will not look for general conditions that
guarantees the existence of the parallel transport. Instead, we
are going to see that the parallel transport can be defined
through some curves on {\em piecewise smooth two-dimensional
Riemannian manifolds} (See Definition
\ref{definicaosuperficiesuaveporpartes} ahead). Before defining
this kind of surface, we introduce some concepts and notations.

\begin{definition}
\label{triangulacao} Let $M$ be a compact differentiable
two-dimensional manifold (eventually with boundary). A {\em
triangulation} of $M$ is a homeomorphism $\Theta:\Upsilon
\rightarrow M$ from a simplicial complex $\Upsilon$ onto $M$.
\end{definition}

As usual, we call the image of each 0-dimensional simplex by
vertex, the image of each 1-dimensional simplex by edge and the
image of each 2-dimensional simplex by face (of the
triangulation). The set of vertices will be denoted by $V$, the
set of points which are in the interior of some edge will be
denoted by $E$ and the set of points which are in the interior of
some face will be denoted by $F$. Then $M$ can be written as the
disjoint union $V \cup E \cup F$.

\begin{definition}
\label{variedadesuaveporpartes} Let $M$ be a compact
differentiable two-dimensional manifold (eventually with boundary)
and $\Theta:\Upsilon \rightarrow M$ be a triangulation of $M$. We
say that $\Theta$ is a piecewise smooth triangulation if $\Theta$
restricted to every two simplex is a diffeomorphism onto its
image. In particular every edge is the image of a regular curve.
\end{definition}

We want give the definition of piecewise smooth Riemannian metric
on compact differentiable two-dimensional manifolds. The following
example indicates what should happen on the edge of a
triangulation.

\begin{example}
\label{superficiepoliedral} Let $\Pi_1=\{(x,y,z)\in\mathbb
R^3;y\leq 0,z=0\}$ and $\Pi_2=\{(x,y,z)\in\mathbb R^3;y\geq
0,z=y\}$. Then $\Pi=\Pi_1\cup\Pi_2$ is a smooth surface outside
the $x$-axis. If we induce the canonical metric of $\mathbb R^3$
on the regular part of $\Pi$ and put the coordinate system
$(x,y)\in\mathbb R^2$ on $\Pi$, then the metric is given by
\[
ds^2= \left\{
\begin{array}{lcl}
dx^2+dy^2 & \;\;\; & \mathrm{if\;}y<0 \\
dx^2+2dy^2 & \;\;\; & \mathrm{if\;}y>0.
\end{array}
\right.
\]
Observe that the metric is not defined on the $x$-axis. This fact
does not represent any problem because the edge is a set of
measure zero. Although the metric can not be extended to the
$x$-axis, it induces a metric $dx^2$ on it.
\end{example}

The behavior found in Example \ref{superficiepoliedral} is defined
as follow.

\begin{definition}
\label{colarbemnovertice} Let $M^n$ be a differentiable manifold
and let $M_1,M_2\subset M$ be two disjoint open sets. Suppose that
$M_3^{n-1}\subset\bar M_1\cap \bar M_2$ is a differentiable
$(n-1)$-dimensional submanifold of $M$, where $\bar M_1$ and $\bar
M_2$ denote the closure of $M_1$ and $M_2$ in $M$ respectively.
Let $g_1$ and $g_2$ be two smooth Riemannian metrics on $M_1$ and
$M_2$ respectively such that they induces smooth Riemannian
metrics on $M_3$. We say that $\bar M_1$ and $\bar M_2$ {\em glues
nicely through $M_3$} if the Riemannian metric induced on $M_3$ by
$g_1$ is equal to the Riemannian metric induced by $g_2$.
\end{definition}

We are now ready to give the definition of piecewise smooth
two-dimensional Riemannian manifold.

\begin{definition}
\label{definicaosuperficiesuaveporpartes} Let $(M,\widehat g)$ be
a compact two-dimensional non-regular Riemannian manifold
(eventually with boundary). We say that $\widehat g$ is a
piecewise smooth Riemannian metric if there exist a piecewise
smooth triangulation $\Theta:\Upsilon \rightarrow M$ such that

\begin{enumerate}

\item $\widehat g\vert_{\Theta(\mathrm{int}(\Xi))}$ is smooth for
every two simplex $\Xi \in C$ and $\widehat
g\vert_{\Theta(\mathrm{int}(\Xi))}$ is smoothly extendable to
$\Theta(\Xi)$;

\item If $\Xi_1$ and $\Xi_2$ are simplexes which have an edge in
common, then $(\mathrm{int}(\Xi_1),\widehat g)$ and
$(\mathrm{int}(\Xi_2),\widehat g)$ glues nicely through $E \cap
\Xi_1\cap\Xi_2$.

\end{enumerate}
In this case, we say that $\Theta$ is a {\em triangulation
associated to $(M,\widehat g)$}.

A compact differentiable two-dimensional manifold (eventually with
boundary) with a piecewise smooth Riemannian metric is called a
{\em piecewise smooth two-dimensional Riemannian manifold}.

\end{definition}

\begin{remark}
\label{observacaosobremetricapoliedral} Observe that $\widehat g$
does not need to satisfy any condition over $V$. Formally it does
not represent any problem because it is enough to define a
non-regular Riemannian metric outside a set of measure zero.
Geometrically, it is quite natural to do so because the vertices
of a polyhedral surface (eventually with non-flat faces and edges)
are quite wild compared to its regular part. What is important
here is that all the geometrical information of a vertex is found
on its neighborhood.
\end{remark}

Here we begin to study the parallel transport in a piecewise
smooth two-dimensional Riemannian manifold $(M,\widehat g)$. The
following theorem tell us that when $\gamma$ crosses an edge
transversally, then the parallel transport of a vector through
$\gamma$ keep its angle with the edge.
\begin{theorem}
\label{passagemporaresta} Let $(M,\widehat g)$ be a piecewise
smooth two-dimensional Riemannian manifold and fix $\mathcal P$.
Let $\nu:[a,b]\rightarrow E$ be a regular parametrization of a
piece of edge. Let $\gamma:[-\epsilon,\epsilon]\rightarrow M$ be a
regular curve such that it intersects $\nu([a,b])$ transversally
at the point $\gamma(0)=\nu(0)$. Suppose also that $(E\cup V)\cap
\gamma([-\epsilon,\epsilon])=\gamma(0)$. Let $(t,s)$ be a
coordinate system of a neighborhood of
$\gamma([-\epsilon,\epsilon])$ such that $\gamma(t)$ has
coordinates $(t,0)$ and $\nu(s)$ has coordinates $(0,s)$. Fix $v
\in T_{\gamma(-\epsilon)}M$. Then  the parallel transport $\hat
\tau^\gamma_{x,y}$ through $\gamma$ with respect to the metric
$\widehat g$ is well defined and
\[
\lim_{t\rightarrow 0^-}\mathrm{angle}
\left(\frac{\partial}{\partial s}(t,0),\hat
\tau^\gamma_{\gamma(-\epsilon),\gamma(t)}
(v)\right)=\lim_{t\rightarrow
0^+}\mathrm{angle}\left(\frac{\partial}{\partial s}(t,0), \hat
\tau^\gamma_{\gamma(-\epsilon),\gamma(t)} (v)\right).
\]
\end{theorem}

\begin{proof}
First of all we remark that $\varepsilon$ and $\epsilon$ are
different variables.

Let $v_\varepsilon(t)$ be the parallel transport of $v$ through
$\gamma$ with respect to the metric $\widehat
g_\varepsilon=\widehat g_{\varepsilon,\mathcal P}$. Let
$\cos_\varepsilon(t)$ be the cosine of the angle between
$v_\varepsilon(t)$ and $\frac{\partial}{\partial s}(t,0)$.

If $\widehat g_{\varepsilon ij}$ are the components of $\widehat
g_\varepsilon$ with respect to $(t,s)$ we have that
\[
\cos_\varepsilon(t) = \frac{v_{\varepsilon 1}(t).\widehat
g_{\varepsilon 12}(t) + v_{\varepsilon 2}(t) .\widehat
g_{\varepsilon 22}(t)}{ \sqrt{\widehat g_{\varepsilon 11}(t).
v_{\varepsilon 1}^2(t) + 2.\widehat g_{\varepsilon
12}(t).v_{\varepsilon 1}(t). v_{\varepsilon 2}(t) + \widehat
g_{\varepsilon 22}(t).v_{\varepsilon 2}^2(t)}. \sqrt{\widehat
g_{\varepsilon 22}(t)}}
\]
where $v_{\varepsilon 1}(t)$ and $v_{\varepsilon 2}(t)$ are the
components of $v_\varepsilon$ with respect to $(t,s)$. Now we use
Eq. (\ref{levicivitacoordenadas}) and
\[
\nabla_{\frac{\partial}{\partial t}} v_{i
\varepsilon}=-(\Gamma_{1i}^1)_\varepsilon v_{\varepsilon
1}(t)-(\Gamma_{1i}^2)_\varepsilon v_{\varepsilon 2}(t)
\]
in order to get
\begin{equation}
\label{variacaocost} \frac{\partial}{\partial t}
cos_\varepsilon(t)=\frac{1}{2}\frac{v_{\varepsilon 1}(t)\left(
\frac{\partial}{\partial s } \widehat g_{\varepsilon 11}(t).
\widehat g_{\varepsilon 22}(t) - \frac{\partial}{\partial t}
\widehat g_{\varepsilon 22}(t). \widehat g_{\varepsilon
12}(t)\right)} {\sqrt{\widehat g_{\varepsilon 11}(t) .
v_{\varepsilon 1}^2(t) + 2 . \widehat g_{\varepsilon 12}(t) .
v_{\varepsilon 1}(t) . v_{\varepsilon 2}(t)+ \widehat
g_{\varepsilon 22}(t). v_{\varepsilon 2}^2(t)}.\left(\widehat
g_{\varepsilon 22}(t)\right)^{\frac{3}{2}}}.
\end{equation}
The term $\frac{\partial}{\partial s} \widehat g_{\varepsilon
11}(t)$ is uniformly bounded with respect to $\varepsilon$ because
the derivation is done parallel to the edge. The term
$\frac{\partial}{\partial t} \widehat g_{\varepsilon 22}(t)$ is
uniformly bounded in terms of $\varepsilon$ because the faces glue
nicely through the edge and $\widehat g_{22}$ is a Lipschitz
function. Then the total variation of $\cos_\varepsilon(t)$ goes
to zero as $\epsilon$ goes to zero. This settles the theorem.

\end{proof}

\begin{theorem}
\label{existenciatransporteparalelo} Let $(M, \widehat g)$ be a
piecewise smooth two-dimensional Riemannian manifold. Let
$\gamma:[a,b] \rightarrow M$ be a piecewise regular curve such
that it does not intercept $V$ and it intercepts $E$
transversally. Then the parallel transport $\hat
\tau^\gamma_{x,y}:T^{m,s}_{\gamma(a)}M \rightarrow
T^{m,s}_{\gamma(b)}M$ through $\gamma$ is well defined.
\end{theorem}

\begin{proof}
The parallel transport of a scalar is trivial. The parallel
transport of a vector through $\gamma\vert_F$ is well defined due
to Theorem \ref{suavizacao2boa} and the parallel transport in a
neighborhood of $E$ is well defined due to Theorem
\ref{passagemporaresta}. The parallel transport of 1-forms through
$\gamma$ is well defined: In fact, given a $\varphi \in
T^*_{\gamma(a)}M$, there exist a unique field of $1$-forms through
$\gamma$ such that its contraction with parallel vector fields is
constant through $\gamma$. It is not difficult to see that this
field of $1$-forms are the parallel transport of $\varphi$ through
$\gamma$. Finally the parallel transport of other types of tensors
through $\gamma$ is also well defined: Analogously to the case of
$1$-forms, it is not difficult to prove that the parallel
transport of general tensors coincide with sums of tensor products
of parallel vector fields and parallel 1-forms.
\end{proof}

\section{The Lipschitz-Killing curvature measure for closed non-regular
Riemannian manifolds}

\label{curvaturalknaoregular}

In this section, we give another application of the mollifier
smoothing with respect to $\mathcal P$. We generalize the
$\kappa$-th Lipschitz-Killing curvature measure for closed
non-regular Riemannian manifolds using the mollifier smoothing
$\widehat g_{\varepsilon,\mathcal P}$. It will be a signed measure
$\mathcal R^\kappa(\cdot,\widehat g):\mathcal {\widetilde
B}(M)\rightarrow \mathbb R$, where $\mathcal {\widetilde B}(M)$ is
a $\sigma$-algebra that contains all Borel sets of $M$.

Denote the $\kappa$-th Lipschitz-Killing curvature measure
(defined on the open subsets of $M$) of $(M,\widehat
g_{\varepsilon,\mathcal P})$ by $\mathcal R^\kappa(\cdot,\widehat
g_{\varepsilon,\mathcal P})$. If we want to extend the the
Lipschitz-Killing curvature measure $\mathcal
R^\kappa(\cdot,\widehat g)$ to a closed non-regular Riemannian
manifold $(M,\widehat g)$, the first natural trial could be:

\

{\rm \noindent {\sc ``Definition''.} The {\em Lipschitz-Killing
curvature measure} of an open set $O$ with respect to $\widehat g$
is defined by
\[
\mathcal R^\kappa(O,\widehat g):=\lim_{\varepsilon\rightarrow
0}\mathcal R^\kappa (O,\widehat g_{\varepsilon, \mathcal P})
\]
if this limit exists and it does not depend on $\mathcal P$. }

\

Unfortunately this approach does not work, even for the example
given in the introduction where a two dimensional sphere converges
to a cube. In fact, consider the open set $M-\{p\}$ where $p$ is a
vertex of the cube. Remember that $\mathcal R^2(\cdot,\widehat
g_{\varepsilon,\mathcal P})$ is proportional to the Gaussian
curvature measure $\mathcal K(\cdot,\widehat
g_{\varepsilon,\mathcal P})$. The ``definition'' above give us
$\mathcal K(M-\{p\},\widehat g)=4\pi$ while the ``right answer''
should be $\mathcal K(M-\{p\},\widehat g)=7\pi/2$. The problem
here is that the curvature of $\{p\}$ has its influence in the
mollifier smoothing. If we want to calculate $\mathcal
R^\kappa(M-\{p\},\widehat g)$ ``correctly'', we should separate
$\{p\}$ from the calculations. This can be done as follows:

Denote the topology of $M$ by $\mathcal T(M)$, the
$\sigma$-algebra of Borel sets of $M$ by $\mathcal B(M)$ and the
collection of all subsets of $M$ by $2^M$.

\begin{definition}
\label{geradorcurvaturemeasure} Let $(M,\widehat g)$ be a closed
non-regular Riemannian manifold. Fix $\kappa \in \{1,\ldots,n\}$.
Let $\mathcal A(M)\subset \mathcal T(M)$ be a family of open sets
such that

\begin{enumerate}

\item $\mathcal A(M)$ contains a basis of the topology $\mathcal
T(M)$;

\item For every $O\in \mathcal T(M)$, there exist a increasing
sequence of open sets $\{\widetilde O_i\}_{i \in \mathbb N}$ in
$\mathcal A(M)$ such that $\widetilde O_i\subset \subset O$ for
every $i\in \mathbb N$ and $\cup_{i=1}^\infty \widetilde O_i=O$;

\item

\[
\mathcal R^\kappa(\widetilde O,\widehat g) :=
\lim_{\varepsilon\rightarrow 0}\mathcal R^\kappa(\widetilde
O,\widehat g_{\varepsilon, \mathcal P})
\]
exists for every $\widetilde O\in \mathcal A(M)$ and this limit
does not depend on $\mathcal P$;

\item $\mathcal R^\kappa(\cdot,\widehat g):\mathcal
A(M)\rightarrow \mathbb R$ can be decomposed as
\[
\mathcal R^\kappa(\cdot,\widehat g)=\mathcal R^{\kappa
+}(\cdot,\widehat g) - \mathcal R^{\kappa -}(\cdot,\widehat g)
\]
where $\mathcal R^{\kappa +}(\cdot,\widehat g):\mathcal A(M)
\rightarrow \mathbb R$ and $\mathcal R^{\kappa -}(\cdot,\widehat
g):\mathcal A(M) \rightarrow \mathbb R$ are two non-negative
functions;

\item Let $\{\widetilde O_j\}_{j \in \mathbb N}$ be a sequence in
$\mathcal A(M)$. Denote $U_i = \cup_{j = i}^\infty \widetilde
O_j$. If $\cap_{i=1}^\infty U_i = \emptyset$, then
\[
\lim_{i \rightarrow \infty} \mathcal R^{\kappa +}(\widetilde
O_i,\widehat g) = \lim_{i \rightarrow \infty} \mathcal R^{\kappa
-}(\widetilde O_i,\widehat g) = 0;
\]

\item If $\widetilde O_1,\widetilde O_2,\widetilde O_3, \widetilde
O_4\in \mathcal A(M)$ are such that $\widetilde O_1\cup \widetilde
O_2 \subset \widetilde O_3 \cup \widetilde O_4$ and $\widetilde
O_1\cap \widetilde O_2 \subset \widetilde O_3 \cap \widetilde
O_4$, then
\[
\mathcal R^{\kappa +}(\widetilde O_1,\widehat g) + \mathcal
R^{\kappa +}(\widetilde O_2,\widehat g) \leq \mathcal R^{\kappa
+}(\widetilde O_3,\widehat g) + \mathcal R^{\kappa +}(\widetilde
O_4,\widehat g)
\]
and
\[
\mathcal R^{\kappa -}(\widetilde O_1,\widehat g) + \mathcal
R^{\kappa -}(\widetilde O_2,\widehat g) \leq \mathcal R^{\kappa
-}(\widetilde O_3,\widehat g) + \mathcal R^{\kappa -}(\widetilde
O_4,\widehat g);
\]

\item If $\{\widetilde O_i\}_{i\in \mathbb N}$ is an increasing
sequence in $\mathcal A(M)$ such that $\cup_{i=1}^\infty
\widetilde O_i = \widetilde O \in \mathcal A(M)$ and $\widetilde
O_i\subset \subset \widetilde O$, then
\[
\lim_{i\rightarrow \infty} \mathcal R^{\kappa +}(\widetilde
O_i,\widehat g)=\mathcal R^{\kappa +}(\widetilde O,\widehat g)
\]
and
\[
\lim_{i\rightarrow \infty} \mathcal R^{\kappa -}(\widetilde
O_i,\widehat g)=\mathcal R^{\kappa -}(\widetilde O,\widehat g).
\]

\end{enumerate}

Then the triple $(\mathcal A(M),\mathcal R^{\kappa +}(\cdot,
\widehat g),\mathcal R^{\kappa -}(\cdot, \widehat g))$ is called a
{\em $\kappa$-th Lipschitz-Killing curvature measure generator}.

\end{definition}

As we said before, a $\kappa$-th Lipschitz-Killing curvature
measure generator can be thought as a set of ``regular'' open sets
with their respective total $\kappa$-th Lipschitz-Killing
curvature. Observe that all the properties given in Definition
\ref{geradorcurvaturemeasure} must be satisfied if we want that
$\mathcal R^\kappa(\cdot,\widehat g)$ is extendable to a signed
measure on $\mathcal B(M)$.

\begin{remark}
\label{rmaisrmenos} The notation $\mathcal R^{\kappa\pm}
(\cdot,\widehat g)$ will be used in order to represent $\mathcal
R^{\kappa +}(\cdot,\widehat g)$ and $\mathcal R^{\kappa
-}(\cdot,\widehat g)$ simultaneously.
\end{remark}

Let $(M,\widehat g)$ be a closed non-regular Riemannian manifold
endowed with a $\kappa$-th Lipschitz-Killing curvature measure
generator $(\mathcal A(M),\mathcal R^{\kappa +}(\cdot, \widehat
g),\mathcal R^{\kappa -}(\cdot, \widehat g))$. Our aim is to
extend $\mathcal R^{\kappa\pm}(\cdot, \widehat g)$ for $\mathcal
B(M)$.

The first step is to extend $\mathcal R^{\kappa\pm}(\cdot,\widehat
g)$ to $\mathcal T(M)$. In order to do it, we prove some lemmas:

\begin{lemma}
\label{relacaodeordememR} Let $(M,\widehat g)$ be a closed
non-regular Riemannian manifold with a $\kappa$-th
Lipschitz-Killing curvature measure generator $(\mathcal
A(M),\mathcal R^{\kappa +}(\cdot,\widehat g),\mathcal R^{\kappa
-}(\cdot,\widehat g))$. Let $\widetilde O, \widetilde U \in
\mathcal A(M)$. Suppose that $\widetilde O\subset \widetilde U$.
Then
\[
\mathcal R^{\kappa \pm}(\widetilde O,\widehat g) \leq \mathcal
R^{\kappa\pm}(\widetilde U,\widehat g).
\]
\end{lemma}

\begin{proof}
Take $\widetilde O_1=\widetilde O$, $\widetilde O_2 = \widetilde
O$, $\widetilde O_3=\widetilde U$ and $\widetilde O_4 = \widetilde
O$ and the lemma is an immediate consequence of Property $(6)$ of
Definition \ref{geradorcurvaturemeasure}.

\end{proof}

\begin{definition}
\label{curvaturaemTM} Let $(M,\widehat g)$ be a closed non-regular
Riemannian manifold with a Lipschitz-Killing curvature measure
generator $(\mathcal A(M),\mathcal R^{\kappa +}(\cdot,\widehat
g),\mathcal R^{\kappa -}(\cdot,\widehat g))$. We define the
Lipschitz-Killing curvature measure of $O\in \mathcal T(M)$ by
\[
\mathcal R^\kappa(O,\widehat g)=\mathcal R^{\kappa +}(O,\widehat
g) - \mathcal R^{\kappa -}(O,\widehat g)
\]
where
\begin{equation}
\label{extensaoraberto} \mathcal R^{\kappa\pm}(O,\widehat
g):=\sup_{\widetilde O\in \mathcal A(M),\widetilde O \subset
\subset O} \mathcal R^{\kappa\pm}(\widetilde O,\widehat g).
\end{equation}
\end{definition}

\begin{lemma}
\label{somasdeabertos} Let $(M,\widehat g)$ be a closed
non-regular Riemannian manifold with a Lipschitz-Killing curvature
measure generator $(\mathcal A(M),\mathcal R^{\kappa
+}(\cdot,\widehat g),\mathcal R^{\kappa -}(\cdot,\widehat g))$.
Let $O_1$ and $O_2$ be two open sets of $M$. Then $\mathcal
R^{\kappa\pm}(O_1\cup O_2,\widehat g) + \mathcal
R^{\kappa\pm}(O_1\cap O_2,\widehat g)=\mathcal
R^{\kappa\pm}(O_1,\widehat g) + \mathcal
R^{\kappa\pm}(O_2,\widehat g)$.
\end{lemma}

\begin{proof}

\

Claim 1: $\mathcal R^{\kappa\pm} (O_1,\widehat g) + \mathcal
R^{\kappa\pm} (O_2,\widehat g) \leq \mathcal R^{\kappa\pm}
(O_1\cup O_2,\widehat g) + \mathcal R^{\kappa\pm} (O_1\cap
O_2,\widehat g)$.

Let $\widetilde O_1, \widetilde O_2 \in \mathcal A(M)$ such that
$\widetilde O_1\subset \subset O_1$ and $\widetilde O_2 \subset
\subset O_2$. Observe that $\widetilde O_1 \cup \widetilde O_2$
and $\widetilde O_1 \cap \widetilde O_2$ are compactly contained
in $O_1 \cup O_2$ and $O_1 \cap O_2$ respectively. Then there
exist $\widetilde O_3, \widetilde O_4 \in \mathcal A(M)$ such that
$\widetilde O_1 \cup \widetilde O_2 \subset \subset \widetilde O_3
\subset \subset O_1 \cup O_2$ and $\widetilde O_1 \cap \widetilde
O_2 \subset \subset \widetilde O_4 \subset \subset O_1 \cap O_2$.
Thus
\[
\mathcal R^{\kappa\pm}(\widetilde O_1, \widehat g) + \mathcal
R^{\kappa\pm}(\widetilde O_2, \widehat g) \leq \mathcal
R^{\kappa\pm}(\widetilde O_3, \widehat g) + \mathcal
R^{\kappa\pm}(\widetilde O_4, \widehat g)
\]
\[
\leq \mathcal R^{\kappa\pm}(O_1 \cup O_2, \widehat g) + \mathcal
R^{\kappa\pm}(O_1 \cap O_2, \widehat g)
\]
due to item $(6)$ of Definition \ref{geradorcurvaturemeasure}.

\

Claim 2: $\mathcal R^{\kappa\pm} (O_1\cup O_2,\widehat g) +
\mathcal R^{\kappa\pm} (O_1\cap O_2,\widehat g)\leq \mathcal
R^{\kappa\pm} (O_1,\widehat g) + \mathcal R^{\kappa\pm}
(O_2,\widehat g)$.

It follows in a similar fashion as Claim 1. Let $\widetilde O_3,
\widetilde O_4 \in \mathcal A(M)$ such that $\widetilde O_3
\subset \subset O_1 \cup O_2$ and $\widetilde O_4 \subset \subset
O_1 \cap O_2$. Then we can find $\widetilde O_1 \subset O_1$ and
$\widetilde O_2 \subset O_2$ in $\mathcal A(M)$ such that
$\widetilde O_3 \subset \subset \widetilde O_1 \cup \widetilde
O_2$ and $\widetilde O_4 \subset \subset \widetilde O_1 \cap
\widetilde O_2$. Thus
\[
\mathcal R^{\kappa\pm}(\widetilde O_3,\widehat g) + \mathcal
R^{\kappa\pm}(\widetilde O_4,\widehat g) \leq \mathcal
R^{\kappa\pm}(\widetilde O_1, \widehat g) + \mathcal
R^{\kappa\pm}(\widetilde O_2, \widehat g)
\]
\[
\leq \mathcal R^{\kappa\pm}(O_1, \widehat g) + \mathcal
R^{\kappa\pm} (O_2, \widehat g).
\]

\end{proof}

\begin{lemma}
\label{abertostendendoavaziomedidanula} Let $\{O_i\}_{i\in\mathbb
N}$ be a decreasing sequence of open sets such that
$\cap_{i=1}^\infty O_i = \emptyset$. Then
\[
\lim_{i\rightarrow \infty}\mathcal R^{\kappa\pm}(O_i,\widehat
g)=0.
\]
\end{lemma}

\begin{proof}
For every $i\in \mathbb N$, there exist a $\widetilde O_i \in
\mathcal A(M)$ such that $\widetilde O_i \subset\subset O_i$ and
$\mathcal R^{\kappa\pm} (\widetilde O_i,\widehat g)> \mathcal
R^{\kappa\pm}(O_i,\widehat g)-1/i$. If we write $V_i =
\cup_{j=i}^\infty \widetilde O_j$, then $\cap _{i=1}^\infty V_i =
\emptyset$. Thus
\[
\lim_{i \rightarrow \infty} \mathcal R^{\kappa\pm}(O_i,\widehat g)
\leq \lim_{i \rightarrow \infty} \left( \mathcal
R^{\kappa\pm}(\widetilde O_i,\widehat g) + \frac{1}{i} \right) = 0
\]
due to Property (5) of Definition \ref{geradorcurvaturemeasure}.

\end{proof}

We extend $\mathcal R^{\kappa\pm}(\cdot,\widehat g)$ to $2^M$
defining
\begin{equation}
\label{extensaormaismenos} \mathcal R^{\kappa\pm}(A,\widehat g)=
\inf_{O\in\mathcal T(M),O\supset A}\{\mathcal
R^{\kappa\pm}(O,\widehat g)\},
\end{equation}
and
\begin{equation}
\label{extensaor} \mathcal R^\kappa(A,\widehat g)=\mathcal
R^{\kappa +}(A,\widehat g) - \mathcal R^{\kappa -}(A,\widehat g).
\end{equation}
Notice that Definition (\ref{extensaor}) coincides with Definition
(\ref{extensaoraberto}) if $A\in\mathcal T(M)$. It is well known
from the classical measure theory that (\ref{extensaormaismenos})
does not define a measure on $2^M$ because it fails to be
countably additive.

Define
\begin{equation}
\label{btilm} \widetilde{\mathcal B}(M) := \{ A\in 2^M;\mathcal
R^{\kappa\pm}(A,\widehat g)=\sup_{C\subset A;
C\mathrm{\;\;compact}} \mathcal R^{\kappa\pm}(C,\widehat g) \}.
\end{equation}

We will prove that $\widetilde{\mathcal B}(M)$ is a
$\sigma$-algebra that contains $\mathcal B(M)$ and that the triple
$(M,\widetilde{\mathcal B}(M),\mathcal R^\kappa(\cdot,\widehat
g))$ is a signed measure space. The full proof of this assertion
is quite long, although not difficult. It uses just elementary
classical measure theory and it will be settled in Theorem
\ref{existenciamedidadecurvatura}.

\begin{lemma} \label{abertosecompactoscontidos} Let $(M,\widehat
g)$ be a closed non-regular Riemannian manifold with a curvature
measure generator $(\mathcal A(M),\mathcal R^{\kappa
+}(\cdot,\widehat g),\mathcal R^{\kappa -}(\cdot,\widehat g))$.
Then $\widetilde{\mathcal B}(M)$ contains every open subset of $M$
and every closed (compact) subset of $M$.
\end{lemma}

\begin{proof}
Let $O\subset M$ be an open set. Then there exist an increasing
sequence of compact sets $\{C_i\}_{i=1 \ldots \infty}$ such that
$\cup_{i=1}^\infty C_i=O$. Fix $\epsilon>0$. Observe that
$\{O-C_i\}_{i=1\ldots\infty}$ is a sequence of open sets such that
their intersection vanishes. Then
$\lim\limits_{i\rightarrow\infty} \mathcal R^{\kappa
\pm}(O-C_i,\widehat g)=0$ due to Lemma
\ref{abertostendendoavaziomedidanula}. This implies that there
exists an $N\in \mathbb N$ such that $\mathcal
R^{\kappa\pm}(O-C_i,\widehat g)<\epsilon/2$ for every $i \geq N$.
We also have that there exists an open set $\tilde O$ such that
$O\supset \supset \tilde O\supset C_N$ and $\mathcal
R^{\kappa\pm}(\tilde O,\widehat g)<\mathcal
R^{\kappa\pm}(C_N,\widehat g)+\epsilon/2$. But Lemma
\ref{somasdeabertos} implies that $\mathcal R^{\kappa\pm}(\tilde
O,\widehat g)+\mathcal R^{\kappa\pm}(O-C_N,\widehat g)=\mathcal
R^{\kappa\pm}(O,\widehat g)+\mathcal R^{\kappa\pm}(\tilde O\cap
(O-C_N),\widehat g)>\mathcal R^{\kappa\pm}(O,\widehat g)$.
Consequently we have that $\mathcal R^{\kappa\pm}(O,\widehat
g)<\mathcal R^{\kappa\pm}(C_N,\widehat g)+\epsilon$ and
$\widetilde{\mathcal B}(M)$ contains every open set of $M$.

Finally notice that every compact subset of $M$ lies in
$\widetilde{\mathcal B}(M)$ by definition, what settles the lemma.

\end{proof}

\begin{lemma}
\label{abertomaisfechadoigualavariedade} Let $(M,\widehat g)$ be a
closed non-regular Riemannian manifold with a Lipschitz-Killing
curvature measure generator $(\mathcal A(M),\mathcal R^{\kappa
+}(\cdot,\widehat g),\mathcal R^{\kappa -}(\cdot,\widehat g))$ and
let $O\subset M$ be an open set. Then $\mathcal R^{\kappa \pm}
(M,\widehat g)=\mathcal R^{\kappa \pm}(O,\widehat g)+\mathcal
R^{\kappa \pm}(M-O,\widehat g)$.
\end{lemma}

\begin{proof}
Let $\tilde O\subset M$ be an open set such that $\tilde O\supset
M-O$. Then
\[
\mathcal R^{\kappa\pm}(M,\widehat g)=\mathcal
R^{\kappa\pm}(O,\widehat g)+\mathcal R^{\kappa\pm}(\tilde
O,\widehat g)-\mathcal R^{\kappa\pm}(O\cap\tilde O,\widehat g),
\]
due to Lemma \ref{somasdeabertos}. If we take a decreasing
sequence of open sets $\{\tilde O_i\}_{i=1\ldots\infty}$ such that
$\bigcap\limits_{i=1}^\infty \tilde O_i=M-O$, then
\[
\begin{array}{ccc}
\mathcal R^{\kappa\pm}(M,\widehat g) & = &
\lim\limits_{i\rightarrow\infty}\left[\mathcal
R^{\kappa\pm}(O,\widehat g)+\mathcal R^{\kappa\pm}(\tilde
O_i,\widehat g)-\mathcal R^{\kappa\pm}(O\cap
\tilde O_i,\widehat g)\right] \\
& = & \mathcal R^{\kappa\pm}(O,\widehat g)+\mathcal
R^{\kappa\pm}(M-O,\widehat g)
\end{array}
\]
due to Lemma \ref{abertostendendoavaziomedidanula}.

\end{proof}

\begin{lemma}
\label{propriedadesaditivas} Let $(M,\widehat g)$ be a closed
non-regular Riemannian manifold with Lipschitz-Killing curvature
measure generator $(\mathcal A(M),\mathcal R^{\kappa
+}(\cdot,\widehat g),\mathcal R^{\kappa -}(\cdot,\widehat g))$.
Let $\{A_i\}_{i=1\ldots\infty}$ be a countable collection of
subsets of $M$. Denote $A=\bigcup\limits_{i=1}^\infty A_i$. Then

\begin{enumerate}

\item $\mathcal R^{\kappa\pm}(A,\widehat g)\leq \sum_{i=1}^\infty
\mathcal R^{\kappa\pm}(A_i,\widehat g)$;

\item If $\{A_i\}_{i=1\ldots N}$ is a finite collection of
pairwise disjoint closed sets, then $\mathcal
R^{\kappa\pm}(A,\widehat g)= \sum_{i=1}^N \mathcal
R^{\kappa\pm}(A_i,\widehat g)$.

\end{enumerate}

\end{lemma}

\begin{proof}

\

\

Assertion (1)

\

Fix $\epsilon>0$. For each $A_i$ there exists an open set
$O_i\supset A_i$ such that $\mathcal R^{\kappa\pm}(O_i,\widehat
g)<\mathcal R^{\kappa\pm}(A_i,\widehat g)+\epsilon.2^{-i}$. \
Denote $O=\bigcup_{i=1}^\infty O_i$. From Lemma
\ref{somasdeabertos}, it follows that $\mathcal
R^{\kappa\pm}(O,\widehat g) \leq \sum\limits_{i=1}^\infty \mathcal
R^{\kappa\pm}(O_i,\widehat g)$. Therefore $O$ is an open set that
contains $A$ and satisfies $\mathcal R^{\kappa\pm}(O,\widehat
g)<\sum\limits_{i=1}^\infty \mathcal R^{\kappa\pm}(A_i,\widehat
g)$ $+ \epsilon$, what proves this assertion.

\

Assertion (2)

\

Fix $\epsilon >0$. Observe that we can choose pairwise disjoint
open sets $\tilde O_i\supset A_i$ such that $\mathcal
R^{\kappa\pm}(\tilde O_i,\widehat g)<\mathcal
R^{\kappa\pm}(A_i,\widehat g)+\epsilon.2^{-i}$. Choose also an
open set $\tilde O\supset A$ such that $\mathcal
R^{\kappa\pm}(\tilde O,\widehat g)<\mathcal
R^{\kappa\pm}(A,\widehat g)+\epsilon$. Now denote $O_i=\tilde
O_i\cap \tilde O$ and $O=\cup_{i=1}^N O_i$. From Lemma
\ref{somasdeabertos}, it follows that $\mathcal
R^{\kappa\pm}(O,\widehat g)= \sum\limits_{i=1}^N \mathcal
R^{\kappa\pm}(O_i,\widehat g)$. Then
\[
\begin{array}{rcccl}
\sum\limits_{i=1}^N \mathcal R^{\kappa\pm}(A_i,\widehat g)
-\epsilon & \leq & \sum\limits_{i=1}^N \mathcal
R^{\kappa\pm}(O_i,\widehat g)-\epsilon & = &
\mathcal R^{\kappa\pm}(O,\widehat g)-\epsilon \\
& < & \mathcal R^{\kappa\pm}(A,\widehat g) & \leq & \mathcal R^{\kappa\pm}(O,\widehat g) \\
& = & \sum\limits_{i=1}^N \mathcal R^{\kappa\pm}(O_i,\widehat g) &
\leq & \sum\limits_{i=1}^N \mathcal R^{\kappa\pm}(A_i,\widehat
g)+\epsilon
\end{array}
\]
what proves this assertion.

\end{proof}

Finally we are able to prove the main theorem of this section:

\begin{theorem}
\label{existenciamedidadecurvatura} Let $(M,\widehat g)$ be a
closed non-regular Riemannian manifold. Let $(\mathcal A(M),$
$\mathcal R^{\kappa +}(\cdot,\widehat g),\mathcal R^{\kappa
-}(\cdot,\widehat g))$ be a Lipschitz-Killing \ curvature measure
generator and consider $\widetilde{\mathcal B}(M)$ as defined in
(\ref{btilm}). Then

\begin{enumerate}

\item $\widetilde{\mathcal B}(M)$ is a $\sigma$-algebra that
contains $\mathcal B(M)$.

\item $(M,\widetilde{\mathcal B}(M),\mathcal
R^\kappa(\cdot,\widehat g))$ is a (signed) measure space.

\end{enumerate}

\end{theorem}

\begin{proof}
First we prove that $\widetilde{\mathcal B}(M)$ is a
$\sigma$-algebra.

It is obvious that $M\in \widetilde{\mathcal B}(M)$.

\

{\em Claim 1:} If $A\in\widetilde{\mathcal B}(M)$, then $M-A\in
\widetilde{\mathcal B}(M)$.

\

Suppose that $A\in\widetilde{\mathcal B}(M)$. Lemma
\ref{abertomaisfechadoigualavariedade} implies that
\[
\begin{array}{ccl}
\mathcal R^{\kappa\pm}(M-A,\widehat g) & = &
\inf\limits_{O\in\mathcal T(M),O\supset
M-A}\mathcal R^{\kappa\pm}(O,\widehat g)\\
& = & \inf\limits_{O\in\mathcal T(M),O\supset
M-A}\left[ \mathcal R^{\kappa\pm}(M,\widehat g)- \mathcal R^{\kappa\pm}(M-O,\widehat g)\right] \\
& = & \mathcal R^{\kappa\pm}(M,\widehat
g)-\sup\limits_{O\in\mathcal T(M),O\supset M-A}
\mathcal R^{\kappa\pm}(M-O,\widehat g) \\
& = & \mathcal R^{\kappa\pm}(M,\widehat
g)-\sup\limits_{M-O\mathrm{\;\;compact},M-O\subset
A}\mathcal R^{\kappa\pm}(M-O,\widehat g)\\
& = &
\mathcal R^{\kappa\pm}(M,\widehat g)-\mathcal R^{\kappa\pm}(A,\widehat g) \\
& = & \mathcal R^{\kappa\pm}(M,\widehat
g)-\inf\limits_{O\in\mathcal T(M),O\supset
A}\mathcal R^{\kappa\pm}(O,\widehat g) \\
& = & \mathcal R^{\kappa\pm}(M,\widehat g)-\inf\limits_{M-O\;\;
\mathrm{compact},M-O\subset
M-A}\mathcal R^{\kappa\pm}(O,\widehat g) \\
& = & \sup\limits_{M-O\;\; \mathrm{compact},M-O\subset
M-A}\mathcal R^{\kappa\pm}(M,\widehat g)-\mathcal R^{\kappa\pm}(O,\widehat g) \\
& = & \sup\limits_{M-O\;\; \mathrm{compact},M-O\subset
M-A}\mathcal R^{\kappa\pm}(M-O,\widehat g).
\end{array}
\]
Thus $M-A\in\widetilde{\mathcal B}(M)$.

\

{\em Claim 2:} If $A_1,A_2\in \widetilde{\mathcal B}(M)$, then
$A_1\cap A_2\in \widetilde{\mathcal B}(M)$.

\

For every $\epsilon>0$, there exist an open set $O_i\supset A_i$
and a compact set $C_i\subset A_i$, $i=1,2$, such that
\[
\mathcal R^{\kappa\pm}(O_i,\widehat g)-\frac{\epsilon}{4}<\mathcal
R^{\kappa\pm}(A_i,\widehat g)<\mathcal R^{\kappa\pm}(C_i,\widehat
g)+ \frac{\epsilon}{4},
\]
what gives
\[
\mathcal R^{\kappa\pm}(O_1,\widehat g)+\mathcal
R^{\kappa\pm}(O_2,\widehat g)-\frac{\epsilon}{2}<\mathcal
R^{\kappa\pm}(C_1,\widehat g)+ \mathcal R^{\kappa\pm}(C_2,\widehat
g)+\frac{\epsilon}{2}.
\]
Lemma \ref{somasdeabertos} implies
\[
\mathcal R^{\kappa\pm}(O_1\cup O_2,\widehat g)+\mathcal
R^{\kappa\pm}(O_1\cap O_2,\widehat g)-\frac{\epsilon}{2} =
\mathcal R^{\kappa\pm}(O_1,\widehat g)+\mathcal
R^{\kappa\pm}(O_2,\widehat g)-\frac{\epsilon}{2},
\]
and using Lemma \ref{abertomaisfechadoigualavariedade} we have
that
\[
\begin{array}{cl}
& \mathcal R^{\kappa\pm}(O_1\cap O_2,\widehat g) - \frac{\epsilon}{2} \\
< & \mathcal R^{\kappa\pm}(C_1,\widehat g) + \mathcal
R^{\kappa\pm}(C_2,\widehat g) - \mathcal R^{\kappa\pm}(O_1\cup
O_2,\widehat g) + \frac{\epsilon}{2}
\\
= & \mathcal R^{\kappa\pm}(M,\widehat g)-\mathcal
R^{\kappa\pm}(M-C_1,\widehat g)+\mathcal R^{\kappa\pm}(M,\widehat
g)-\mathcal R^{\kappa\pm}(M-C_2,\widehat g)
\\
& -\mathcal R^{\kappa\pm}(O_1\cup O_2,\widehat g)+\frac{\epsilon}{2} \\
= & 2 \mathcal R^{\kappa\pm}(M,\widehat g)-\mathcal
R^{\kappa\pm}(M-(C_1\cap C_2),\widehat g)-\mathcal
R^{\kappa\pm}(M-(C_1\cup
C_2),\widehat g) \\
& -\mathcal R^{\kappa\pm}(O_1\cup O_2,\widehat g)+\frac{\epsilon}{2} \\
= & \mathcal R^{\kappa\pm}(C_1\cap C_2,\widehat g)+\mathcal
R^{\kappa\pm}(C_1\cup C_2,\widehat g)-\mathcal
R^{\kappa\pm}(O_1\cup
O_2,\widehat g)+\frac{\epsilon}{2} \\
< & \mathcal R^{\kappa\pm}(C_1\cap C_2,\widehat
g)+\frac{\epsilon}{2}.
\end{array}
\]
Now notice that $(C_1\cap C_2)\subset (A_1\cap A_2)$ and
\[
\mathcal R^{\kappa\pm}(A_1\cap A_2,\widehat g)<\mathcal
R^{\kappa\pm}(O_1\cap O_2,\widehat g)<\mathcal
R^{\kappa\pm}(C_1\cap C_2,\widehat g)+\epsilon.
\]
Thus $A_1\cap A_2\in\widetilde{\mathcal B}(M)$.

\

{\em Claim 3:} If $\{A_i\}_{i=1\ldots\infty}$ is a {\em disjoint}
countable collection of subsets in $\widetilde{\mathcal B}(M)$,
then $A:=\bigcup\limits_{i=1}^\infty A_i\in \widetilde{\mathcal
B}(M)$.

\

Fix $\epsilon>0$. Then for every $i$ there exist a compact set
$C_i\subset A_i$ such that

\[
\mathcal R^{\kappa\pm}(A_i,\widehat g)<\mathcal
R^{\kappa\pm}(C_i,\widehat g)+\epsilon.2^{-i-1}
\]
what implies
\begin{equation}
\label{aproximacaofechado} \sum_{i=1}^\infty \mathcal
R^{\kappa\pm}(A_i,\widehat g) < \sum_{i=1}^\infty \mathcal
R^{\kappa\pm}(C_i,\widehat g)+\epsilon/2.
\end{equation}
Observe that there exist $N\in\mathbb N$ such that
\begin{equation}
\label{aproximacaofinita} \sum_{i=1}^\infty \mathcal
R^{\kappa\pm}(C_i,\widehat g)<\sum_{i=1}^N \mathcal R^{\kappa
\pm}(C_i,\widehat g)+\epsilon/2.
\end{equation}
Denote $C=\bigcup\limits_{i=1}^N C_i$. Eqs.
(\ref{aproximacaofechado}) and (\ref{aproximacaofinita}) and Lemma
\ref{propriedadesaditivas} gives
\begin{equation}
\begin{array}{rcccl}
\mathcal R^{\kappa\pm}(A,\widehat g) & \leq & \sum_{i=1}^\infty
\mathcal R^{\kappa\pm}(A_i,\widehat g) & \leq &
\sum_{i=1}^\infty \mathcal R^{\kappa\pm}(C_i,\widehat g)+\epsilon/2 \\
& \leq & \sum_{i=1}^N \mathcal R^{\kappa\pm}(C_i,\widehat
g)+\epsilon & = & \mathcal R^{\kappa\pm}(C,\widehat
g)+\epsilon.\label{cadeiamajoracaofechado}
\end{array}
\end{equation}
We found a compact set $C\subset A$ such that the inequality above
is satisfied, what implies that $A\in\widetilde{\mathcal B}(M)$.

\

{\em Claim 4:} If $\{A_i\}_{i=1\ldots\infty}$ is {\em any}
countable collection of subsets in $\widetilde{\mathcal B}(M)$,
then $\bigcup\limits_{i=1}^\infty A_i\in \widetilde{\mathcal
B}(M)$.

\

Let $A_1,A_2\in\widetilde{\mathcal B}(M)$. Using Claims 1, 2 and
3, we can prove that $(A_1-A_2)$ and $(A_2-A_1)$ lie in
$\widetilde{\mathcal B}(M)$. Observe that $A_1\cup A_2$ is the
disjoint union $(A_1\cap A_2)\cup (A_1-A_2)\cup (A_2-A_1)$ and it
also lie in $\widetilde{\mathcal B}(M)$. Similarly we can prove
that $\bigcup\limits_{i=1}^\infty A_i$ can be written as the union
of disjoint sets in $\widetilde{\mathcal B}(M)$. Consequently
Claim 3 implies that $\bigcup\limits_{i=1}^\infty A_i\in
\widetilde{\mathcal B}(M)$. Thus $\widetilde{\mathcal B}(M)$ is a
$\sigma$-algebra that contains $\mathcal T(M)$, what implies that
it also contains $\mathcal B(M)$.

\

{\em Claim 5:} $(M,\widetilde{\mathcal B}(M), \mathcal
R^\kappa(\cdot,\widehat g))$ is a (signed) measure space.

\

Take a pairwise disjoint sequence $\{A_i\}_{i=1\ldots\infty}$,
$A_i\in \widetilde{\mathcal B}(M)$. Eq.
(\ref{cadeiamajoracaofechado}) gives
\[
\mathcal R^{\kappa\pm}(A,\widehat g) \leq \sum_{i=1}^\infty
\mathcal R^{\kappa\pm}(A_i,\widehat g) \leq \mathcal
R^{\kappa\pm}(C,\widehat g)+\epsilon\leq \mathcal
R^{\kappa\pm}(A,\widehat g)+\epsilon
\]
and the equality
\[
\mathcal R^{\kappa\pm}(A,\widehat g)=\sum_{i=1}^\infty \mathcal
R^{\kappa\pm}(A_i,\widehat g)
\]
follows.

\end{proof}

\begin{definition}
\label{definicaomedidadecurvatura} Let $(M,\widehat g)$ be a
closed non-regular Riemannian manifold. Let $(\mathcal A(M),$
$\mathcal R^{\kappa +}(\cdot,\widehat g),\mathcal R^{\kappa
-}(\cdot,\widehat g))$ be a Lipschitz-Killing \ curvature \
measure generator and consider $\widetilde{\mathcal B}(M)$ as
defined in (\ref{btilm}). The function $\mathcal
R^\kappa(\cdot,\widehat g): \widetilde{\mathcal B}(M) \rightarrow
\mathbb R$ is the {\em Lipschitz-Killing curvature measure} of
$(M,\widehat g)$ (with respect to $(\mathcal A(M),$ $\mathcal
R^{\kappa +}(\cdot,\widehat g),\mathcal R^{\kappa
-}(\cdot,\widehat g))$).
\end{definition}

\begin{remark}
\label{approachnaounico} We defined the $\kappa$-th
Lipschitz-Killing curvature measure from a collection $(\mathcal
A(M),\mathcal R^{\kappa +} (\cdot, \widehat g),\mathcal R^{\kappa
-}(\cdot, \widehat g))$ because it is useful in some situations.
For instance, it is adequate to study the Gaussian curvature
measure for piecewise smooth two-dimensional Riemannian manifolds
(See Section \ref{gaussbonnetnaoregular}). However we can consider
other approaches if we want to study other problems. The most
important thing here is to show that the mollifier smoothing with
respect to $\mathcal P$ is a useful tool to study non-regular
Riemannian manifolds.
\end{remark}

\section{The Gaussian curvature measure for piecewise
smooth two-dimensional Riemannian manifolds}

\label{gaussbonnetnaoregular}

In this section we prove the existence of a ``natural'' Gaussian
curvature measure generator for closed oriented piecewise smooth
two-dimensional Riemannian manifolds. We get a Gaussian curvature
measure $\mathcal K(\cdot, \widehat g):\mathcal B(M) \rightarrow
\mathbb R$ for this kind of surfaces and we show that it has the
expected values for some subsets of $M$ (See Theorem
\ref{curvaturasuavizacaoconverge}). As a direct consequence, we
get an alternative proof of the (well known) generalization of the
Gauss-Bonnet theorem for this class of surfaces.

Let $(M^2,\widehat g)$ be a closed oriented piecewise smooth
Riemannian manifold. Let $\Theta:\Upsilon \rightarrow M$ be a
triangulation associated to $(M,\widehat g)$. Fix a vertex $x\in
V$. The point $\Theta^{-1}(x)$ lies to some simplexes, let us say
$\Xi_1,\ldots,\Xi_k$. For every triangle $\Theta(\Xi_i)\subset M$,
we can associate the internal angle $\varsigma_i(x)$ of $x$. {\em
We suppose that $\varsigma_i(x)\neq 0$ for every $i=1, \ldots,
k$}. Consider the number $K_0(x)=2\pi-\sum_{i=1}^k
\varsigma_i(x)$. This is intuitively the ``zero-dimensional
curvature'' of a vertex as explained in the example of the cube in
the introduction.

Now consider a point $x\in E$. The point $\Theta^{-1}(x)$ lies to
two simplexes $\Xi_1$ and $\Xi_2$. The sum of the geodesic
curvature of $E$ at $x$ with respect to $\Theta(\Xi_1)$ and
$\Theta(\Xi_2)$ will be denoted by $K_1(x)$. This is intuitively
the ``one-dimensional curvature'' of a point that lies on a edge
of $M$.

Let $x\in F$. We denote the Gaussian curvature at $x$ by $K_2(x)$.

Consider an open set $O\subset M$ such that its boundary $\partial
O$ is the image of a piecewise regular curve. Let $\alpha: [a,b]
\rightarrow M$ be a positive regular parametrization of $\partial
O$ such that $\alpha(a)=\alpha(b)$. The external angle of
$\partial O$ at $\alpha(t_0)\in
\partial O - V - E$ is the angle in the interval $(-\pi,\pi)$
formed by $\lim_{t\rightarrow t_0^-} \alpha^\prime(t)$ and
$\lim_{t\rightarrow t_0^+} \alpha^\prime(t)$ (in this order). We
will not allow the angles $\pm\pi$.

Although the following theorem is well known, we prove it here for
the sake of completeness.

\begin{theorem}
\label{curvaturaemabertossuavesporpartes} Let $(M^2,\widehat g)$
be an closed and oriented piecewise smooth Riemannian manifold.
Let $O\subset M$ be an open set with a piecewise smooth boundary
$\partial O$ such that:

\begin{enumerate}

\item $\partial O$ is smooth outside $\{q_1,\ldots,q_k\}$.

\item $\partial O$ does not intercept $V$.

\item $\{q_1,\ldots,q_k\}$ does not intercept $E$.

\item $\partial O$ intercepts $E$ transversally (this includes the
case $\partial O \cap E=\emptyset$).

\end{enumerate}

Suppose that the external angles at $q_i$, i=1,\ldots,k, are given
by $\varrho(q_i)$. Then

\[
\sum_{i=1}^k \varrho(q_i) + \int_{\partial O} k_{\widehat
g}(x).ds_{\widehat g}
\]
\begin{equation}
\label{mudancadeangulo} = 2.\pi.\chi(O \cup \partial O) - \sum_{
O\cap V} K_0 (x) - \int_{O \cap E} K_1(x).ds_{\widehat g} -
\int\int_{O \cap F} K_2(x).dV_{\widehat g}
\end{equation}
where $k_{\widehat g}$ is the geodesic curvature of $\partial O$
with respect to $O$, $ds_{\widehat g}$ denotes the length element
and $\chi(O \cup
\partial O)$ denotes the Euler characteristic of $O \cup \partial
O$.

\end{theorem}

\begin{proof}
This proof is analogous to the proof of the classical Gauss-Bonnet
Theorem (For instance, see \cite{20}).

First of all, we can take a triangulation $\Theta^\prime:
\Upsilon^\prime \rightarrow (O\cup
\partial O)$ associated to $(O \cup \partial O,\widehat g\vert_{O \cup \partial
O})$. Let $\# F$, $\# E$ and $\# V$ be the number of faces, edges
and vertices of $\Upsilon^\prime$. Let $\# B$ be the number of
vertices (and edges) on the boundary of $\Upsilon^\prime$. Formula
(\ref{mudancadeangulo}) holds for a simplex $\Xi^\prime$, that is,
\[
\int \int_{\Theta^\prime(\Xi^\prime)} K_2(x).dV_{\widehat g} =
-\pi - \int_{\Theta^\prime(\partial \Xi^\prime)} k_{\widehat
g}(x).ds_{\widehat g}+\sum_{i=1}^3 \varsigma_i.
\]
Summing up the formula above for every simplex, we have that
\[
\int\int_{O} K_2(x).dV_{\widehat g}
\]
\[
= -3\pi(\# F) + 2\pi(\# F) - \int_{O\cap E} K_1(x).ds_{\widehat g}
- \int_{\partial O} k_{\widehat g}(x).ds_{\widehat g}
\]
\begin{equation}
\label{primeiraformulagauss} + \sum_{V\in O} \sum \varsigma_i +
\sum_{V\in
\partial O} \sum \varsigma_i=(*).
\end{equation}
But $3(\# F)=2(\# E)-(\# B)$, what implies that
\[
(*)=-2\pi(\# E) + 2\pi(\# F) + 2\pi(\# V) - \int_{O\cap E}
K_1(x).ds_{\widehat g} - \int_{\partial O} k_g(x).ds_{\widehat g}
\]
\[
- \sum_{V\in O} \left( 2\pi - \sum \varsigma_i \right) -
\sum_{V\in
\partial O} \left( \pi - \sum \varsigma_i \right)
\]
\[
=2.\pi.\chi(O\cup \partial O) - \int_{O\cap E} K_1(x).ds_{\widehat
g} - \int_{\partial O} k_g(x).ds_{\widehat g}
\]
\[
- \sum_{V\in O} K_0(x) - \sum_{V\in
\partial O} \varrho(q_i)
\]
what settles the Theorem.

\end{proof}

The following generalization of the classical Gauss-Bonnet Theorem
is an immediate consequence of the proof of Theorem
\ref{curvaturaemabertossuavesporpartes}.

\begin{theorem}
\label{gaussbonnetgeneralizado} Let $(M^2,\widehat g)$ be a closed
and oriented piecewise smooth Riemannian manifold. Then
\[
2.\pi.\chi(M)=\sum_{V} K_0(x)+\int_{E} K_1(x).ds_{\widehat g} +
\int\int_F K_2(x).dV_{\widehat g}
\]
where $ds_{\widehat g}$ is the length element of the edge.
\end{theorem}

The next theorem states that for piecewise smooth two-dimensional
Riemannian manifold, $K_0(x)$, $K_1(x)$ and $K_2(x)$ can be really
interpreted as the zero dimensional curvature of the vertex, the
one dimensional curvature of the edge and the two dimensional
curvature of the face respectively (See Eq.
(\ref{convergenciacurvatura})).

\begin{theorem}
\label{curvaturasuavizacaoconverge} Let $(M^2,\widehat g)$ be an
closed and oriented piecewise smooth Riemannian manifold. Let
$O\subset M$ be an open set with a piecewise smooth boundary
$\partial O$ such that:

\begin{enumerate}

\item $\partial O$ is smooth outside $\{q_1,\ldots,q_k\}$.

\item $\partial O$ does not intercept $V$.

\item $\{q_1,\ldots,q_k\}$ does not intercept $E$.

\item $\partial O$ intercepts $E$ transversally (this includes the
case $\partial O \cap E=\emptyset$).

\end{enumerate}

Let $\mathcal P$ be a locally finite covering $\{(U_\omega \subset
\subset O_\omega,\widetilde e_\omega)\}_{\omega \in \Lambda}$
together with the partition of unity $\{\psi_\omega\}_{\omega \in
\Lambda}$ subordinated to $\{U_\omega \}_{\omega \in \Lambda}$.
Let $\widehat g_{\varepsilon,\mathcal P}$ be the mollifier
smoothing of $\widehat g$ with respect to $\mathcal P$. Suppose
that the external angles of $\partial O$ at $q_i\in (M,\widehat
g_{\varepsilon,\mathcal P})$, i=1,\ldots,k, is given by
$\varrho_\varepsilon(q_i)$. Then

\begin{equation}
\label{deflexaoangular} \lim_{\varepsilon \rightarrow 0} \left(
\sum_{i=1}^k \varrho_\varepsilon(q_i) + \int_{\partial O}
k_{g_{\varepsilon,\mathcal P}}(x).ds_{\widehat g} \right) =
\sum_{i=1}^k \varrho(q_i) + \int_{\partial O} k_{\widehat
g}(x).ds_{\widehat g}
\end{equation}
and
\begin{equation}
\label{convergenciacurvatura} \lim_{\varepsilon \rightarrow 0}
\left( \int\int_{O} K_2(x).dV_{\widehat g_{\varepsilon, \mathcal
P}} \right) = \sum_{ O\cap V} K_0 (x) + \int_{O \cap E}
K_1(x).ds_{\widehat g} + \int\int_{O \cap F} K_2(x).dV_{\widehat
g}.
\end{equation}

\end{theorem}

\begin{proof}
We have that
\[
\lim_{\varepsilon \rightarrow 0} \left( \sum_{i=1}^k
\varrho_\varepsilon(q_i) \right) = \sum_{i=1}^k \varrho(q_i)
\]
because $\widehat g_{\varepsilon, \mathcal P}$ converges to
$\widehat g$ on $F$.

In order to see that
\[
\lim_{\varepsilon \rightarrow 0}  \int_{\partial O}
k_{g_{\varepsilon,\mathcal P}}(x).ds_{\widehat g} = \int_{\partial
O} k_{\widehat g}(x).ds_{\widehat g}
\]
we split the integral through $\partial O$ in two parts: the
integral near the edges and the integral far from the edges.

Let $y\in E\cap \partial O$. Theorem \ref{passagemporaresta}
implies that
\[
\lim_{\varepsilon \rightarrow 0} \int_{\partial O\cap V_y}
k_{\widehat g_{\varepsilon,\mathcal P}}(x).ds_{\widehat g}
\]
can be made as small as we want if we choose a sufficiently small
neighborhood $V_y$ of $y$.

Theorem \ref{suavizacao2boa} implies that the Levi-Civita
connection of $\widehat g_{\varepsilon,\mathcal P}$ converges
uniformly to the Levi-Civita connection of $\widehat g$ on compact
subsets of $F$. Then
\[
\lim_{\varepsilon \rightarrow 0} \int_{\partial O \cap \widetilde
O} k_{\widehat g_{\varepsilon,\mathcal P}}(x).ds_{\widehat g} =
\int_{\partial
 O \cap \widetilde O} k_{\widehat g}(x).ds_{\widehat g}
\]
for every open set $\widetilde O \subset \subset F$.

Therefore
\[
\lim_{\varepsilon \rightarrow 0} \int_{\partial O} k_{\widehat
g_{\varepsilon,\mathcal P}}(x).ds_{\widehat g} = \int_{\partial
 O} k_{\widehat g}(x).ds_{\widehat g}.
\]
and Eq. (\ref{deflexaoangular}) follows. Equation
(\ref{convergenciacurvatura}) is a direct consequence of
(\ref{deflexaoangular}) and Theorem
\ref{curvaturaemabertossuavesporpartes}.

\end{proof}

Let $(M^2,\widehat g)$ be a closed and oriented piecewise smooth
Riemannian manifold. Let us define a Gaussian curvature measure
generator. Define $\mathcal A(M)$ as the family of open sets that
satisfies the conditions given in Theorem
\ref{curvaturaemabertossuavesporpartes}. Define $K^-_0(x)=\max
(-K_0(x),0)$, $K^+_0(x)=\max (K_0(x),0)$, $K^-_1(x)=\max
(-K_1(x),0)$, $K^+_1(x)=\max (K_1(x),0)$, $K^-_2(x)=$ $\max
(-K_2(x),0)$, $K^+_2(x)$, $=\max (K_2(x),0)$. Define $\mathcal
K^\pm(\cdot,\widehat g):\mathcal A(M) \rightarrow \mathbb R$ as
\begin{equation}
\label{rmaismenospoliedro} \mathcal K^\pm(\widetilde O,\widehat
g):= \sum_{ \widetilde O\cap V} K^\pm_0 (x) + \int_{\widetilde O
\cap E} K^\pm_1(x).ds_{\widehat g} + \int\int_{\widetilde O \cap
F}K^\pm_2(x).dV_{\widehat g},
\end{equation}
where the superscripts of $K^\pm_i$ and $\mathcal K^\pm$ are
explained in Remark \ref{rmaisrmenos}.

\begin{theorem}
The triple $(\mathcal A(M),$ $\mathcal K^-(\cdot,\widehat
g),\mathcal K^+(\cdot,\widehat g))$ defined by Eq.
(\ref{rmaismenospoliedro}) is a Gaussian curvature measure
generator.
\end{theorem}

\begin{proof}
Properties (1), (3), (4), (5), (6) and (7) of Definition
\ref{geradorcurvaturemeasure} hold due to the definition of
$\mathcal A(M)$, (\ref{convergenciacurvatura}) and
(\ref{rmaismenospoliedro}).

Let us see that Property (2) holds: Let $O\in \mathcal T(M)$.
Denote $C=M-O$. Put an arbitrary smooth Riemannian metric $\check
g$ on $M$. Define
\[
C(n)=\{x\in O;d_{\check g}(x,C)\geq 1/n \}
\]
for every $n \in \mathbb N$. Notice that $C(n)$ is a compact set.
Cover $C(n)$ by elements of $\mathcal A(M)$ such that they are
compactly contained on $O$. Take a finite subcover
$\{O_1,\ldots,O_k\}$ of this covering. Now notice that it is
possible to make small deformations in $O_i$, $i=1,\ldots,k$, in
such a way that
\begin{enumerate}

\item $C(n) \subset \subset \bigcup\limits_{i=1}^k O_i \subset
\subset O$.

\item The vertices of $O_j$ do not intercept vertices and edges of
other open sets $O_l$, $l \neq j$.

\item All the edges of the open sets $O_i$, $i=1,\ldots, k$,
intercept themselves transversally.

\item All the edges of the open sets $O_i$, $i=1,\ldots, k$,
intercept $E$ transversally.

\item $\partial O_i$ do not intercept $V$ for every
$i=1,\ldots,k$.

\end{enumerate}

Observe that $\bigcup\limits_{i=1}^k O_i$ is an open set
satisfying all the conditions given in Theorem
\ref{curvaturaemabertossuavesporpartes} except, eventually, (3).
But we can make a further perturbation on $\widetilde O(n) :=
\bigcup\limits_{i=1}^k O_i$ in order to satisfy (3). Then, for
every $n \in \mathbb N$, we can find an open set $\widetilde
O(n)\in \mathcal A(M)$ such that $C(n)\subset \subset \widetilde
O(n) \subset \subset O$ and these sets can be built in such a way
that $\widetilde O(n) \subset \widetilde O(n+1)$. Therefore
Property (2) of Definition \ref{geradorcurvaturemeasure} holds and
the triple $(\mathcal A(M),\mathcal K^+(\cdot,\widehat g),\mathcal
K^-(\cdot,\widehat g))$ is a Gaussian curvature measure generator.

\end{proof}

Theorem \ref{existenciamedidadecurvatura} implies that the
Gaussian curvature measure is extendable for every Borel set. If
$O\in \mathcal T(M)$, then it is not difficult to see that
\begin{equation}
\label{curvaturaemabertosparapoliedros} \mathcal K^\pm(O,\widehat
g)= \sum_{O \cap V} K^\pm_0 (x) + \int_{O \cap E}
K^\pm_1(x).ds_{\widehat g} + \int\int_{O \cap F}
K^\pm_2(x).dV_{\widehat g}.
\end{equation}

It is not difficult to see either that if $x\in V$, then
\begin{equation}
\label{curvaturadosvertices} \mathcal K^\pm(\{x\},\widehat g) =
K_0(x)
\end{equation}
and if $\upsilon \subset E$, then
\begin{equation}
\label{curvaturadasarestas} \mathcal K^\pm(\upsilon,\widehat g) =
\int_{\upsilon} K_1(x).ds_{\widehat g}.
\end{equation}

Therefore $\mathcal K(\cdot,\widehat g):\mathcal B(M)\rightarrow
\mathbb R$ is a Gaussian curvature measure which gives the
expected geometrical values for some subsets of $M$.

\section{Curvature dimension and an instructive non-regular example}

\label{secaocurvaturadimensaofracionaria}

In Section \ref{gaussbonnetnaoregular}, we studied piecewise
smooth two-dimensional Riemannian manifolds. In particular,
(\ref{curvaturadosvertices}) and (\ref{curvaturadasarestas}) show
that the ``curvature integra'' can be calculated in subsets with
dimension less than two. Moreover $V$, $E$ and $F$ are the
``natural'' subsets where the ``zero-dimensional curvature'',
``one-dimensional curvature'' and the ``two-dimensional
curvature'' arise respectively. A natural question here is whether
the ``curvature integra'' can arise naturally in subsets with
non-integer Hausdorff dimension (For Hausdorff dimension, see, for
instance, \cite{19}). In this Section, we present a natural
candidate to answer this question affirmatively.

When we study spaces with non-integer Hausdorff dimension, the
natural starting point is the classical Cantor ternary set, which
we denote by $\mathcal{C}$. It is the subset of $[0,1]$ which is
composed by numbers that have some ternary representation without
the number `$1$' (See, for instance, \cite{17}). It is well known
that the Hausdorff dimension of $\mathcal C$ is equal to $\ln
2/\ln 3$ (See, for instance, \cite{13}). One easy way to
understand $\mathcal C$ is to take out the open middle third
interval (that is, $(1/3,2/3)$) from $[0,1]$. After that, we take
out the open middle third intervals from every remaining connected
set (At this point, the connected sets we are talking about are
$[0,1/3]$ and $[2/3,1]$). We do it recursively and the remaining
set is $\mathcal C$.

We have the Cantor ternary function associated with $\mathcal C$,
which we denote by $f_{\mathcal C}:[0,1]\rightarrow \mathbb R$.
Let $x=0.a_1a_2a_3\ldots$ be a ternary representation of $x\in
[0,1]$. Denote by $N:[0,1]\rightarrow \mathbb N$ the first $i$
such that $a_i=1$. If the number `1' is absent from the ternary
representation of $x$, make $N(x)=\infty$. Now define
\[
f_{\mathcal C}(x)=\frac{1}{2}\sum_{i=1}^{N(x)}\frac{a_i}{2^i}.
\]
It is not difficult to prove that $f_{\mathcal C}$ is well
defined. Moreover it is a continuous increasing function such that
its derivative is equal to zero outside $\mathcal C$.

We will built a curve of class $C^1$ in $\mathbb R^2$ which is
homeomorphic to a circle using $f_{\mathcal C}$. It will be flat
almost everywhere and the curvature will be concentrated in a set
with non-integer dimension. Define $\theta_{\mathcal
C}:[0,1]\rightarrow\mathbb R$ by
\[
\theta_{\mathcal C}(t)=\left\{
\begin{array}{ccc}
 2\pi f_{\mathcal C}(\frac{4t}{3})  & \mathrm{if} &
 t\in[0,\frac{3}{4}]\\
 2\pi & \mathrm{if} & t \in(\frac{3}{4},1]
\end{array}
\right.
\]

and $\gamma_{\mathcal C}:[0,1]\rightarrow \mathbb R^2$ by
\[
\gamma_{\mathcal C}(s)=\left(\int_0^s\cos\theta_{\mathcal
C}(t)dt,\int_0^s\sin\theta_{\mathcal C}(t)dt\right).
\]
$\gamma_{\mathcal C}$ is a closed curve parametrized by arclength
and the angle that its tangent vector does with the vector
$(1,0)\in\mathbb R^2$ is $\theta_{\mathcal C}$. Notice that it is
flat outside $\mathcal C$, that is, its curvature is concentrated
on $\mathcal C$. Now we compare the nature of the curvature of
$\gamma_{\mathcal C}$ and the curvature of some well known
examples.

Consider a piecewise regular closed simple planar curve
$\alpha:[a,b]\rightarrow \mathbb{R}^2$. Denote the points where
the curve fails to be smooth by $p_1,\ldots,p_k\in [a,b]$. Suppose
that $\alpha$ is oriented positively and that it is parametrized
by arclength outside $p_1,\ldots,p_k\in [a,b]$. At each $p_i$,
$1\leq i\leq k$, we can associate an angle $K_0(p_i)$ which is the
angle between the vectors $\lim\limits_{t\rightarrow
p_i^-}\alpha^\prime(t)$ and $\lim\limits_{t\rightarrow
p_i^+}\alpha^\prime(t)$. If $K_1$ is the geodesic curvature of
$\alpha$, it is well known that
\[
\sum_{i=1}^k K_0(p_i)+\int_\alpha K_1(t)dt=2\pi,
\]
where the integral is done with respect to the arclength of
$\alpha$.

$K_0$ and $K_1$ can be seen as the zero dimensional and one
dimensional curvature respectively, that is, the curvature that is
concentrated in zero dimensional and one dimensional subsets
respectively. Alternatively, they can be defined as
\[
K_{d}(p)= \lim_{t \rightarrow p^+,\; s \rightarrow
p^-}\frac{\theta(t)-\theta(s)}{(t-s)^{d}},
\]
where $\theta(t)$ is the angle that the tangent vector at
$\alpha(t)$ makes with some fixed vector in $\mathbb R^2$ (If
$p\in \{a,b\}$, then these expressions can be adapted in a obvious
way). The curvature dimension of a curve at a point $p$ can be
defined as the supremum of $d$ such that
\[
\lim_{t \rightarrow p^+,\; s \rightarrow
p^-}\frac{\theta(t)-\theta(s)}{(t-s)^{d}}
\]
is bounded. Observe that this definition has some similarities
with the classical definition of Hausdorff dimension. It is not
difficult to prove that the curvature dimension of
$\gamma_{\mathcal C}$ at $t\in \mathcal C$ is $\left( \ln 2/\ln 3
\right)$.

Now we present a surface which curvature is concentrated in a
non-integer (Hausdorff) dimensional set. Translate the image of
$\gamma_{\mathcal C}$ in such a way that its center of mass is
located at the origin of plane $xy$. Denote such a curve by
$S^1_{\mathcal C}$. Observe that this curve is symmetric with
respect to both axes. Now put the plane $xy$ in the space $xyz$.
Rotate $S^1_{\mathcal C}$ around the $y$-axis. We have a surface
$S^2_{\mathcal C}$ homeomorphic to a bidimensional sphere.

Consider $S^2_{\mathcal C}$ with the induced metric $\widehat g$
of $\mathbb R^3$. It is a $C^1$ metric. Therefore it is
non-regular in the sense of this work. Observe that the points
outside the orbit of $\gamma_{\mathcal C}(\mathcal C)$ are flat,
and that the Hausdorff dimension of the orbit of $\gamma_{\mathcal
C}(\mathcal C)$ is $1+\ln 2/\ln 3$. Hence $(S^2_{\mathcal
C},\widehat g)$ should have its curvature ``concentrated'' in a
$(1+(\ln 2/\ln 3))$ dimensional set.

This sphere indicates that the curvature of a $C^1$ Riemannian
manifold can have a strange behavior. We intend to study this
sphere and other non-regular Riemannian manifolds in future works.

\section*{Acknowledgments}
The author would like to thank Professor Armando Caputi. His
comments and ideas were very valuable at several points of this
work.


\begin{thebibliography}{99}

\bibitem{1} A. D. Alexandrov, {\it Die innere Geometrie der konvexen
Fl\"achen}, Akademie-Verlag, Berlin, 1955.

\bibitem{2} M. Berger, {\it Riemannian Geometry During the Second Half of
the Twentieth Century}, AMS, 2000.

\bibitem{3} G. de Cecco and G. Palmieri, Distanza Intrinseca
su una Variet\`a Riemanniana di Lipschitz, {\it Rend. Sem. Mat.
Univers. Politecn. Torino}, Vol. 46, 2 (1998), 157--170.

\bibitem{4} G. de Cecco and G. Palmieri, Integral distance on
a Lipschitz Riemannian manifold, {\it Math. Z.} 207, (1991),
223--243.

\bibitem{5} J. Cheeger, W. M\"uller and R.
Schrader, On the Curvature of Piecewise Flat Spaces, {\it Comm.
Math. Phys.} 92 (1984), no. 3, 405--454.

\bibitem{6} B. A. Dubrovin, A. T. Fomenko and S. P. Novikov, {\it Modern
Geometry - Methods and Applications, Part. 1}, Graduate Texts in
Mathematics, Vol. 93, Springer-Verlag 1984.

\bibitem{7} L. C. Evans, {\it Partial Differential Equations},
Graduate Studies in Mathematics, Vol. 19, AMS, 1998.

\bibitem{8} G. de Rham, {\it Differentiable Manifolds, Forms, Currents,
Harmonic Forms}, Springer-Verlag, 1984.

\bibitem{9} M. P. do Carmo, {\it Riemannian geometry}.
Translated from the second Portuguese edition by Francis Flaherty.
Mathematics: Theory \& Applications. Birkh\"auser Boston, Inc.,
Boston, MA, 1992.

\bibitem{10} H. Federer, Curvature measures,
{\it Trans. Amer. Math. Soc.} 93 1959 418--491.

\bibitem{11} H. Karcher, Riemannian Center of Mass
and Mollifier Smoothing, {\it Communications on Pure and Applied
Mathematics}, Vol. XXX, 509--541 (1977).

\bibitem{12} S. Kobayashi and K. Nomizu, {\it Foundations of
Differential Geometry}, Vol. 1, Interscience Publishers, 1969.

\bibitem{13} B. B. Mandelbrot, {\it The Fractal Geometry of
Nature}, W. H. Freeman and Co., San Francisco, Calif., 1982.

\bibitem{14} J. Nash, The imbedding problem for
Riemannian manifolds, {\it Ann. of Math.} (2) 63 (1956), 20--63.

\bibitem{15} I. G. Nikolaev, On the parallel displacement
of vectors in spaces with bilaterally bounded curvature in the
sense of A.D. Aleksandrov, {\it Sib. Math. J.} 24 (1983),
106--119; translation from Sib. Mat. Zh. 24:1 (137) (1983),
130--145.

\bibitem{16} A. Petrunin, Parallel transportation
for Alexandrov space with curvature bounded below, {\it Geom.
Funct. Anal.} 8 (1998), no. 1, 123--148.

\bibitem{17} H. L. Royden, {\it Real Analysis}, The
Macmillan Company, 1970.

\bibitem{18} W. Rudin, {\it Real \& complex analysis}, Tata
McGraw-Hill Publishing Co. Limited, New Delhi, 1978.

\bibitem{19} L. Simon, {\it Lectures on Geometric Measure Theory},
Proceedings of the Centre for Mathematical Analysis, Australian
National University, 1983.

\bibitem{20} M. Spivak, {\it A Comprehensive Introduction to
Differential Geometry}, Vol. 1, 2, 3, Second Edition, Publish or
Perish, Inc. 1999.

\end{thebibliography}
\end{document}